\documentclass{amsart} 
\usepackage{amssymb, amscd} 
\usepackage[all]{xy}
\usepackage{epic}
\usepackage{url}

\makeindex

\makeatletter
\def\@doGroups@#1#2#3{#1{#3}#2{#3}}
\def\doGroups#1#2#3#4{%
\def\@doGroups{\@ifnextchar#1\@doGroups@@{\@doGroups@#3#4}}%
\def\@doGroups@@#1##1#2{#3{##1}%
\@ifnextchar#1\@doGroups@@{#4}}%
\@doGroups}

\let\@DefHandle\index

\newcommand\Definition{\doGroups[]\@DefHandle\emph}

\makeatother

\numberwithin{equation}{section}

\newtheorem{theorem}{Theorem}[section]
\newtheorem{lemma}[theorem]{Lemma}
\newtheorem{proposition}[theorem]{Proposition}
\newtheorem{corollary}[theorem]{Corollary}

\newtheorem{definition-lemma}[theorem]{Definition-Lemma}

\theoremstyle{definition}
\newtheorem{definition}[theorem]{Definition}
\newtheorem{formula}[theorem]{Formula}
\newtheorem{example}[theorem]{Example}

\theoremstyle{remark} 
\newtheorem{remark}[theorem]{Remark}

\newenvironment{conventions}
{\trivlist\item[\hskip \labelsep{\itshape Conventions.}] }
{\endtrivlist}

\let\Emph\textit

\def\AA{{\mathbb A}}
\def\CC{{\mathbb C}}  
\def\DD{{\mathbb D}}
\def\EE{{\mathbb E}}    
\def\FF{{\mathbb F}}  
\def\GG{{\mathbb G}} 
\def\HH{{\mathbb H}}

\def\PP{{\mathbb P}}
\def\QQ{{\mathbb Q}} 
\def\RR{{\mathbb R}}

\def\ZZ{{\mathbb Z}} 

\def\k{{\mathbf{k}}}

\def\Db{\mathbf{D}}

\def\Acal{{\mathcal A}}  
\def\Dcal{{\mathcal D}} \def\Ecal{{\mathcal E}} \def\Fcal{{\mathcal F}}
 \def\Hcal{{\mathcal H}} \def\Ical{{\mathcal I}}
 \def\Lcal{{\mathcal L}} \def\Mcal{{\mathcal M}}
 \def\Ocal{{\mathcal O}} 
\def\Rcal{{\mathcal R}} \def\Scal{{\mathcal S}} \def\Tcal{{\mathcal T}}
 \def\Ucal{{\mathcal U}} \def\Vcal{{\mathcal V}}
\def\Xcal{{\mathcal X}}  
 \def\Fl{{\mathcal F\ell}} \def\Uo{{\mathcal U}}
\def\Uc{\overline{\mathcal U}} \def\Vo{{\mathcal V}} \def\Vc{\overline{\mathcal
V}}

\newcommand\CH{\operatorname{CH}}

\newcommand\Flag{\operatorname{Flag}} 
 
\newcommand\Hom{\operatorname{Hom}}
\newcommand\Id{\operatorname{Id}}
\newcommand\Ker{\operatorname{Ker}}
\newcommand\Lie{\operatorname{Lie}} 
\newcommand\Spec{\operatorname{Spec}}
\newcommand\Sp{\operatorname{Sp}} 
\newcommand\SL{\operatorname{SL}} 

\newcommand\area{\operatorname{area}}

\newcommand\codim{\operatorname{codim}}
\newcommand\height{\operatorname{ht}} 
\newcommand\im{\operatorname{Im}} 

\newcommand\rk{\operatorname{rk}}

\makeatletter
\newcommand\symb[1]{\mathop {\operator@font #1}\nolimits}
\makeatother

\newcommand\Dsum{\bigoplus}
\newcommand\dividesnot{\mathop{\!\not|}}

\let\iso\cong
\let\co\colon
\let\l\ell
\newcommand\map[3]{\ensuremath{{#1}\co{#2}\to{#3}}}

\newcommand\set[3][:]{{\ensuremath{\,\{#2 #1 #3\,\}}}}

\newcommand\riso{\mathrel{\hskip2pt\raise-2.5pt\hbox{$\widetilde{\phantom{xx}}$}
\kern-16pt\longrightarrow}}

\newcommand\scprod[2]{\langle #1,#2\rangle}

\newcommand\scpr[2]{\langle#1,#2\rangle}
\newcommand\disjunion{\setbox0\hbox{$\cup$}%
\dimen0\wd0
\setbox1=\hbox to\dimen0{\hfil.\hfil}
\dimen0=\dp1
\mathbin{\vbox{\offinterlineskip\box1\box0}}}

\newcommand\Disjunion{\setbox0\hbox{$\bigcup$}%
\dimen0\wd0
\setbox1=\hbox to\dimen0{\hfil.\hfil}
\dimen0=\dp1
\mathbin{\vbox{\offinterlineskip\box1\box0}}}

\newcommand\Tensor{\bigotimes}
\newcommand\tensor{\otimes}
\newcommand\proofsquare{\nobreak\hfill \hbox{%
\vrule height 5pt 
\kern-.4pt
 \vbox{%
\hrule width 5pt depth0pt height.4pt
 \kern4.6pt \hrule  }
\kern-3.75pt 
\vrule height 5pt}\kern1pt
\par}

\countdef\Yxpos255
\countdef\Yypos254
\countdef\Yxht 253
\countdef\Yxold252
\def\oneYoung#1{\ifnum#1=0
\put(\Yxpos,\Yypos){\line(0,1){0.2}\line(1,0){1}}%
\else
\Yxht=#1
{\ifnum#1<\Yxold
\Yxht=\Yxold
\fi
\put(\Yxpos,\Yypos){\line(0,1){\Yxht}}}%
\advance\Yxht by 1
\multiput(\Yxpos,\Yypos)(0,1){\Yxht}{\line(1,0){1}}%
\fi
\Yxold=#1
}

\def\Young(#1,#2)(#3){%
\def\nxt##1,{\advance\Yxpos by 1
\ifx!##1\else\oneYoung{##1}\expandafter\nxt\fi}%
\Yxold=0
\Yxpos=0
\advance\Yxpos by -1
\Yypos=0
\put(#1,#2){\nxt#3,!,
\put(\Yxpos,\Yypos){\line(0,1){\Yxold}}}}

\begin{document}

\title[Cycle Classes on the Moduli of 
Abelian Varieties]{Cycle Classes of the E-O
Stratification \\
on the Moduli of Abelian Varieties
\footnote{\tt classes.tex\today}} 

\author{Torsten Ekedahl} 

\address{Department of Mathematics\\
Stockholm University\\ 
SE-106 91 Stockholm\\ 
Sweden}
\email{teke@math.su.se} 

\author{Gerard van der Geer}

\address{Faculteit Wiskunde en Informatica, University of
Amsterdam, Plantage Muidergracht 24, 1018 TV Amsterdam, The Netherlands}

\email{geer@science.uva.nl}


\begin{abstract}
We introduce a stratification on the space of symplectic flags on the de Rham
bundle of the universal principally polarized abelian variety in positive
characteristic and study its geometric properties like irreducibility of the
strata and we calculate the cycle classes. When the characteristic $p$ is
treated as a formal variable these classes can be seen as a deformation of the
classes of the Schubert varieties for the corresponding classical flag variety
(the classical case is recovered by putting $p$ equal to $0$). We relate our
stratification with the E-O stratification on the moduli space of principally
polarized abelian varieties of a fixed dimension and derive properties of the
latter. Our results are strongly linked with the combinatorics of the Weyl group
of the symplectic group.
\end{abstract}
\maketitle

\begin{section}{Introduction}

The moduli space ${\Acal}_g$ of principally polarized abelian varieties of
dimension $g$ is defined over the integers. For the characteristic zero fibre
${\Acal}_g \otimes \CC$ we have an explicit description as an orbifold
$\Sp_{2g}(\ZZ) \backslash {\Hcal}_g$ with ${\Hcal}_g$ the Siegel upper half
space of degree $g$.  It is a recent insight, though, that perhaps the positive
characteristic fibres ${\Acal}_g \otimes \FF_p$ are more accessible than the
characteristic zero one.  A good illustration of this is provided by the E-O
stratification of ${\Acal}_g \otimes \FF_p$, a stratification consisting of
$2^g$ strata, each of which is quasi-affine. It was originally defined by
Ekedahl and Oort (see \cite{O1}) by analyzing the structure of the kernel of
multiplication by $p$ of an abelian variety.  It turns out that this group
scheme can assume $2^g$ forms only, and this led to the strata. For $g=1$ the
two strata are the locus of ordinary and that of supersingular elliptic
curves. Some strata possess intriguing properties. For example, the stratum of
abelian varieties of $p$-rank $0$ is a complete subvariety of ${\Acal}_g \otimes
\FF_p$ of codimension $g$, the smallest codimension possible. No analogue in
characteristic $0$ of this stratum nor of the stratification is known, and in
fact, Keel and Sadun (\cite{K-S}) proved that complete subvarieties of
${\Acal}_g \otimes \CC$ of codimension $g$ do not exist for $g\geq 3$.

While trying to find cycle classes for the E-O strata we realised that the
strata could be described as degeneration loci for maps between vector bundles
and as such loci are indexed by Young diagrams our attention was turned towards
the combinatorics of the Weyl group. When considered in this light it is clear
that much of the combinatorics of \cite{O1} is closely related to the Weyl group
$W_g$ of $\Sp_{2g}$, which is the group relevant for the analytic description of
$\Acal_g\otimes \CC$. The main idea of this paper is to try to make this
connection more explicit. More precisely, the combinatorics of the E-O strata is
most closely related to the combinatorics associated to $W_g$ and the sub-Weyl
group corresponding to the maximal parabolic subgroup $P$ of elements of
$\Sp_{2g}$ fixing a maximal isotropic subspace in the $2g$-dimensional
symplectic vector space. Indeed, this sub-Weyl group is $S_g$, the group of
permutations on $g$ letters (embedded as a sub-Weyl group in $W_g$), and the E-O
strata are in bijection with the cosets $W_g/S_g$; we shall use the notation
$\Vo_\nu$ for the (open) stratum of ${\Acal}_g \otimes \FF_p$ corresponding to
$\nu \in W_g/S_g$ (and $\Vc_\nu$ for its closure). The coset space $W_g/S_g$ is
also in bijection with the set of Bruhat cells in the space of maximal totally
isotropic flags $\Sp_{2g}/P$ and we believe this to be no accident. (The formal
relation between $\Acal_g$ and $\Sp_{2g}/P$ is that $\Sp_{2g}/P$ is the
\emph{compact dual} of $\Hcal_g$.)

In order to push the analogy further we introduce a ``flag space'' $\Fcal_g \to
\Acal_g$ whose fibres are isomorphic to the fibres of the map $\Sp_{2g}/B \to
\Sp_{2g}/P$, where $B$ is a Borel subgroup of $P$. In positive characteristic we
define (and this definition clearly makes sense only in positive characteristic)
a stratification of $\Fcal_g$ whose open strata $\Uo_w$ and corresponding
closures $\Uc_w$ are parametrised by the elements of $W_g$. This is very similar
to the Bruhat cells of $\Sp_{2g}/B$ and their closures, the Schubert strata,
which are also parametrised by the elements of $W_g$. Our first main result is
that this is more than a similarity when one is working locally; we show that
for each point of $\Fcal_g$ there is a stratum preserving local isomorphism (in
the \'etale topology) taking the point to some point of $\Sp_{2g}/B$. As a lot
is known about the local structure of the Schubert varieties we immediately get
a lot of information about the local structure of our strata. The first
consequence is that the dimension of $\Uc_w$ is everywhere equal to the length
of $w$. A very important consequence is that the $\Uc_w$ are all normal; this
situation
differs markedly from the case of the closed E-O strata which in general are not
normal. Another consequence is that the inclusion relation between the strata is
given exactly by the Bruhat-Chevalley order on $W_g$. (A much more sophisticated
consequence is that the local structure of the $\l$-adic intersection complex
for a closed stratum is the same as for the Schubert varieties and in particular
that the dimensions of its fibres over the open strata of the closed stratum are
given by the Kazhdan-Lusztig polynomials. We shall not however pursue that in
this article.)

We make several applications of our results on the structure of the strata
$\Uc_w$. The first, and most important, is that by construction the strata
$\Uo_w$ are defined as the locus where two symplectic flags on the same vector
bundle are in relative position given by $w$. As we have shown that they have
the expected codimension and are reduced we can use formulas of Fulton, and
those of Pragacz and Ratajski as crystallized in the formulas of Kresch and
Tamvakis to get formulas for the cycle classes of the strata. A result of Fulton
give such formulas for all strata but in terms of a recursion formula that we
have not been able to turn into a closed formula; even though these formulas
should have independent interest we can use them to get formulas for the E-O
strata as follows. If $w \in W_g$ is minimal for the Bruhat-Chevalley order in
its coset $wS_g$, then $\Uo_w$ maps by a finite \'etale map to the open E-O
stratum $\Vo_\nu$ corresponding to the coset $\nu:=wS_g$. We can compute the
degree of this map in terms of the combinatorics of the element $w$ and we then
can push down our formula for $\Uc_w$ to obtain formulas for the cycle classes
of the E-O strata. Also the formulas of Kresch and Tamvakis can be used to give
the classes of E-O strata. One interesting general consequence is that each
class is a polynomial in the Chern classes $\lambda_i$ of the Hodge bundle whose
coefficients are polynomials in $p$. This is a phenomenon already visible in the
special cases of our formula that were known previously; the oldest such example
being Deuring's mass formula for the number of supersingular elliptic curves
(weighted by one over the cardinalities of their automorphism groups) that says
that this mass is $(p-1)/12$. This appears in our context as the combination of
the formula $(p-1)\lambda_1$ for the class of the supersingular locus and the
formula $\deg\lambda_1=1/12$. We interpret these results as giving rise to
elements in the \Emph{$p$-tautological ring}; this is the ring obtained from the
usual \Emph{tautological ring}, the ring generated by the Chern classes of the
Hodge bundle, by extending the scalars to $\ZZ\{p\}$, the localisation of the
polynomial ring $\ZZ[p]$ at the polynomials with constant coefficient $1$. Hence
we get elements parametrised by $W_g/S_g$ in the $p$-tautological ring and we
show that they form a basis for the $p$-tautological ring. Putting $p$ equal to
$0$ maps these elements to elements of the ordinary tautological ring which can
identified with the Chow ring of $\Sp_{2g}/P$ and these elements are the usual
classes of the Schubert varieties. It seems that these results call for a
$p$-Schubert calculus in the sense of a better understanding of these elements
of the $p$-tautological ring and for instance their behaviour under
multiplication.

However, there seems be a more intriguing problem. We have for each $w \in W_g$
a stratum in our flag space and they push down to elements of the
$p$-tautological ring under the projection map to $\tilde\Acal_g$ (a toroidal
compactification of $\Acal_g$). When setting $p$ to $0$ these elements
specialise to the classes of the images of the Schubert varieties of
$\Sp_{2g}/B_g$ in $\Sp_{2g}/P$ and for them the situation is very simple, either
$w$ is minimal in its $S_g$ coset and then the Schubert variety maps
birationally to the corresponding Schubert variety of $\Sp_{2g}/P$ or it is not
and then it maps to $0$. However, when it comes to the elements of the
$p$-tautological ring this only allows us to conclude -- in the non-minimal case
-- that the coefficients are divisible by $p$ and indeed in general they are not
zero. We show that unless they map to $0$ they will always map to a multiple of
a class of an E-O stratum. When the element is minimal in its $S_g$ coset this
stratum is indexed by the coset spanned by the element, but our considerations
give an extension of this map from elements minimal in their cosets to a larger
class of elements. We give some examples of this map but in general it seems a
very mysterious construction.

Another application is to the irreducibility of our strata (and hence also to
the strata of the E-O stratification as they are images of some of our
strata). As the strata are normal this is equivalent to the connectedness of a
stratum and this connectedness can sometimes be proved by an arithmetic
argument. It is natural to ask if this method produces all the irreducible
strata and for the characteristic large enough (the size depending on $g$) we
can show that indeed it does. This is done using a Pieri type formula for our
strata obtained by applying a result of Pittie and Ram. A Pieri type formula for
multiplying the class of a connected cycle by an ample line bundle has as a
consequence that a part of boundary is supported by an ample line bundle and
hence that it is connected. Exploiting that together with the fact that
$\lambda_1$ is an ample line bundle on $\Acal_g$ allows to show the desired
converse. We are forced to assume that the characteristic is large (and are
unable to specify how large) as we need to know how high a power of $\lambda_1$
one needs to twist the exterior powers of the dual of the Hodge bundle to make
it generated by global sections.

There is a particular element of $w_\emptyset \in W_g$ which is the largest of
the elements that are minimal in their right $S_g$-cosets which has the property
that $\Uc_{w_\emptyset}$ maps birationally onto $\Acal_g$. It is really the
strata that are contained in this stratum that seem geometrically related to
$\Acal_g$ and indeed the elements $w \in W_g$ lying below $w_\emptyset$ are the
ones of most interest to us. 
(The rest of $\Fcal_g$ appears mostly as a technical
device for relating our strata to the Schubert varieties.) It should be of
particular interest to understand the map $\Uc_{w_\emptyset} \subset \Fcal_g \to
\Acal_g$. It follows from a result of Oort on Dieudonn\'e modules that the
inverse image of an open E-O stratum under this map is a locally constant
fibration which focuses particular interest on the fibres of the map, this fibre
depending only on the element of $\nu \in W_g/S_g$ that specifies the E-O
stratum. We call these fibres \emph{punctual flag spaces}, using the notation
$\Fcal_\nu$ for the one associated to $\nu$, and establish some basic properties
of them: We determine their connected components, showing in particular that two
points in the same connected component can be connected by a sequence of quite
simple rational curves. We also show that knowing which strata $\Uo_w$ have
non-empty intersections with a given punctual flag space would determine the
inclusion relations between the E-O strata. We intend to make a more detailed
study of the punctual flag spaces and show how that knowledge will give more
information on these inclusion relations in a sequel to this paper.

The geometric points of the stratum $\Uc_{w_\emptyset}$ correspond to symplectic
flags of subgroup schemes of the kernel of multiplication by $p$ on the
principally polarized abelian variety that is given by the image of the point in
$\Acal_g$. This is of course is reminiscent of de Jong's moduli stack
$\Scal(g,p)$ of $\Gamma_0(p)$-structures. The major difference (apart from the
fact that $\Uc_{w_\emptyset}$ only makes sense in positive characteristic) is
that the $g$-dimensional element of the flag is determined by the abelian
variety in our case. We shall indeed identify $\Uc_{w_\emptyset}$ with the
component of $\Scal(g,p)$ that is the closure of the ordinary abelian varieties
provided with a flag on the local part of the kernel of multiplication by
$p$. As a consequence we get that that component of $\Scal(g,p)$ is normal and
Cohen-Macaulay.

This paper is clearly heavily inspired by \cite{O1}. The attentive reader will
notice that we reprove some of the results of that paper, sometimes with proofs
that are very close to the proofs used by Oort. We justify such duplications by
our desire to emphasize the relations with the combinatorics of $W_g$ and the
flag spaces. Hence, we start with (a rather long) combinatorial section where
the combinatorial aspects have been separated from the geometric ones. We hope
that this way of presenting the material will be as clarifying to the reader as
it has been to us. We intend to continue to exploit the relations with the flag
spaces in a future paper that will deal with K3-surfaces and one on Dieudonn\'e
modules.  Since its announcement in \cite{vdG1} our idea of connecting the E-O
stratification on ${\Acal}_g$ with the Weyl group and filtrations on the de Rham
cohomology has been taken up in other work. In this connection we want to draw
attention to papers by Moonen and Wedhorn, cf.,
\cite{moonen01::group+weyl,moonen04::discr}.

We would like to thank Piotr Pragacz for some useful  comments.

\begin{conventions}
We shall exclusively work in positive characteristic $p>0$. After having
identified final types and final elements in Section \ref{sec:combinatorics} we
shall often use the same notation for the final type (which is a function on
$\{1,\dots,2g\}$) and the corresponding final element (which is an element of
the Weyl group $W_g$). In Sections \ref{sec:Pieri formulas} and
\ref{sec:Irreducibility} our strata will be considered in flag spaces over not
just $\Acal_g$ and $\tilde\Acal_{g}$ but also over the corresponding moduli
stacks with a level structure.
\end{conventions}
\end{section}
\tableofcontents
\begin{section}{Combinatorics}
\label{sec:combinatorics}

This section is of a preparatory nature and deals with the combinatorial aspects
of the E-O stratification. The combinatorics is determined by the Weyl group of
the symplectic group of degree $g$. A general reference for the combinatorics of
Weyl groups is \cite{B-L}. We start by recalling some general notations and
facts about $W_g$ and its Bruhat-Chevalley order. We then go on to give various
descriptions of the minimal elements in the $S_g$ cosets (which we presume are
well-known). The short subsection on shuffles will be used to understand the
r\^ole that the multiplicative and \'etale part of the Barsotti-Tate group
play in our stratification in the case of positive $p$-rank.
\begin{subsection}{Final Elements in The Weyl Group}
The Weyl group $W_g$ of type $C_g$ in Cartan's terminology is isomorphic to the
semi-direct product $S_g \ltimes (\ZZ / 2\ZZ)^g$, where the symmetric group
$S_g$ on $g$ letters acts on $(\ZZ / 2\ZZ)^g$ by permuting the $g$
factors. Another description of this group, and the one we shall use here, is as
the subgroup of the symmetric group $S_{2g}$ of elements which map any symmetric
$2$-element subset of $\{ 1, \ldots, 2g\}$ of the form $\{ i , 2g+1-i \}$ to a
subset of the same type:
$$
W_g = \{ \sigma \in S_{2g} \colon \sigma(i) + \sigma(2g+1-i) = 2g+1 { \hbox {
\rm for } } i=1,\ldots, g \}.
$$
The function $i \mapsto 2g+1-i$ on the set $\{ 1, \dots, 2g\}$ will occur
frequently. We shall sometimes use the notation $\overline{\imath}$ for
$2g+1-i$. Using it we can say that $\sigma \in S_{2g}$ is an element of
$W_g$ precisely when
$\sigma(\overline{\imath})=\overline{\sigma(i)}$ for all $i$. This makes the
connection with another standard description of $W_g$, namely as a group
of \Emph{signed
permutations}. An element in this Weyl group has a \Definition{length} and a
\Definition{codimension} defined by:
$$
\ell(w) = \#\{i < j \le g : w(i) > w(j) \} + \# \{ i \leq j \le g : w(i)+w(j)>
2g+1\}
$$
and
$$
\mathrm{codim}(w) =\#\{i < j \le g \co  w(i) < w(j) \} + \# \{ i \leq j \le g\co 
w(i)+w(j) < 2g+1\}
$$
and these satisfy the equality
$$
\ell(w)+ {\rm codim}(w) = g^2.
$$
We shall use the following notation for elements in $W_g$. By $[a_1,a_2,\ldots,
a_{2g}]$ we mean the permutation of $\{1,2,\ldots,2g\}$ with $\sigma(i)=a_i$.
Since $\sigma(i)$ determines $\sigma(2g+1-i)$ for $1 \leq i \leq g$ sometimes we
use the notation $[a_1,\ldots,a_g]$ instead (when $g \le 4$ we shall even
dispence with the commas and write $[a_1\ldots a_g]$ which should cause no
confusion as all the $a_i$ must be single digit numbers. We shall also use
cycle notation for permutations. In particular, for $1\leq i < g$ we let $s_i\in
S_{2g}$ be the permutation $(i, i+1)(2g-i, 2g+1-i)$ in $W_g$ which interchanges
$i$ and $i+1$ and we let $s_g=(g, g+1) \in S_{2g}$.  Then
$(W=W_g,S=\{s_1,\ldots,s_g\})$ is a Coxeter system.

Let $(W,S)$ be a Coxeter system and $a \in W$. If $X$ is a subset of $S$ we
denote by $W_X$ the subgroup of $W$ generated by $X$.  It is well-known that for
any subset $X$ of $S$ there exists precisely one element $w$ of minimal length
in $aW_X$ and it has the property that every element $w' \in aW_X$ can be
written in the form $w'= wx$ with $x \in W_X$ and
$\ell(w')=\ell(w)+\ell(x)$. Such an element $w$ is called a
\Definition{$X$-reduced element}. (Cf.\ \cite{Bou}, Ch.\ IV, Excercises \S 1, No
3).)

Let $W=W_g$ be the Weyl group and $S$ be the set of simple reflections. If we
take $X=S\backslash \{ s_g\}$ then we find
$$
W_X=\{ \sigma \in W_g\co  \sigma \{1,2,\ldots,g\} = \{1,2,\ldots,g\}\}\cong S_g.
$$

There is a natural partial order on $W_g$ with respect $W_X$, the Bruhat-Chevalley order.
It is defined in terms of Schubert cells $X(w_i)$ by:
$$w_1 \geq w_2 \iff
X(w_1) \supseteq X(w_2).$$ Equivalently, if we define
\begin{equation}\label{rw definition}
r_w(i,j) :=\# \{ a \leq i \co  \, w(a) \leq j \}
\end{equation}
then we have the combinatorial characterization
$$
w_1 \leq w_2 \iff r_{w_1}(i,j) \geq r_{w_2}(i,j) \quad \hbox{ \rm for all $1\leq
i,j \leq 2g$}.
$$
(Indeed it is easy to see that it is enough to check this for all $1\le i \le g$
and $1\le j \le 2g$.)  Chevalley has shown that $w_1 \geq w_2 $ if and only if
any (hence every) $X$-reduced expression for $w_1$ contains a subexpression
(obtained by just deleting elements) that is a reduced expression for $w_2$;
here reduced means that $w_2$ is written as a product of $\ell(w_2)$ elements
from $S$.  Again, a reference for these facts is \cite{B-L}.

We now restrict to the following case. Let $V$ be a symplectic vector space over
$\QQ$ and consider the associated algebraic group $G=\Sp(V)$.  If $E \subset V$
is a maximal isotropic subspace then the stabilizer of the flag $(0) \subset E
\subset V$ is a parabolic subgroup conjugate to the standard parabolic
corresponding to $H:=S\backslash \{ s_g\} \subset S$. Hence $W_H$ will denote
the subgroup of $W_g$ generated by the elements of $H$ and we will also use the
notation $P_H$ for the parabolic subgroup corresponding to $W_H$, i.e., the
subgroup of the symplectic group stabilizing a maximal totally isotropic
subgroup. As $W_H$ consists of the permutations of $W_g$ that stabilise the
subsets $\{1,\dots,g\}$ and $\{g+1,\dots,2g\}$ and the restriction of its action
to $\{1,\dots,g\}$ determine the full permutation we may identify $W_H$ with
$S_g$, the group of permutations of $\{1,\dots,g\}$ and we shall do so without
further mention. (This is of course compatible with the fact that $H$ spans an
$A_{g-1}$-subdiagram of the Dynkin diagram of $G$.) There are $2^g=|W_g|/|W_H|$
elements in $W_g$ which are $H$-reduced elements. These $2^g$ elements will be
called \Definition{final} elements of $W_g$. The Bruhat-Chevalley order between elements
in $W_g$ as well as the condition for being $H$-reduced can be conveniently
expressed in terms of the concrete representation of elements of $W_g$ as
permutations in the following way.

Let $A, B$ two finite subsets of $\{1,2,\ldots,g\}$ of the same cardinality. We
shall write $A \prec B$ if the $i$'th largest element of $A$ is $\le$ the $i$'th
largest element of $B$ for all $1 \le i \le |A|$.
\begin{lemma}
\label{bruhat}
i) If $w = [a_1a_2\dots a_g]$ and $w' = [b_1b_2\dots b_g]$ are two elements of
$W_g$, then $w \le w'$ in the Bruhat-Chevalley order precisely when for all $1 \le d \le
g$ we have $\{a_1,a_2,\dots,a_d\} \prec \{b_1,b_2,\dots,b_d\}$.

ii) Let $w = [a_1a_2\dots a_g]$ and $w^f$ be the final element of $wW_H$ and let
$w' = [b_1b_2\dots b_g]$ then $w^f \le w'$ in the Bruhat-Chevalley order precisely when
$\{a_1,a_2,\dots,a_g\} \prec \{b_1,b_2,\dots,b_g\}$.

iii) An element $\sigma \in W_g$ is $H$-reduced (or final) if and only if
$\sigma(i)<\sigma(j)$ for all $1\leq i < j \leq g$.  Also, $\sigma$ is
$H$-reduced if and only if $\sigma$ sends the first $g-1$ simple roots into
positive roots.
\begin{proof}
See for instance \cite[p.\ 30]{B-L}.
\end{proof}
\end{lemma}

\end{subsection}
\par
\begin{subsection}{Final Types and Young Diagrams}
\noindent
There are other descriptions of final elements that are sometimes equally
useful. These involve maps of $\{1,2,\ldots,2g\}$ to $\{1,2,\ldots,g\}$ and
certain Young diagrams. We begin with the maps.
\begin{definition}
A \Definition{final type} (of degree $g$) is an increasing surjective map
$$
\nu\co  \{0,1,2,\ldots,2g\} \to \{ 0,1,2,\ldots,g\}
$$
satisfying
$$
\nu(2g-i) = \nu(i)-i+g \qquad \hbox{\rm for} \quad 0\leq i \leq g.
$$
\end{definition}
Note that we either have $\nu(i+1)=\nu(i)$ and then $\nu(2g-i)=\nu(2g-i-1)+1$ or
$\nu(i+1)=\nu(i)+1$ and then $\nu(2g-i)=\nu(2g-i-1)$. A final type is determined
by its values on $\{0,1,\ldots,g\}$. There are $2^g$ final types of degree $g$
corresponding to the vectors $(\nu(i+1)-\nu(i))_{i=0}^{g-1} \in \{0,1\}^g$.
The notion of a final type was introduced by Oort (cf.\ \cite{O1}).

To an element $w \in W_g$ we can associate the final type $\nu_w$ defined by
$$
\nu_w(i)=i-r_w(g,i).
$$
This is a final type because of the rule $r_w(g,2g-i)-r_w(g,i)=g-i$ 
which follows by induction on $i$ 
from the fact that $w(2g+1-a)=2g+1-w(a)$. It depends only on the coset
$wW_H$ of $w$ since a permutation of the $a\leq g$ does not change
the definition of $r_w(g,i)=\# \{a \leq g \colon w(a)\leq i \}$.

Conversely, to a final type $\nu $ we now
associate the following element $w_{\nu}$ of the Weyl
group, a permutation of $\{ 1,2, \ldots, 2g\}$ as follows. Let
$$
\beta = \{ i_1,i_2,\dots,i_k \}= \{ 1 \leq i \leq g \co  \nu(i)=\nu(i-1)\}
$$
with $i_1 < i_2 <\ldots $ given in \Emph{increasing} order and let
$$
\beta^c= \{ j_1, j_2, \dots,j_{g-k} \}
$$
be the elements of $\{ 1,2, \ldots, g\}$ not in $\xi$, in \Emph{decreasing}
order.  We then define a permutation $w_\nu$ by mapping $1 \le s \le k$ to
$i_s$ and $k+1 \le s \le g$ to $2g+1-j_{s-k}$. The requirement that $w_\nu$
belong to $W_g$ now completely specifies $w_\nu$ and by construction
$w_{\nu}(i)<w_{\nu}(j)$ if $1\leq i<j \leq g$ so it is a final element of
$W_g$. It is clear from \ref{bruhat} that we get in this way all final elements
of $W_g$.  
The Bruhat-Chevalley order for final elements can also be read off from the final type
$\nu$. We have $w \geq w'$ if and only if $\nu_w \geq \nu_{w'}$.
This follows from \ref{bruhat}, ii).

We summarize:

\begin{lemma} By associating to a final type $\nu$ the element $w_{\nu}$
and to a final element $w \in W_g$ the final type $\nu_w$
we get an order preserving 
bijection between the set of $2^g$ final types and the set of final
elements of $W_g$.
\end{lemma}

The final types are in bijection with certain Young diagrams: Our Young diagrams
will be put in a position that is opposite to the usual positioning, i.e.,
larger rows will be below smaller ones and the rows will be lined up to the
right. Furthermore, we shall make Young diagrams correspond to partitions by
associating to one the parts that are the lengths of the rows of the diagram. We
shall say that a Young diagram is \Definition{final of degree $g$} if its parts
are $\le g$ and no two parts are equal. They therefore correspond to subsets
$\xi=\{g\geq \xi_1 > \xi_2 > \ldots > \xi_r \}$ of $\{1,2,\dots,g\}$.

To a final type $\nu$ we now associate the Young diagram $Y_{\nu}$ whose
associated subset $\xi$ is defined by
$$
\xi_j = \# \{ i \co  1 \leq i \leq g , \, \nu(i)\leq i-j \}.
$$
A pictorial way of describing the Young diagram
is by putting a stack of $i-\nu(i)$ squares in 
vertical position $i$ for $1\leq i \leq g$.
\begin{example}
This example corresponds to 
$$
\{ \nu(i) \co  i=1,...,g\} = \{1,2,\dots,g-5, g-5, g-4, g-4, g-3,
g-3 \}
$$
and hence $\xi= \{5 , 3, 1\}$.

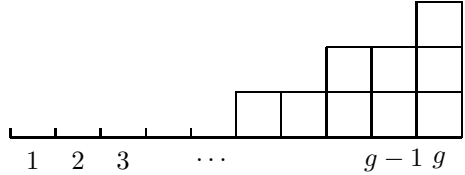
\begin{figure}[htb]
\begin{picture}(10,4)
\Young(0,1)(0,0,0,0,0,1,1,2,2,3) \put(0,0){\makebox(1,1){1}}
\put(1,0){\makebox(1,1){2}} \put(2,0){\makebox(1,1){3}}
\put(3,0){\makebox(1,1){}} \put(4,0){\makebox(1,1){$\dots$}}
\put(8,0){\makebox(1,1){$g-1$}} \put(9,0){\makebox(1,1){$g$}}
\end{picture}
\end{figure}
\par
\end{example}
The final elements $w$ in $W_g$ are in $1-1$-correspondence with the 
elements of $W_g/W_H$. The group $W_g$ acts on $W_g/W_H$ by multiplication
on the left, i.e., by the permutation representation. Therefore $W_g$
also acts on the set of final types and the set of final Young diagrams.
To describe this action we need the notion of a break point.

By a \Definition{break point} of a final type $\nu$ we mean an integer $i$ with
$1 \leq i \leq g$ such that either
\begin{enumerate}
\item $\nu(i-1)=\nu(i)\neq \nu(i+1)$, or
\item $\nu(i-1)\neq \nu(i)=\nu(i+1)$.
\end{enumerate}
If $i$ is not a break point of $\nu=\nu_w$ then $\nu(i+1)=\nu(i-1)$ or
$\nu(i+1)=\nu(i-1)+2$ and then $\nu_{s_iw}=\nu_w$. In particular,
$g$ is always a break point. The set of break points of $\nu$
equals
$$
\{ 1\leq i \leq g \colon  \nu_{s_iw}\neq \nu_{w} \}.
$$
Since $\nu=\nu_w$ determines a coset $wW_H$ we have that $i$ is not a 
break point of $\nu$ if and only if $w^{-1}s_iw\in W_H$, i.e., if and only if
$wW_H$ is a fixed point of $s_i$ acting on $W_g/W_H$.
The action of $s_i$ on a final type $\nu$ is as follows: if $i$ is not a break
point then $\nu$ is fixed, otherwise replace the value of $\nu$ at $i$ by
$\nu(i)+1$ if $\nu(i-1)=\nu(i)$ and $\nu(i)-1$ otherwise.

If $w$ is a final element given by the permutation $[a_1,a_2,\ldots,a_g]$
then it defines a second final element, called the \Definition{complementary
permutation}, defined by the permutation $[b_1,b_2,\ldots,b_g]$,
where $b_1 < b_2 < \ldots < b_g$ are the elements of the complement
$\{1,2,\ldots,2g\} 
\backslash \{a_1,\ldots,a_g\}$. If $\xi$ is the partition defining
the Young diagram of $w$ then $\xi^c$ defines the Young diagram of the
complementary permutation. The set of break points of $w$ and its
complementary element are the same.

\begin{lemma}\label{complementary}
Let $w$ be a final element with associated final type $\nu$ and complementary
element $v$.

i) We have that $v=\sigma_1w\sigma_0=w\sigma_1\sigma_0$, where $\sigma_0$
(resp.\ $\sigma_1$) is the element of $S_g$ (resp.\ $W_g$) that maps $1 \le i
\le g$ to $g+1-i$ resp.\ to $2g+1-i$.

ii) If $1 \leq i \leq g$ has $\nu(i-1)\neq \nu(i)$ then $v^{-1}(i)=\nu(i)$ and
if it has $\nu(i-1)= \nu(i)$ then $v^{-1}(i)=2g+1-\nu(2g+1-i)$.
\begin{proof}
As $w$ maps $i$ to  $a_i$ we have that $2g+1-a_i$ is not among the $a_j$ and
hence $b_i=2g+1-a_{g+1-i}$ (using that both the $a_i$ and $b_i$ are increasing
sequences). This gives
$\sigma_1w\sigma_0(i)=\sigma_1(a_{g+1-i})=2g+1-a_{g+1-i}=b_i$ but we note that
as $w \in W_g$ it commutes with $\sigma_1$.

If $\nu(i-1)\neq \nu(i)$ and, say, $\nu(i)=i-k$ then we have $k$ natural numbers
$1\leq i_1 < i_2 < \cdots < i_k< i$ such that $\nu(i_j-1)=\nu(i_j)$. By the
definition of $v$ we then have $v(i-k)= (i-k)+k$ since the $k$ values $i_j$
($j=1,\ldots,k$) are values for $v$, hence not of $w$. The second part is
similar.
\end{proof}
\end{lemma}
\begin{remark}
The permutations $\sigma_1$ and $\sigma_0$ of course have clear root theoretic
relevance, they are respectively the longest element of $W_g$ and
$S_g$. Multiplication by them reverses the Bruhat-Chevalley order. Similarly it is clear
that going from a final element to its complementary element also reverses the
Bruhat-Chevalley order among the final elements and the first part of our statement says
that that operation is obtained by multiplying by $\sigma_1$ and $\sigma_0$, a
fact that we presume generalises. Somewhat curiously our use of the
complementary permutation seems unrelated to these facts.
\end{remark}

In terms of Young diagrams the description is analogous and gives us a way to
write the element $w_{\nu}$ as a reduced product of simple reflections. To each
$s_i$ we can associate an operator on final Young diagrams.  If $Y$ is a final
diagram, $s_i$ is defined on $Y$ by adding or deleting a box in the $i$'th
column if this gives a final diagram (only one of the two can give a final
diagram) and then $s_iY$ will be that new diagram; if neither adding nor
deleting such a box gives a final Youg diagram we do nothing. In terms of the
description as subsets $\xi$ adding a box corresponds to $g+1-i \in \xi$ and
$g+2-i \notin \xi$ and then $s_i\xi = (\xi\setminus \{g+1-i\})\cup\{g+2-i\}$.
It is then clear that for any final Young diagram $Y$ there is a word
$s_{i_1}s_{i_2}s_{i_3}\cdots s_{i_k}$ such that $Y=s_{i_1}s_{i_2}s_{i_3}\cdots
s_{i_k}\emptyset$, where $\emptyset$ denotes the empty Young diagram. Comparison
with the action of $s_i$ on final types and the correspondence between final
types and Young diagrams shows that the action of $s_i$ on diagrams is indeed
obtained from that on final types. If we now have a word $t=s_{i_1}\dots
s_{i_k}$ in the $s_i$ we can make it act on Young diagram by letting each
individual $s_i$ act as specified. Note that this action only depends on the
image of $t$ in $W_g$ but for the moment we want to consider the action by
words. We define the \Definition{area} of a Young diagram $Y$ to be the number
of boxes it contains. We shall say that the word $t$ is \Definition{building} if
the area of $t\emptyset$ is equal to $k$, the length of the \Emph{word} (not the
resulting element). This is equivalent to the action of $s_{i_r}$ adding a
box to $s_{i_{r+1}}\dots s_{i_k}\emptyset$ for all $r$.
\begin{lemma}\label{final factoring}
i) If $\nu$ is a final type and $t$ is a word in the $s_i$ such that
$Y_\nu^c=t\emptyset$ then $w_{\nu}=w$, where $w$ is the image of $t$ in
$W_g$ and $\ell(w_{\nu})=g(g+1)/2-\textrm{area}(Y_{\nu})$.

ii) $t$ is $H$-reduced if and only if $t$ is building.
\begin{proof} 
To prove i) we begin by noting that $t\emptyset$ only depends on the image of
$t$ in $W_g$ so that i) is independent of the choice of $t$. Hence we may prove
it by choosing a particular $t$ using induction on the area of $Y_{\nu}^c$. Note
that $g(g+1)/2-\area(Y_{\nu})=\area(Y_\nu^c)$ so that the last part of i) says
that $\ell(w_\nu)=\area(Y_\nu^c)$.  The final type $\nu$ with $\nu(i)=0$ for
$i\leq g$ corresponds to final diagram $Y_{\nu}$ with empty complementary
diagram.  We have $w_{\nu}= 1 \in W_g$, the empty product and it has length
$0$. This proves the base case of the induction.  Suppose we proved the
statement for diagram $Y_{\nu}$ with $\area(Y_{\nu}^c) \le a$.  Adding one block
to $Y_{\nu}^c$ to obtain $Y_{\nu '}^c$ means that for some $i$ we have $g+1-i
\in \xi^c$ and $g+2-i \notin \xi^c$, where $\xi^c$ is the subset corresponding
to $Y_\nu^c$, the new subset is $(\xi ')^c = (\xi\setminus
\{g-i\})\cup\{g-i+1\}$. This means that if $i<g$ there are $b < a \le g$ such
that $w_\nu(b)=i$ and $w_\nu(a)=2g-i$ and $w_{\nu '}(b)=i+1$ and $w_{\nu
'}(a)=2g+1-i$ and the rest of the integers between $1$ and $g$ mapped to the
same elements. (The case when $i=g$ is similar and left to the reader.) This
makes it clear that we have $w_{\nu '}=s_iw_\nu$ so by the induction $t$ maps to
$w_\nu$. It remains to show the formula for $\ell(w_\nu)$. In the definition of
$\ell(w)$ only the second term contributes as $w_\nu(i) < w_\nu(j)$ if $i<j\le
g$. Now, the only difference in the collections of sums $w(i)+w(j)$ for $i\le j$
and $w$ equal to $w_\nu$ and $w_{\nu '}$ appears for $(i,j)=(b,a)$ and we have
$w_\nu(b)+w_\nu(a)=2g$ and $w_{\nu '}(b)+w_{\nu '}(a)=2g+2$ so that the length
of $w_{\nu '}$ is indeed one larger than that of $w_{\nu '}$.

As for ii), we have that $t\emptyset=Y_\nu^c$, where $\nu$ is the final type of
$w$ and then ii) is equivalent to $t$ being $H$-reduced if and only
$\area(Y_\nu^c)$ is equal to the length of $t$. However by i) we know that
$\area(Y_\nu^c)$ is equal to $\l(w_\nu)$ and $t$ is indeed $H$-reduced precisely
when its length is equal to $\l(w_\nu)$.
\end{proof}
\end{lemma}
\begin{example}
Consider again the Young diagram of the previous example but now for $g=5$:
\begin{figure}[htb]
\begin{picture}(5,4)
\Young(0,1)(1,1,2,2,3)
\put(0,0){\makebox(1,1){1}}
\put(1,0){\makebox(1,1){2}}
\put(2,0){\makebox(1,1){3}}
\put(3,0){\makebox(1,1){4}}
\put(4,0){\makebox(1,1){5}}
\end{picture}
\end{figure}
We have $\xi=\{5,3,1\}$ and thus $\xi^c=\{2,4\}$, so
$w_\nu=[13579]$ and $w_\nu$ can be written as $s_4s_5s_2s_3s_4s_5$ (notice that
permutations act from the left on diagrams).
\end{example}

\bigskip
We now characterize final types. Besides the function $\nu_w$ 
defined by
$$
\nu(i)= i- \#\{ a \leq g \co  w(a) \leq i\}=i-r_w(g,i)
$$
and extended by $\nu(2g-i)=\nu(i)-i+g$ for $i=0,\ldots,g$, we define
a function $\mu=\mu_w$ on the integers $1 \leq i \leq 2g$ by
$$
\mu(i) := \left(\max\{ w^{-1}(a)\co  1\leq a \leq i\} -g\right)^+,
$$
where $(x)^+:=\max(x,0)$. Alternatively, we have 
\begin{displaymath}
\mu(i)= \min\set{0\le j\le g}{r_w(g+j,i)=i}.
\end{displaymath}
Note that both $\mu$ and $\nu$ are non-decreasing
functions taking values between $0$ and $g$. Also $\nu$ can increase by $1$
only and $\nu(2g)=\mu(2g)=g$. If $w$ is final then $\nu_w$ is the final type
associated to $w$. For an arbitrary $w$ the function $\nu$ is the final type of
the final element in the coset $wS_g$. 

\begin{lemma} We have $\mu_w(i) \geq \nu_w(i)$ for $1 \leq i \leq 2g$ with
equality precisely when $w$ is a final element and then $\nu_w$ is the final
type of $w$.
\begin{proof} We first prove the inequality $\mu \geq \nu$. Let
$1\leq i \leq g$. Suppose that $\mu(i)=m$, i.e., the maximal $j$ with $w(j) $ in
$[1,i]$ is $g+m$. Then there are at most $m$ elements from $[g+1,2g] $ which map
into $[1,i]$ and so there are at least $i-m$ elements from $[1,g]$ with their
image under $w$ in $[1,i]$, so $i-\nu(i)\geq i-m$, in other words $\nu(i)\leq
\mu(i)$. For $i$ in the interval $[g+1,2g]$ we consider $\nu(2g-i)= \#\{ a \leq
g \co w(a) > i\}$. If $\mu(2g-i)=m$ then there are at least $g-m$ elements from
$[1,g]$ mapping into $[1,2g-i]$ so $\nu(2g-i)$ is at most equal to $m$.

If $w$ is final then $w$ respects the order on $[1,g]$ and this implies that if
$\# \{ a \leq g\co  w(a) \leq t\} =n$ then $t-n$ elements from $[g+1,2g]$ map to
$[1,t]$, so the maximum element from $[g+1,2g]$ mapping into $[1,t]$ is
$g+t-n$. Hence $\mu(t)=t-n=\nu(t)$.

Conversely, if $\mu(i)=\nu(i)$ then this guarantees that $w(i)< w(j)$ for all
pairs $1\leq i < j \leq g$.
\end{proof}
 \end{lemma}
We get as an immediate corollary.
\begin{corollary}\label{r-fulfillment}
Let $w \in W_g$. We have that $r_w(g+\nu_w(i),i)=i$ for all $1 \le i \le g$
precisely when $w$ is a final element.
\begin{proof}
The lemma says that if $w$ is final then we have 
\begin{displaymath}
\nu(i)=\mu(i)=\min\set{0\le j \le g}{r_w(g+j,i)=i}
\end{displaymath}
and in particular that $\nu(i) \in \set{0\le j \le
g}{r_w(g+j,i)=i}$ which gives one direction.

Conversely, if we have $r_w(g+\nu_w(i),i)=i$, then $\nu_w(i) \ge \mu_w(i)$ and
then the lemma gives that $w$ is final.
\end{proof}
\end{corollary}
\end{subsection}
\begin{subsection}{Canonical Types}

We now deal with an iterative way of constructing the function $\nu$ 
starting from its values on the endpoints and applying it repeatedly.

A final type $\nu$ is given by specifying $\nu(j)$ for $j=1,\ldots,2g$.  But it
suffices to specify the values of $\nu$ for the break points of $\nu$. Under
$\nu$ an interval $[i_1,i_2]$ between two consecutive break points of $\nu$ is
mapped to an interval of length $i_2-i_1$ or is mapped to one point. However,
the image points $\nu(i_1)$ and $\nu(i_2)$ need not be break points of
$\nu$. Therefore we enlarge the set of break points to a larger set $C_{\nu}$,
called the \Definition{canonical domain}. We define $C_{\nu}$ to be the smallest
subset of $\{0,1,\ldots,2g\}$ containing $0$ and $2g$ such that if $j \in
C_{\nu}$ then also $2g-j \in C_{\nu}$ and if $j\in C_{\nu}$ then $\nu(j) \in
C_{\nu}$. It is obtained by starting from $R=\{0,2g\}$ and adding the values
$\nu(k)$ and $\nu(2g-k)$ for $k\in R$ and continuing till this stabilizes.  The
restriction of $\nu$ to $C_{\nu}$ is called a \Definition{canonical type}. We
wish to see that the canonical domain $C_{\nu}$ contains the break points of
$\nu$ and hence that we can retrieve $\nu$ from the canonical type of $\nu$. To
see this we need a technical lemma (its formulation is somewhat obscured by the
fact that we also want to use it in another slightly different context).
\begin{definition-lemma}\label{0-1 slope}
We shall say that a subset $S \subseteq \{0,1,\dots,2g\}$ is \Definition{stable}
if it has the property that it contains $0$ and is stable under $i \mapsto
i^\perp:=2g-i$. For a stable subset $S$ a map \map{f}{S}{S\cap\{0,1,\dots,g\}}
is \Definition{adapted to $S$} if $f(0)=0$ and $f(2g)=g$, if it is
\Definition{contracting}, i.e., it is increasing and $f(j)-f(i)\le j-i$ for $i<j$ and
if it fulfills the following \Definition{complementarity condition}: For any two
consecutive $i,j \in S$ (i.e., $i<j$ and there are no $k\in S$ with $i<k<j$) we
have $f(j)-f(i)=j-i \implies f(j^\perp)=f(i^\perp)$.

i) If $S$ is stable and $f$ is a non-surjective function adapted to $S$ then
there is a proper subset $T \subset S$ such that $f_{|T}$ is adapted to $T$.

ii) If $S$ is stable and $f$ is a surjective function adapted to $S$ then for
any two consecutive $i,j \in S$ we have either $f(i)=f(j)$ or $f(j)-f(i)=j-i$.

iii) We say that $(S,f)$ is \Definition{minimally stable} if $S$ is stable and
$f$ is adapted to $S$ and furthermore there is no proper stable subset $T\subset
S$ for which $f_{|T}$ is adapted to it, then the function
\map{\nu}{\{1,2,\dots,2g\}}{\{1,2,\dots,g\}} obtained from $f$ by extending it
linearly between any two consecutive $i,j \in S$ is a final type, $S=C_\nu$ and
$\nu$ is the unique final extension of $f$. Conversely, if $f$ is the canonical
type of a final type $\nu$ then $(C_\nu,f)$ is minimally stable and in
particular $\nu$ is the linear extension of its canonical type.
\begin{proof}
For i) consider $T=f(S)\cup (f(S))^\perp$. It is clearly stable under $f$ and
$\perp$ and contains $0$. If $f$ is not surjective $T$ is a proper subset of
$S$.

Assume now that we are in situation of ii). We show that if $i<j \in S$ are
consecutive then either $f(j)-f(i)=j-i$ or $f(i)=f(j)$ by descending induction
on $j-i$.

By induction we are going to construct a sequence $i_k<j_k \in S$ $k=1,2,\dots$
of consecutive elements such that either
$(i_{k-1},j_{k-1})=(j_k^\perp,i_k^\perp)$ or $(f(i_k),f(j_k))=(i_{k-1},j_{k-1})$
but not both $(i_{k-1},j_{k-1})=(j_k^\perp,i_k^\perp)$ and
$(i_{k-2},j_{k-2})=(j_{k-1}^\perp,i_{k-1}^\perp)$ and in any case $j_k-i_k=j-i$.
We start by putting $i_1:=i, j_1:=j$. Assume now that $i_k<j_k$ have been
constructed. If we do not have $i_k,j_k\le g$, then as $g \in S$, we must have
$j_k^\perp,i_k^\perp\le g$ and then we put
$(i_{k+1},j_{k+1})=(j_k^\perp,i_k^\perp)$. If we do have $i_k,j_k\le g$ then by
the surjectivity of $f$ there are $i_{k+1},j_{k+1} \in S$ such that
$f(i_{k+1})=i_k$ and $f(j_{k+1})=j_k$. As $f$ is increasing $i_{k+1} < j_{k+1}$
and by choosing $i_{k+1}$ to be maximal and $j_{k+1}$ to be minimal we may
assume that they are neighbours. We must have that
$j_{k+1}-i_{k+1}=j_1-i_1$. Indeed, we have $f(j_{k+1})-f(i_{k+1})\le j_k-i_k$ as
$f$ is contracting. If we have strict inequality we have
$j-i=j_k-i_k<j_{k+1}-i_{k+1}$ and hence by the induction assumption we have
either that $j_k-i_k=f(j_{k+1})-f(i_{k+1})=j_{k+1}-i_{k+1}$ which is a
contradiction or $j_k=f(j_{k+1})=f(i_{k+1})=i_k$ which is also a
contradiction. Hence we have $j_{k+1}-i_{k+1}= j_k-i_k=j-i$ and we have verified
the required properties of $(i_{k+1},j_{k+1})$.

There must now exist $1\le k<\l$ such that $(i_k,j_k)=(i_{\l},j_{\l})$ and we pick
$k$ minimal for this property. If $k=1$ we have either
$j-i=j_{\l-1}-i_{\l-1}=f(j_{\l})-f(i_{\l})=f(j)-f(i)$ or
$j-i=j_{\l-2}-i_{\l-2}=f(j_{\l-1})-f(i_{\l-1})=f(i^\perp)-f(j^\perp)$ which
implies that $f(i)=f(j)$ by assumptions on $f$. We may hence assume that $k>1$.
We can not have both $(i_{k-1},j_{k-1})=(j_{k}^\perp,i_{k}^\perp)$ and
$(i_{\l-1},j_{\l-1})=(j_{\l}^\perp,i_{\l}^\perp)$ as that would contradict the
minimality of $k$. If $(i_{k-1},j_{k-1})=(j_{k}^\perp,i_{k}^\perp)$ and
$i_{\l-1},j_{\l-1})=(f(i_{\l}),f(j_{\l}))$ then we get
$j_{k-1}-i_{k-1}=j-i=j_{\l-1}-i_{\l-1}=f(i_{k-1}^\perp)-f(j_{k-1}^\perp)$ which
implies $f(j_{k-1})=f(i_{k-1})$ which is either what we want in case $k=2$ or a
contradiction. Similarly the case $(i_{\l-1},j_{\l-1})=(j_{\l}^\perp,i_{\l}^\perp)$ and
$i_{k-1},j_{k-1})=(f(i_k),f(j_k))$ leads to a contradiction as does the case
$i_{\l-1},j_{\l-1})=(f(i_{\l}),f(j_{\l}))$ and $i_{k-1},j_{k-1})=(f(i_k),f(j_k))$.

Finally, to prove iii) we note that by ii) for $i<j \in S$ consecutive we either
have $f(i)=f(j)$ or $f(j)-f(i)=j-i$. This means that the linear extension $\nu$
has the property that for $1\le i\le 2g$ we have either $\nu(i)=\nu(i-1)$ or
$\nu(i)=\nu(i-1)+1$ and if $\nu(i)=\nu(i-1)+1$ we get by the conditions on $f$
that $\nu(2g-i+1)=\nu(2g-i)$. If for some $i$ $\nu(i)=\nu(i-1)$ and
$\nu(2g-i+1)=\nu(2g-i)$ we get that $g=f(2g)=\nu(2g)<g$ which is impossible by
assumption and hence $\nu$ is indeed a final type. It is clear that $S$ fulfills
the defining property of $C_\nu$ so that $S=C_\nu$. The conditions on a final
element implies that $\nu(j)-\nu(i)\le j-i$ for $j<i$ which implies that $f$ has
a unique final extension.

Conversely, if $\nu$ is a final type then $C_\nu$ clearly fulfills the required
conditions and we have just noted that $\nu(j)-\nu(i)\le j-i$ for $j<i$. The
complementarity condition follows from the condition $\nu(i)=\nu(i-1) \iff
\nu(2g-i+1)=\nu(2g-i)+1$.
\end{proof}
\end{definition-lemma}

We now give an interpretation of the canonical domain in terms of
the Weyl group.
Let $v \in W_g$ be a final element. A \Definition{canonical fragment}
of $v$ is a maximal interval $\{i< a\le j\}=:]i,j] \subseteq
\{1,2,\dots,2g\}$ such that $v^k(]i,j])$ remains an interval for all
$k$.
\begin{proposition}\label{Fragment sorites}
Let $v \in W_g$ be a final element, $w$ its the complementary element, and $\nu$
the final type of $w$.

i) $\{1,2,\dots,2g\}$ is the disjoint union of the canonical fragments of $v$
and they are permuted by $v$.

ii) If $]i,j]$ is a canonical fragment for $v$, and if $\nu(j) \ne \nu(j-1)$,
then $\nu$ maps $]i,j]$ bijectively to $]\nu(i),\nu(j)]$.

iii) If $]i,j]$ is a canonical fragment, then so is
$]\overline{\jmath},\overline{\imath}]$.

iv) The upper endpoints of the canonical fragments together with $0$ comprise
exactly the canonical domain for $w$.
\begin{proof}
If two canonical fragments $I$ and $J$ meet, their union $K$ will be an interval
and as $v^k(K)=v^k(I)\cup v^k(J)$, we see that $v^k(K)$ will be an interval for
all $k$. By the maximality we get that $I=J$. On the other hand $]i-1,i]$
fulfills the stability condition so that $i$ lies in a fragment. Hence
$\{1,2,\dots,2g\}$ is the disjoint union of fragments.

Let now $R$ be the set of upper endpoints of fragments together with $0$.  As
$\{1,2,\dots,2g\}$ is the disjoint union of the fragments of $v$ it follows that
if $]i,j]$ is a fragment, then $i$ is also the upper end point of a
fragment. Thus it follows from iii) that $R$ is stable under $i \mapsto
\overline{\imath}$. Let now $i$ be an upper endpoint of a fragment. We want to
show that $\nu(i) \in R$ and we may certainly assume that $\nu(i)\ne 0$ and we
may also, by way of contradiction assume that $i$ is a minimal upper endpoint
for which $\nu(i)$ is not an upper endpoint. If $\nu(i)\ne \nu(i-1)$, then
$v^{-1}(i)=\nu(i)$ and hence $\nu(i)$ is an upper endpoint of a fragment. Hence
we may pick $j < i$ such that $\nu(i)=\nu(i-1)=\dots=\nu(j)\ne \nu(j-1)$. Then
$j$ can not belong to the same fragment as $i$ and thus there must be an upper
endpoint $j\le k <i$. Then $\nu(k)=\nu(i)$ and by minimality of $i$ we see that
$\nu(k)$ is an upper endpoint which is a contradiction.

We therefore have shown that $R$ contains $0$ and is stable under $i \mapsto
\overline{\imath}$ and $\nu$. Hence it contains the canonical domain. Let now $j
\in C_\nu \setminus\{0\}$ and let $i$ be the largest $j \in C_\nu$ such that
$i<j$. We now want show by induction on $k$ that $v^{-k}(I)$, $I:=]i,j]$,
remains an interval for all $k$ and that also $v^{-k}(j)$ is one of its
endpoints. Now, it follows from Lemma \ref{complementary} that
$C_\nu\setminus\{0\}$ is stable under $v$ and hence $v^{-k}(j)$ will be the only
element of $C_\nu$ in $v^{-k}(I)$. Under the induction assumption $v^{-k}(I)$ is an
interval with $v^{-k}(j)$ as one of its endpoints and hence $\nu$ is constant on
$v^{-k}(I)$ by Lemma \ref{0-1 slope}. By Lemma \ref{complementary} $v^{-1}$ maps
$v^{-k}(I)$ to interval with $v^{-k-1}(j)$ as one of its endpoints. This means
that $I$ is contained in a fragment and of the elements of $I$ it is only $j$
that can be contained in $R$. This means that there are no elements of $R$
between consecutive elements of $C_\nu$ and hence $R \subseteq C_\nu$.
\end{proof}
\end{proposition}

\begin{corollary}
By associating to a final type $\nu$ its canonical type, its Young diagram and
the element $w_{\nu}$ we obtain a 1-1 correspondence between the following sets
of cardinality $2^g$: the set of canonical types, the set of final types, the
set of final Young diagrams and the set of final elements of $W_g$.
\end{corollary}
\end{subsection}

\begin{subsection}{Admissible Elements}

The longest final element of $W_g$ is the element 
$$
w_{\emptyset}:=s_gs_{g-1}s_gs_{g-2}s_{g-1}s_g\ldots s_gs_1s_2s_3\ldots s_g,
$$
which as a permutation equals $[g+1,g+2,\dots,2g]$. Elements of $W_g$ that
satisfy $w \leq w_{\emptyset}$ are called \Definition{admissible}. We now
characterize these.
\begin{lemma}\label{below empty set}
i) An element $w \in W_g$ fulfills $w\leq w_{\emptyset}$ if and only
if $w(i)\leq g+i$ for all $1\leq i \leq g$.
                                                                                
ii) The condition that $w \le w_\emptyset$ is equivalent to $r_w(i,g+i)=i$ for
all $1 \le i \le g$.
\begin{proof}
The first part follows immediately from the description of the Bruhat-Chevalley order
(\ref{bruhat}) and the presentation of $w_\emptyset$.
                                                                                
For the second part one easily shows that $w(i)\leq g+i$ for all $1\leq i \leq
g$ is equivalent to $r_w(i,g+i)=i$ for all $1\leq i \leq g$ which gives the
first equivalence.
\end{proof}
\end{lemma}

\begin{remark}
The number of elements $w \in W_g$ with $w \leq w_{\emptyset}$  seems to be given by
\begin{displaymath}
\left(x\frac{d}{dx}\right)^g\left(\frac{1}{1-x}\right)_{\big|x=1/2}
\end{displaymath}
This is the number of necklaces of sets of $g$ labeled beads, see \cite[Seq.\
A000629]{S}. We leave this as a problem for the combinatorially inclined reader.
\end{remark}
We give an illustration of the various notions for the case $g=2$.
\noindent
\begin{example} $g=2$.
\noindent
The Weyl group $W_2$ consists of $8$ elements. We list the element, a reduced
expression as a word (i.e.\ a decomposition $w=s_{i_1}\cdots c_{i_k}$ with
$k=\ell(w)$), its length, the functions $\nu$ and $\mu$ and for final elements
we also give the partition defining the Young diagram.

$$ 
\begin{matrix}
&&  w & s & \ell & \nu & \mu & Y \\
\noalign{\hrule} \\
&& [4,3,2,1] & s_1s_2s_1s_2 & 4 & \{1,2\} & \{2,2\} \\
&&[4,2,3,1] & s_1s_2s_1 & 3&  \{1,1\} & \{ 2,2\}\\
&&[3,4,1,2] & s_2s_1s_2 & 3 & \{1,2\} & \{1,2\} & \emptyset\\
&&[2,4,1,3] & s_1s_2 & 2 & \{1,1\}  & \{1,1\}& \{ 1 \} \\
&&[3,1,4,2] & s_2s_1 & 2 & \{0,1\} & \{0,2\} \\
&&[2,1,4,3] & s_1 & 1 & \{0,0\} & \{0,0\} \\
&&[1,3,2,4] & s_2 & 1 & \{0,1\} & \{0,1\}& \{ 2 \} \\
&&[1,2,3,4] & 1 & 0 & \{0,0\} & \{0,0\} & \{ 1,2\} \\
\end{matrix}
$$
\end{example}
The orbits of the complementary element will be play an important r\^ole in our
discussion of the canonical flag. Here we introduce some definitions related to
them.
\begin{definition}
Let $w \in W_g$ be a final and let $v$ be its complementary element. Assume that
$S$ is an orbit of the action of $v$ on its fragments. As $v$ commutes with $i
\mapsto \overline{\imath}$ we have that $S$ is either invariant under $i \mapsto
\overline{\imath}$ in which case we say that it is an \Definition{odd orbit} or
that $\overline{S}$ is another orbit in which case we say $\{S,\overline{S}\}$
is an \Definition{even orbit pair}.
\end{definition}
\end{subsection}
\begin{subsection}{Shuffles}
Recall that a \Definition{$(p,q)$-shuffle} is a permutation $\sigma$ of
$\{1,2,\dots,p+q\}$ for which $\sigma(i) < \sigma(j)$ whenever $i<j \le p$ or $p
<i<j$. It is clear that for each subset $I$ of $\{1,2,\dots,g\}$ there is a
unique $(|I|,g-|I|)$-shuffle $\sigma^I$ such that
$I=\{\sigma^I(1),\sigma^I(2),\dots,\sigma^I(|I|)\}$ and we will call it the
\Definition{shuffle associated to $I$}. We will use the same notation for the
corresponding element in $W_g$ (i.e., fulfilling
$\sigma^I(2g+1-i)=2g+1-\sigma^I(i)$ for $1\le i \le g$). By doing the shuffling
from above instead of from below we get another shuffle $\sigma_I$ given by
$\sigma_I(i)=g+1-\sigma^I(g+1-i)$. We will use the same notation for its
extension to $W_g$. Note that $\sigma_I$ will shuffle the elements
$\{g+1,g+2,\dots,2g\}$ in the same way that $\sigma^I$ shuffles
$\{1,2,\dots,g\}$, i.e., $\sigma^I(g+i)=g+\sigma_I(i)$, which is the relation
with $\sigma^I$ that motivates the definition. Note that if
$I=\{i_1<\dots<i_r\}$ and if we assume that $i_r>r$ (if not then $\sigma^I$ and
$\sigma_I$ are the identity elements) and we let $k$ be the smallest index such
that $i_k>k$, then $\sigma^I=s_{i_k-1}\sigma^{I'}$ and
$\sigma_I=s_{g+1-(i_k-1)}\sigma_{I'}$, where
$I'=\{i_1,\dots,i_k-1,\dots,i_r\}$. We call the $s_{i-1} w s_{g+1-(i-1)}$ for $w
\in W_g$ the \Definition{$i$'th elementary shuffle} of $w$, say that $I'$ is the
\Definition{elementary reduction} of $I$ whose \Definition{reduction index} is
$i_k$.

We define the \Definition{height} of a shuffle associated to a subset
$\{i_1,i_2,\dots,i_k\} \subseteq \{1,2,\dots,g\}$ to be $\sum_s(i_s-s)$. 
Using $w'=s_iws_{g+1-t} \iff s_iw's_{g+1-t}=w$, we see that
starting with a shuffle $w$ of a final element we arrive at a final element
after $\height(w)$ elementary shuffles.
\begin{definition}
Let $Y$ be a final Young diagram of degree $g$. The \Definition{shuffles of $Y$}
are the elements of $W_g$ of the form $\sigma^Iw_Y\sigma_I^{-1}$ for $I
\subseteq \{1,2,\dots,g\}$.
\end{definition}
If $w \le w_\emptyset$ we say that $1 \le i \le g$ is a \Definition{semi-simple
index} for $w$ if $w(i)=g+i$ (note that as $w \le w_\emptyset$ we always have
$w(i)\le g+i$). The set of semi-simple indices will be called the
\Definition{semi-simple index set} and its cardinality the
\Definition{semi-simple rank}. We say that $w$ is \Definition{semi-simply final}
if the semi-simple index set has the form $[g-f+1,g]$ (where then $f$ is the
semi-simple rank). This is equivalent to $w$ having the form
$[\dots,2g-f+1,2g-f+2,\dots,2g]$. If $w=w_Y$, $Y$ a final Young diagram, then
$w$ is semi-simply final and the semi-simple rank is equal to $g$ minus the
length of the largest row of $Y$ (defined to be zero if $Y$ is empty).
\begin{proposition}\label{Etale shuffling}
Let $w \le w_\emptyset$ be a semi-simply final element of semi-simple rank $f$
and let $I \subseteq \{1,2,\dots,2g\}$ be a subset with $\#I = f$. Put
$\tilde I:= \set{g+1-i}{i \in I}$. Then $w':=\sigma^Iw\sigma_I^{-1}$ is an
element with $w' \le w_\emptyset$ of semi-simple rank $f$ and semi-simple index
set $\tilde I$. Conversely, all $w' \le w_\emptyset$ whose semi-simple index set
is equal to $I$ are of this form.
\begin{proof}
Put $j:=\sigma_I^{-1}(i)$. Note that $j > g-f \iff i \in \tilde I$. If $j > g-f$
we have $w(j)=g+j$ and hence
$\sigma^Iw_Y\sigma_I^{-1}(i)=\sigma^I(g+j)=g+\sigma_I(j)=g+i$. If on the other
hand $j\le g-f$, then if $w_Y(j)\le g$ there is nothing to prove, if not we may
write $w_Y(j)=g+k$ and as the semi-simple rank of $Y$ is $f$ we have $k<j$. Then
$\sigma^Iw_Y\sigma_I^{-1}(i)=\sigma^I(g+j)=g+\sigma_I(k)$ and as $k<j\le g-f$ we
have $\sigma_I(k)<\sigma_I(j)=i$ which gives $\sigma^Iw_Y\sigma_I^{-1}(i)<g+i$.

The converse is easy and left to the reader.
\end{proof}
\end{proposition}
Finally, we define the \Definition{$a$-number} for any element $w \in W_g$ by
$$
a(w):= r_w(g,g).
$$
If $w$ is final with associated Young diagram $Y$ then its $a$-number, 
also denoted $a_Y$, is the largest integer $a$ with
$0\leq a \leq g$ such that $Y$ contains the diagram with
$\xi=\{a,a-1,a-2,\ldots,1\}$.


\end{subsection}
\end{section}
\begin{section}{The Flag Space}
\label{sec: flag}
\begin{subsection}{The Flag space of the Hodge bundle}
In this section we introduce the flag space of a principally polarized
abelian scheme over a base scheme of characteristic $p$. We use the Frobenius
morphism to produce from a chosen flag on the de Rham cohomology a second flag,
the position of which with respect to the first flag will be the object of 
study. 

We let $S$ be a scheme (or Deligne-Mumford stack) in characteristic $p$ and let
${\Xcal}\to S$ be an abelian variety over $S$ with principal polarization
(everything would go through using a polarization of degree prime to $p$ but we
shall stick to the principally polarized case). We consider the de Rham
cohomology sheaf ${\Hcal}_{dR}^1({\Xcal}/S)$. It is defined as the hyper-direct
image $ {\Rcal}^1\pi_*(O_{\Xcal} \to \Omega_{{\Xcal}/S}^1)$.  It is a locally
free sheaf of rank $2g$ on $S$.  The polarization (locally in the \'etale
topology given by a relatively ample line bundle on ${\Xcal}/S$) provides us
with a symmetric homomorphism $\rho \colon {\Xcal} \to \hat {\Xcal}$ and the
Poincar\'e bundle defines a perfect pairing between $ {\Hcal}_{dR}^1({\Xcal}/S)
$ and ${\Hcal}_{dR}^1(\hat{\Xcal}/S)$ and thus ${\Hcal}_{dR}({\Xcal}/S)$ comes
equipped with a non-degenerate alternating form (cf.~\cite{O})
$$
\langle \, , \, \rangle \colon {\Hcal}_{dR}^1({\Xcal}/S) \times
{\Hcal}_{dR}^1({\Xcal}/S) \to O_S.
$$
Moreover, we have an exact sequence of locally free sheaves on $S$
$$
0 \to \pi_*(\Omega^1_{{\Xcal}/S}) \to {\Hcal}_{dR}^1({\Xcal}/S) \to
R^1\pi_*O_{\Xcal} \to 0.
$$
We shall write \Definition{$\protect\HH$} for the sheaf ${\Hcal}_{dR}^1({\Xcal}/ S)$ and
\Definition{$\protect\EE$} for the \Definition{Hodge bundle}
$\pi_*(\Omega^1_{{\Xcal}/S})$. We thus have an exact sequence
$$
0 \to \EE \to \HH \to \EE^{\vee} \to 0
$$
of locally free sheaves on $S$.
The relative Frobenius $F\co {\Xcal} \to {\Xcal}^{(p)}$ and the Verschiebung
$V\co {\Xcal}^{(p)} \to {\Xcal}$ satisfy $F\cdot V = p\cdot {\rm
id}_{{\Xcal}^{(p)}}$ and $V \cdot F= p \cdot {\rm id}_{\Xcal}$ and they induce
maps, also denoted $F$ resp.$V$, in cohomology:
$$
F\co \HH^{(p)} \to \HH \qquad {\hbox {\rm and } } \quad V\co  \HH \to \HH^{(p)}.
$$
Of course, we have $FV=0$ and $VF=0$ and $F$ and $V$ are adjoints (with respect
to the alternating form).  This implies that $ \im(F)= \ker (V)$ and $\im(V) =
\ker (F) $ are maximally isotropic subbundles of $\HH$ and
$\HH^{(p)}$. Moreover, since $dF=0$ on $\Lie({\Xcal})$ it follows that
$F=0$ on $\EE$ and thus $\im(V)=\ker(F) = \EE^{(p)}$.  Verschiebung thus provides
us with a bundle map (again denoted by $V$): $V\co \HH \to \EE^{(p)}$.
\par
\smallskip
Consider the space $\Fcal= \Flag(\HH)$ of symplectic flags on the bundle $\HH$
consisting of flags of subbundles $\{ \EE_i\}_{i=1}^{2g}$ satisfying
$\rk(\EE_i)=i$, $\EE_{g+i}=E_{g-i}^{\bot}$, and $\EE_g=\EE$. This space is a
scheme over $S$ and it is fibred by the spaces ${\Fcal}^{(i)}$ of partial flags
$$
\EE_i \subsetneq \EE_{i+1} \subsetneq \ldots \subsetneq \EE_g.
$$
So ${\Fcal}={\Fcal}^{(1)}= {\rm Flag}(\HH)$ and ${\Fcal}^{(g)}= S$ and there are
natural maps
$$
\pi_{i,i+1}\colon {\Fcal}^{(i)} \to {\Fcal}^{(i+1)},
$$
the fibres of which are Grassmann varieties of dimension $i$. So the relative
dimension of ${\Fcal} $ is $g(g-1)/2$. The space ${\Fcal}^{(i)}$ is equipped
with a universal partial flag. On ${\Fcal}$ the Chern classes of the bundle
$\EE$ decompose into their roots:
$$
\lambda_i= \sigma_i(\ell_1,\ldots, \ell_g)\qquad {\rm with} \qquad \ell_i =
c_1(\EE_{i}/ \EE_{i-1}),
$$
where $\sigma_i$ is the $i$-th elementary symmetric function.

On ${\Fcal}^{(i)}$ we have the Chern classes $\ell_{i+1}, \ldots,\ell_g$ and
$$
\lambda_j(i):= c_j(\EE_i) \qquad j=0,1,\ldots,i.
$$
Its Chow ring is generated over that of ${\Acal}_g$ by the monomials
$\ell_1^{m_1}\cdots \ell_g^{m_g}$ with $0 \leq m_j \leq j-1$.  For later use we
record the following Gysin formula.

\begin{formula}\label{Gysin} We have $(\pi_{i,i+1})_* \ell_{i+1}^k = s_{k-i}(i+1)$,
where $s_j(i+1)$ denotes the $j$th Segre class of $\EE_{i+1}$ ($j$th complete
symmetric function in the Chern roots $\ell_1, \ldots, \ell_{i+1}$).

\end{formula}

Given an arbitrary flag of subbundles
$$
0=\EE_0 \subsetneq \EE_1 \subsetneq\ldots \subsetneq \EE_g=\EE 
$$
with rank$(\EE_i)=i$ we can extend this uniquely to a symplectic filtration on
$\HH$ by putting
$$
\EE_{g+i}= (\EE_{g-i})^{\bot }.
$$
By base change we can transport this filtration to $\HH^{(p)}$.

We introduce a second filtration by starting with the isotropic subbundle
$$
\DD_g :=\ker(V)= V^{-1}(0) \subset \HH
$$
and continuing with
$$
\DD_{g+i}= V^{-1}(\EE_i^{(p)}).
$$
We extend it to a symplectic filtration by setting $\DD_{g-i}=
(\DD_{g+i})^{\bot}$. We thus have two filtrations $\EE_{\bullet}$ and
${\DD}_{\bullet}$ on the pullback of $\HH$ to ${\Fcal}$.

We shall use the following notation
$$
{\Lcal}_i= \EE_i/\EE_{i-1}\quad \text{and} \quad {\Mcal}_i= \DD_i/\DD_{i-1}
\qquad \text{ \rm for $1 \leq i \leq 2g$.}
$$
For ease of reference we formulate a lemma which follows immediately from
definitions.
\begin{lemma} 
We have ${\Mcal}_{g+i} \cong {\Lcal}_i^{(p)}$, ${\Lcal}_{2g+1-i} \cong
{\Lcal}_i^{\vee}$ and ${\Mcal}_{2g+1-i} \cong {\Mcal}_i^{\vee}$.
\end{lemma}
We shall more generally for a family $X \to S$ of principally polarized abelian
varieties say that a \Definition{Hodge flag} for the family is a complete
symplectic flag $\{\EE_i\}$ of $\HH$ for which $\EE_g$ is equal to the Hodge
bundle. By construction this is the same thing as a section of $\Fcal_g \to
S$. We shall also call the associated flag $\{\DD_i\}$ the \Definition{conjugate
flag} of the Hodge flag.
\end{subsection}

\begin{subsection}{The canonical flag of an abelian variety}
\label{ssec:canonical flag}

In this section we shall confirm that the \Emph{canonical filtration} of $X[p]$
(kernel of multiplication by $p$) by subgroup schemes of a principally polarized
abelian variety $X$ as defined by Ekedahl and Oort \cite{O1} has its analogue
for de Rham cohomology. Just as in \cite{O1} we do this in a family $\Xcal \to
S$. It is the coarsest flag which is isotropic (i.e., if $\DD$ is a member of
the flag then so is $\DD^\perp$) and stable under $F$ (i.e., if $\DD$ is a
member of the flag then so is $F(\DD^{(p)})$). The existence of such a minimal
flag is proven by adding elements $\Fcal^\perp$ and $F(\DD^{(p)})$ for $\DD$
already in the flag in a controlled fashion. We start by adding $0$ to the
flag. We then insist on three rules: 
\begin{enumerate}
\item If we added $\DD \subseteq \DD_g$, then we immediately add
$\DD^\perp$ (unless it is already in the flag constructed so far).

\item If we added $\DD_g \subseteq \DD$, then we immediately add
$F(\DD^{(p)})$ (unless it is already in the flag constructed so far).

\item If neither rule 1) nor rule 2) applies, then we add $F(\DD^{(p)})$ for the
largest element $\DD$ of the flag for which $F(\DD^{(p)})$ is not already in the
flag.
\end{enumerate}
We should not however do this construction on $S$; we want to insure that we get
a filtration by vector bundles: At each stage when we want to add the image
$F(\DD^{(p)})$, we have maps \map{F}{\DD^{(p)}}{\HH} of vector bundles and we
then have a unique minimal decomposition of the base as a disjoint union of
subschemes such that on each subscheme this map has constant rank. At the same
time as we add $F(\DD^{(p)})$ to the flag we replace the base by this disjoint
union. On this disjoint union $F(\DD^{(p)})$ then becomes a subbundle of $\HH$
and whether or not it is equal to the one of the previously defined subbundles
is a locally constant condition. A simple induction then shows that we get a
flag, i.e., for any two elements constructed one is included in the other, on a
disjoint union of subschemes of $S$. As each element added is either the image
under $V$ of an element previously constructed or the orthogonal of such an
element it is clear that this flag is the coarsening of any isotropic flag
stable under $V$ and it is equally clear that the decomposition of $S$ is the
coarsest possible decomposition. We shall call the (partial) flag obtained in
this way the \Definition{canonical flag} of $\Xcal/S$ and the decomposition of
$S$ the \Definition{canonical decomposition} of the base.

To each stratum $S'$ of the canonical decomposition of $S$ we associate a
canonical type as follows: Let $T \subseteq \{1,2,\dots,2g\}$ be the set of
ranks of the elements of the canonical flag and let
\map{f}{T}{T\cap\{1,\dots,g\}} be the function which to $t$ associates
$\rk(F(\DD^{(p)}))$, where $\DD$ is the element of the canonical flag of rank
$t$. We now claim that $T$ and $f$ fulfills the conditions of Lemma \ref{0-1
slope}. Clearly $T$ contains $0$ and by construction it is invariant under $i
\mapsto 2g-i$. Again by construction $f$ is increasing and has $f(0)=0$ and
$f(2g)=g$. Furthermore, if $i,j \in T$ with $i<j$ then $F$ induces a surjective
map $(\DD/\DD')^{(p)} \to F(\DD)/F(\DD')$, where $\DD$ resp.\ $\DD'$ are the
elements of the canonical flag for which the rank is $j$ resp.\ $i$ and hence
$f(j)-f(i)=\rk(F(\DD)/F(\DD'))\le\rk(\DD/\DD')=j-i$. Finally, assume that
$f(j)-f(i)=j-i$ and let $\DD$ and $\DD'$ be as before. Putting
$\DD_1:=F(\DD^{(p)})$ and $\DD'_1:=F(\DD'^{(p)})$ these are also elements of the
canonical filtration and we by assumption $F$ induces an isomorphism $F\co
(\DD/\DD')^{(p)} \to \DD_1/\DD'_1$. The fact that it is injective means that
$\DD\cap\ker F=\DD'\cap \ker F$ which by taking annihilators, and using that
$\ker F$ is its own annihilator, gives $\DD^\perp+\ker F=\DD'\perp+\ker F$ which
implies $F(\DD^\perp)=(\DD'^\perp)$ and hence that $f(2g-i)=f(2g-j)$. Now, if
$f$ is not surjective then by Lemma \ref{0-1 slope} there is a proper subset of
$T$ fulfilling the conditions of the lemma. This is not possible as $T$ by
construction is a minimal subset with these conditions. Hence $T$ and $f$
fulfills the conditions of Lemma \ref{0-1 slope} and hence by it we get that
$(f,T)$ is a canonical type. Let $\nu$ be its associated final type. If
$0=\DD_0\subset \DD_{i_1}\subset\dots\subset \DD_{2g}=\HH$ is the canonical flag
with $\rk\DD_i=i$ we have also proved that $F$ induces an isomorphism
$(\DD_j/\DD_i)^{(p)} \to \DD_{\nu(j)}/\DD_{\nu(i)}$ which can be rephrased as an
isomorphism
\begin{displaymath}
F\co \DD_{v(I)}^{(p)} \riso \DD_I,
\end{displaymath}
where we have used the notation $\DD_J:=\DD_i/\DD_j$ for an interval $J=]j,i]$
and $v \in W_g$ is the complementary element of (the final element of) $\nu$. We
shall say that $\nu$ (or more properly $f$) is the \Definition{canonical type}
of the principally polarised abelian variety $\Xcal_{S'} \to S'$. (We could
consider the canonical type as a locally constant function on the canonical
decomposition to the set of canonical (final) types.)
\begin{remark}
Note that the canonical flag is a flag containing $\DD_g$ and not $\EE_g$. That
will later mean that the canonical flag will be a coarsening of a conjugate flag
which is derived from a Hodge flag. On the one hand this is to be expected. As
the canonical flag is just that; it will be constructed in a canonical fashion
from the family of principally polarized abelian varieties. Hence it is to be
expected, and it is clearly true, that the canonical flag is horizontal with
respect to the Gauss-Manin connection. On the other hand we do not want to just
consider conjugate flags (or make constructions starting only with conjugate
flags). The reason is essentially the same; as $\DD_g$ (or more generally
the elements of the canonical flag) is horizontal it will not reflect first
order deformations whereas $\EE_g$ isn't and does. This will turn out to be of
crucial importance to us and is the reason why the Hodge flags will be the
primary objects and the conjugate flags derived. On the other hand, when working
pointwise, over an algebraically closed field say, we may recover the Hodge
flag from the conjugate flag and then it is usually most convenient to work with
the conjugate flag.
\end{remark}
\begin{example}\label{p-rank f example}
Let $X$ be an abelian variety with $p$-rank $f$ and $a(X)=1$ (equivalently, on $
G_g$ the operator $V$ has rank $g-1$ and semi-simple rank $g-f$). Then the
canonical type is given by the numbers $\{\rk(C_i)\}$, i.e.,
$$
\{ 0,f,f+1,\dots ,2g-f-1,2g-f,2g\}
$$
and $\nu$ is given by $\nu(f)=f$, $\nu(f+1)=f$, $\nu(f+2)=f+1, \,\dots,
\,\nu(g)=g-1,\dots,\, $ $\nu(2g-f-1)=g-1$, $\nu(2g-f)=g$ and $\nu(2g)=g$.  The
corresponding element $w \in W_g$ is $[f+1,g+1,\ldots,2g-f-1,2g-f+1,
\ldots,2g]$.
\end{example}
\end{subsection}
\end{section} 
\begin{section}{Strata on the Flag Space}

\begin{subsection}{The Stratification}

The respective positions of two symplectic flags are encoded by a combinatorial
datum, an element of a Weyl group.  We shall now define strata on the flag space
$\Fcal$ over the base $S$ of a principally polarized abelian scheme $X \to S$
that mark the respective position of the two filtrations $\EE_{\bullet}$ and
$\DD_{\bullet}$ that we have on the de Rham bundle over $\Fcal$.

Intuitively, the stratum ${\Uc}_w$ is defined as the locus of points $x$ such
that at $x$ we have
$$
\dim (\EE_i \cap \DD_j) \geq r_w(i,j) = \# \{ a \leq i\co  w(a) \leq j \} \quad
{\hbox {\rm for all } } \quad 1 \leq i,j \leq 2g.
$$
A more precise definition would be as degeneracy loci for some appropriate
bundle maps. While this definition would work fine in our situation where we are
dealing with flag spaces for the symplectic group it would not quite work when
the symplectic group is replaced by the orthogonal group on an even-dimensional
space (cf., \cite{F-P}).  With a view on future extensions of the ideas of this
paper to other situations we therefore adopt the definition that would work in
general. Hence assume that we have a semi-simple group $G$, a Borel group $B$ of
it, a $G/B$-bundle $T \to Y$ (with $G$ as structure group) over some scheme $Y$,
and two sections $s,t \co Y \to T$ of it. Then for any element $w$ of the Weyl
group of $G$ we define a (locally) closed subscheme $\Uo_w$ resp.~$\Uc_w$ of $Y$
in the following way. We choose locally (possibly in the \'etale topology) a
trivialization of $T$ for which $t$ is a constant section. Then $s$ corresponds
to a map $Y \to G/B$ and we let $\Uo_w$ (resp.~$\Uc_w$) be the inverse image of
the $B$-orbit $BwB$ (resp.~of its closure). Another trivialization will differ
by a map $Y \to B$ and as $BwB$ and its closure are $B$-invariant these
definitions are consistent and hence give global subschemes on $Y$. If $s$ and
$t$ have the property that $X=\Uo_w$, then we shall say that $s$ and $t$ are in
\Definition{relative position $w$} and if $Y=\Uc_w$ we shall say that $s$ and
$t$ are in \Definition{relative position $\le w$}.
\begin{remark}
The notation is somewhat misleading as it suggests that $\Uc_w$ is the closure
of $\Uo_w$ which may not be the case in general. In the situation that we shall
meet it will however be the case.
\end{remark}
The situation to which we will apply this construction is when the base scheme
is the space $\Fcal$ of symplectic flags $\EE_{\bullet}$ as above, $s$ is the
tautological section of the flag space of $\HH$ over $\Fcal$, and $t$ is the
section given by the conjugate flag $\DD_\bullet$. From now on we shall, unless
otherwise mentioned, let $\Uo_w$ and $\Uc_w$ denote the subschemes of $\Fcal$
coming from the given $s$ and $t$ and $w \in W_g$. In this case it is actually
often more convenient to use the language of flags rather than sections of
$G/B$-bundles and we shall do so without further mention. We shall also say that
a Hodge flag $\EE_\bullet$ is of \Definition{stamp $w$} resp.\ \Definition{stamp
$\le w$} if $\EE_\bullet$ and its conjugate flag $\DD_\bullet$ are in relative
position $w$ resp.\ $\le w$.
\begin{lemma}\label{bundle-isom} 
Over $\Uo_w$ we have an isomorphism ${\Lcal}_{i} \cong {\Mcal}_{w(i)}$ for all
$1\leq i \leq 2g$.
\end{lemma}
\begin{proof} By the definition of the strata we have that the
image of $\EE_{i} \cap \DD_{w(i)}$ has rank one greater than the ranks of
$\EE_{i-1} \cap \DD_{w(i)}$, $\EE_{i} \cap \DD_{w(i)-1}$, and
$\EE_{i-1} \cap \DD_{w(i)-1}$. So the maps $ \EE_i/\EE_{i-1} \leftarrow \EE_i
\cap \DD_{w(i)}/\EE_{i-1} \cap \DD_{w(i)-1} \to \DD_{w(i)}/\DD_{w(i)-1} $ give
the isomorphism.
\end{proof}
When the base of the principally polarized abelian scheme is $\Acal_g$ we shall
use the notation \Definition{$\protect\Fcal_g$} for the space of Hodge flags. Note that
a Hodge flag with respect to $X \to S$ is the same thing as a lifting over
$\Fcal_g \to \Acal_g$ of the classifying map $S \to \Acal_g$. The conjugate flag
as well as the strata $\Uo_g$ and $\Uc_g$ on $S$ are then the pullbacks of the
conjugate flag resp.\ the strata on $\Fcal_g$.
\end{subsection}
\begin{subsection}{Some Important Strata}
\label{ssec:Some Important Strata}

We now give an interpretation for some of the most important strata. To begin
with, if one thinks in terms instead of filtrations of $X[p]$ by subgroup schemes
it becomes clear that the condition $F(\DD_i) \subseteq \DD_i$ should be of
interest. It can almost be characterized in terms of the strata $\Uc_w$.
\begin{proposition}\label{V- and F-stability}
Let $X \to S$ be a family of principally polarized abelian varieties and
$\EE_\bullet$ a Hodge flag such that the flag is of stamp $\le w$ and $w$ is the
smallest element with that property.

i) For $j \le g$ we have that $r_w(i,g+j)=i$ precisely when $V(\EE_i) \subseteq
\EE_j^{(p)}$.

ii) For $j \le g$ we have that $r_w(g+j,i)=i$ implies that $F(\DD_i^{(p)})
\subseteq \DD_j$ and the converse is true if $S$ is reduced.

iii) We have that $V(\EE_i) \subseteq \EE_i^{(p)}$ for all $i$ precisely when $w
\le w_\emptyset$. If $S$ is reduced $F(\DD_i^{(p)}) \subseteq \DD_i$ for all $i$ precisely when $w
\le w_\emptyset$.
\begin{proof}
We have that $V(\EE_i) \subseteq \EE_j^{(p)} $ if and only if $\EE_i \subseteq
V^{-1}(\EE_j^{(p)})=\DD_{g+j}$. On the other hand, by definition $\rk \EE_i \cap
\DD_{g+j} \le r_w(i,g+j)$ with equality for at least one point of $S$. As $\rk
\EE_i \cap \DD_{g+j}=i \iff \EE_i \subseteq \DD_{g+j}$ we get the first part.

For the second part we start by claiming that $\EE_i^{(p)} \subseteq \DD_j$ is
implied by $F(\DD_i^{(p)}) \subseteq \DD_j$. Indeed, $F(\DD_i^{(p)}) \subseteq
\DD_j$ is equivalent to $F(\DD_i^{(p)})$ being orthogonal to $\DD_{2g-j}$, i.e.,
to the condition that for $u \in \DD_i^{(p)})$ and $v \in \DD_{2g-j}$ we have
$\scpr{Fu}{v}=0$. This implies that $0=\scpr{Fu}{v}={\scpr{u}{Vv}}^p$ and hence
$\scpr{u}{Vv}=0$ as $S$ is reduced which means that $\DD_i^{(p)} \subseteq
(V(\DD_{2g-j}))^\perp=(\EE_{g-j}^{(p)})^\perp=\EE_{g+j}^{(p)}$. As $S$ is reduced
this implies that $\DD_i\subseteq\EE_{g+j}$ and this in turn is equivalent to
$r_w(g+j,i)=i$. The argument can be reversed and then does not require $S$ to be
reduced.

Finally, we have from the first part that $V(\EE_i) \subseteq \EE_i^{(p)}$ for
all $i\le g$ precisely when $r_w(i,g+i)=i$ for all $i\le g$ but by induction on
$i$ that is easily seen to be equivalent to $w(i)\le g+i$ for all $i \le g$
which by definition means that $w\le w_\emptyset$. As $\EE_g^{(p)}=V(\HH)$ the
condition for $i>g$ is trivially fulfilled.

The proof of the second equivalence is analogous in that using ii) the condition
that $F(\DD_i^{(p)}) \subseteq \DD_i$ is equivalent to $r_w(g+i,i)=i$. In
general, $r_u(i,j)=r_{u^{-1}}(j,i)$ so that this condition is equivalent to
$r_{w^{-1}}(i,g+i)=i$ and hence by the same argument as before this condition
for all $i$ is equivalent to $w^{-1}\le w_\emptyset$. Now, Chevalley's
characterisation of the Bruhat-Chevalley order makes it clear that $u \le v \iff u^{-1}
\le v^{-1}$ and hence get $w^{-1}\le w_\emptyset \iff w \le
w_\emptyset^{-1}$. However, $w_\emptyset$ is an involution.
\end{proof}
\end{proposition}
\begin{remark} 
i) As we shall see the strata $\Uc_w$ in the universal case of $\Fcal_g$ are
reduced. 

ii) Flags of stamp $w \leq w_{\emptyset}$ are called admissible.
\end{remark}
We can also show that the relations between final and canonical types are
reflected for flags. We say that a Hodge flag is a \Definition{final flag} if it
is of stamp $w$ for a final element $w$. Also if $I=]i,j] \subseteq
\{1,2,\dots,2g\}$ is an interval and $\FF_\bullet$ is a complete flag of a
vector bundle of rank $2g$ then we define $\FF_I$ to be $\FF_j/\FF_i$.
\begin{proposition}\label{final and canonical flag}
Let $X \to S$ be a principally polarized abelian scheme over $S$ and
$\EE_\bullet$ is a final flag for it of stamp $w$.

i) The conjugate flag $\DD_\bullet$ is a refinement of the canonical flag. In
particular $w$ is determined by $X \to S$. More directly, we have that the final
type $\nu$ associated to $w$ is given by
\begin{displaymath}
\rk(\EE_g\cap\DD_i) = i-\nu(i)
\end{displaymath}
for all $i$. In particular the canonical decomposition of $S$ with respect to
$X \to S$ consists of a single stratum and its canonical type is the canonical
type associated to $w$.

ii) Conversely, assume that $S$ is reduced and that canonical decomposition of
$S$ consists of a single stratum, and let $\nu$ be the final type associated to
the canonical type of the canonical flag. Then any Hodge flag $\EE_\bullet$ whose
conjugate flag $\DD_\bullet$ is a refinement of the canonical flag and for which
we have $F(\DD_i^{(p)}) \subseteq \DD_{\nu(i)}$ for all $i$, is a final flag.

ii) If $I$ is a canonical fragment for $v$, the complementary element to $w$, then $F$
induces a bijection $(\DD_{v(I)})^{(p)} \riso \DD_I$.
\begin{proof}
We start by showing that $F(\DD_i)=\DD_{\nu(i)}$ for all $i$. Indeed, this is
equivalent to $F(\DD_i) \subseteq \DD_{\nu(i)}$ and
$\rk\bigl((\ker(F)=\EE_g^{(p)})\cap\DD_i^{(p)}\bigl)=i-\nu(i)$ as the second
condition says that $F(\DD_i)$ has rank $\nu(i)$. Now, the condition $F(\DD_i)
\subseteq \DD_{\nu(i)}$ is by Proposition \ref{V- and F-stability} implied by
$r_w(g+\nu(i),i)=\nu(i)$ which is true for a final element by
Corollary \ref{r-fulfillment}. On the other hand, the condition
$\rk(\EE_g^{(p)}\cap\DD_i^{(p)})=i-\nu(i)$ which is implied by
$r_w(g,i)=\rk(\EE_g\cap\DD_i)=i-\nu(i)$ which is true by definition of
$\nu$. Now, the fact that $F(\DD_i)=\DD_{j}$ for all $i$ and some $j$ (depending
on $i$) implies by induction on the steps of the construction of the canonical
flag that $\DD$ is a refinement of the canonical flag. The rest of the first
part then follows from what we have proved.

As for ii), assume first that $\EE_\bullet$ has a fixed stamp $w'$ and let $\nu
'$ be its final type. As $\DD_\bullet$ is an extension of the Hodge flag we get
that when $i$ is in the canonical domain of $\nu$ that $i-\nu '(i)=\rk(\EE_g\cap
\DD_i)=i-\nu(i)$ so that $\nu$ and $\nu '$ coincide on the canonical domain of
$\nu$ and hence they coincide by Lemma \ref{0-1 slope}.
The assumption that $F(\DD_i^{(p)}) \subseteq \DD_{\nu(i)}=\DD_{\nu '(i)}$ for
all $i$ is by Proposition \ref{V- and F-stability} equivalent to $r_w(g+\nu
'(i),i)=i$ for all $i$ and hence by Corollary \ref{r-fulfillment} gives that $w$
is final of type $\nu '=\nu$.

Finally, assume that $I=]i,j]$. The induced map $\DD_j^{(p)}/\DD_i^{(p)} \to
\DD_{\nu(j)}/\DD_{\nu(i)}$ is always surjective but it follows from Lemma
\ref{0-1 slope} that either the right hand side has dimension $0$ or it has the
same dimension as the right hand side. If they have the same dimension it
induces an isomorphism $\DD_I^{(p)} \riso \DD_{v^{-1}(I)}$. If the right hand
side has dimension zero, then again from Lemma \ref{0-1 slope} the two sides of
$\DD_{\overline{\jmath}}^{(p)}/\DD_{\overline{\imath}}^{(p)} \to
\DD_{\nu(\overline{\jmath})}/\DD_{\nu(\overline{\imath})}$ have the same
dimension and hence this map is an isomorphism and again is an isomorphism
$\DD_I^{(p)} \riso \DD_{v^{-1}(I)}$. As $I$ is an arbitrary fragment we
conclude.
\end{proof}
\end{proposition}
The number of final extensions of a canonical flag will now be expressed in the
familiar terms of the number of flags resp.\ self-dual flags in a vector space
over a finite field (resp.\ a vector space with a unitary form). Hence we let
\Definition{$\gamma_n^e(m)$} be the number of complete $\FF_{p^m}$-flags in $\FF_{p^m}^n$ and let
\Definition{$\gamma_n^o(m)$} be the number of complete $\FF_{p^{2m}}$-flags self-dual under
the unitary form
$\scpr{(u_1,\dots,u_n)}{(v_1,\dots,v_n)}:=u_1v_1^{p^m}+\cdots+u_nv_n^{p^m}$.
\begin{lemma}\label{canonical to final}
Let $X$ be a principally polarized abelian variety over an algebraically closed
field and $w \in W_g$ the element whose canonical type is the canonical type of
$X$. Put
\begin{displaymath}
\gamma(w)=\gamma_g(w):=\prod_{S=\overline{S}}\gamma^o_{\#I}(\#S/2)\prod_{\{S,\overline{S}\}}\gamma^e_{\#I}(\#S),
\end{displaymath}
where the first product runs over the odd orbits and the second over the even
orbit pairs and in both cases $I$ is a typical member of $S$.

The number of final flags for $X$ is then equal to $\gamma(w)$.
\begin{proof}
As we are over a perfect field, any symplectic flag extending $\DD_g$ is the
conjugate flag of a unique Hodge flag. Hence we get from Proposition \ref{final
and canonical flag} that a final flag is the same thing as a flag $\DD_\bullet$
extending the canonical flag and for which $F(\DD_i^{(p)})\subseteq
\DD_{\nu(i)}$. The condition that $\dim(\EE_g\cap \DD_i)=i-v(i)$ then gives that
we actually have $F(\DD_i^{(p)})= \DD_{\nu(i)}$. However, as $\DD_\bullet$
refines the canonical flag it is determined by the induced flags of the $\DD_I$
for $v$-fragments $I$ and the stability condition $F(\DD_i^{(p)})= \DD_{\nu(i)}$
transfers into the stability under the isomorphisms $F\co\DD_{v(I)}^{(p)} \riso
\DD_I$ of Subsection \ref{ssec:canonical flag}.

Hence the problem splits up into separate problems for each orbit under $v$ and
$i \mapsto \overline{\imath}$ (the map $i \mapsto \overline{\imath}$ transfers
into the isomorphism $\DD_{\overline{I}}\iso \DD_I^\vee$ induced by the
symplectic form). Given any fragment $I \in S$, $S$ a $v$-orbit of fragments
under $v$, the flag of $\DD_{v^k(I)}$ for any $k$ is then determined by the flag
corresponding to $I$ by the condition of $F^k$-stability takes the $k$'th
Frobenius pullback of the flag on $\DD_{v^k(I)}$ to the one of
$\DD_I$. Furthermore, the flag on $\DD_I$ then has to satisfy the consistency
condition of being stable under $F_S:=F^{\#S}$.

If now $\{S,\overline{S}\}$ is an even orbit pair the self-duality requirement
for the flag means that the flags for the elements of $\overline{S}$ are
determined by those for the elements of $S$ and given $I \in S$ there is no
other constraint on the flag on $\DD_I$ than the condition of stability
under $F_S$. Hence we have the situation of a vector bundle $\Dcal$ over our
base field $\k$ and an isomorphism $F_S\co \Dcal^{(p^m)} \riso \Dcal$, where
$m=\#S$, and we want to count the number of flags stable under $F_S$. Now, as
$\k$ is algebraically closed $\Dcal_p:=\set{v \in \Dcal}{F_S(v)=v}$ is an
$\FF_{p^m}$-vector space for which the inclusion map induces an isomorphism
$\k\Tensor\Dcal_p\riso \Dcal$. It then follows that $F_S$-stable flags
correspond to $\FF_{p^m}$-flags of $\Dcal_p$.

If instead $S$ is an odd orbit we will together with $I$ also have
$\overline{I}$ in $S$ and then, if $\#S=2m$, we will have that the flag on
$\DD_I$ must be mapped to the dual flag on $\DD_{\overline{I}}$ by
$F_S:=F^m$. The situation will be similar to the even orbit pair situation but
with a ``unitary twist'' as in Proposition \ref{p-unitary vb and local unitary
systems} and we get instead a correspondence with self-dual flags.
\end{proof}
\end{lemma}
\begin{example}
For the canonical type associated to the final type of example \ref{p-rank f
example} we have
$$
\gamma_g(C_{\bullet})= (p+1)(p^2+p+1) \cdots (p^{f-1}+p^{f-2}+\ldots +1).
$$
\end{example}
\begin{definition} We let $w=u_f$ be the final element
$$
u_f= s_g s_{g-1}s_g \cdots s_{g-f+1}\cdots s_gs_{g-f-1} \cdots s_g \cdots
s_1\cdots s_g,
$$
i.e., if we introduce $\tau_j=s_j s_{j+1}\cdots s_g$ we have $u_f= \tau_g
\tau_{g-1} \cdots \hat{\tau}_f \cdots \tau_1 $.  It corresponds to the Young
diagram consisting of one row with $g-f$ blocks and equals the element $w$ given
in \ref{p-rank f example}.
\end{definition}
Recall the notion of $a$-number $a(X)$ for an abelian variety $X$ of an
(algebraically closed) field $k$ of characteristic $p$, cf.\ \cite{O1}. It
equals the dimension over $k$ of the vector space ${\rm Hom}(\alpha_p,X)$ of
maps of the group scheme $\alpha_p$ to $X$. Equivalently, $a(X)$ equals the
dimension of the kernel of $V$ on $H^0(X,\Omega_X^1)$. In our terms
$H^0(X,\Omega_X^1)=\EE_g$ and $\ker V=\DD_g$ so that $a(X)=\dim \EE_g\cap
\DD_g$. The $p$-rank or \Definition{semi-simple rank} $f$ on the other hand can
be characterised by the condition that the $\dim_{\k} \cap_iV^{i}\HH$ of $V$ is
equal to $f$.
\begin{lemma}\label{prank=floci} 
i) Let $Y$ be a final Young diagram and $w \in W_g$ its final element and
assume that $x=(X,\EE_{\bullet},\DD_{\bullet})\in \Uo_w(\k)$, $\k$ a field. Then
the $p$-rank of $X$ equals the $p$-rank of $Y$.

ii) Let $w=u_f$ and $x=(X,\EE_{\bullet},\DD_{\bullet})\in {\Fcal}_g(\k)$. Then
we have $x \in \Uo_w$ (resp.\ ${\Uc_w}$) if and only if the filtration is $V$-stable and
the $p$-rank of $X$ is $f$ and the $a$-number of $X$ is $1$ (resp.\ the $p$-rank
of $X$ is $\leq f$). The image of $\Uo_w$ in ${\Acal}_g$ is the locus of abelian
varieties of $p$-rank $f$ and $a$-number $1$.
\begin{proof} 
By definition we have $r_w(i,g+i)=i$ for all $1 \le i \le g$ and
$r_w(i,g+i-1)=i$ precisely when $i \le g-f$. Hence by Proposition \ref{V- and
F-stability} $V(\EE_i)\subseteq \EE_i^{(p)}$ for all $1 \le i \le g$,
$V(\EE_i)\subseteq \EE_{i-1}^{(p)}$ for $i\le g-f$ and $V(\EE_i)\subsetneq
\EE_{i-1}^{(p)}$ for $i>g-f$. The first and last conditions means that $V$
induces an isomorphism $V\co \EE_{\{i\}} \riso  \EE_{\{i\}}$. On the other hand,
the second condition gives $V^{g-f}(\EE_{g-f})=0$. Together this gives that the
semi-simple rank of $X$ is $f$.

For the second part, since we have $w \leq w_\emptyset$ we must check the
condition on the $p$-rank and the $a$-number. By the definition on $\Uo_w$ we
have $x \in \Uo_w$ if and only if $\rk(\EE_g \cap \DD_j)=1$ for $f+1\leq j \leq g$ and
$0$ for $1\leq j \leq f$. This implies that the kernel of $V$ ($=\DD_g\cap
\EE_g$) has rank $1$ and the semi-simple rank of $V$ on $\EE_g$ is $f$. For $x
\in \Uc_w$ we get instead $\rk(\EE_g \cap \DD_j)\le 1$ for $f+1\leq j \leq g$
and $0$ for $1\leq j \leq f$.
\end{proof}
\end{lemma}
Also the strata $\Uo_w$ with $w \in S_g$ admit a relatively simple
interpretation. Recall that an abelian variety is called
\Definition{superspecial} if its $a$-number is equal to its dimension. This
happens if and only if the abelian variety (without polarization) is isomorphic
to a product of supersingular elliptic curves.

\begin{lemma} Let $x$ be a point of ${\Fcal}_g$ lying over $[X] \in {\Acal}_g$.
 The following are equivalent:
\begin{enumerate}
\item $x \in \cup_{w \in S_g} \Uo_w $.  

\item $\dim (\EE_g \cap \DD_g) \geq
g$. 

\item $\ker (V) =\EE_g$.  

\item The underlying abelian variety $X$ is superspecial.
\end{enumerate}
\begin{proof} It is well-known that $X$ is superspecial if and only if $V$
vanishes on $\EE_g$. But this is equivalent to $\dim (\EE_g \cap \DD_g) \geq
g$. This explains the equivalences of (2), (3) and (4).  Clearly, if $x \in
\Uo_w$ with $w \in S_g$ then $r_w(g,g)=g$, hence (2) holds.  Conversely, if $X$
is superspecial then any filtration $\EE_{\bullet}$ on $\EE_g$ is $V$-stable and
can be extended to a symplectic filtration.  Since the degeneracy strata for $w
\in S_g$, the Weyl group of ${\rm GL}_g$, cover the flag space of flags on
$\EE_g$ the lemma follows.
\end{proof}
\end{lemma}
\begin{lemma} Let $x$ be a point of $\Uo_w$ with underlying abelian variety $X$.
Then the $a$-number of $X$ equals $a(w)$.
Moreover, if $Y=\{1,2,\ldots,a\}$ with corresponding final element
$w_Y \in W_g$ then the image of $\Uo_{w_Y}$ in ${\Acal}_g$ is the locus $T_a$
of abelian varieties with $a$-number $a$.
\begin{proof}
The $a$-number of an abelian variety is by definition the dimension of the
kernel of $V$ on $H^0(X,\Omega_X^1)$. But this is equal to $r_w(g,g)=a(w)$.
The condition that $a(X)=a$ implies that $r_w(g,g)=a$, hence $\nu(g)=g-a$.  This
implies that $\nu(g-a+i)\geq i$ for $i=1,\ldots,a$. Therefore the `smallest'
$\nu$ satisfying these conditions is $\nu_{w_Y}$.
\end{proof}
\end{lemma}
\end{subsection}
\begin{subsection}{Shuffling flags}
\label{ssec:Shuffling flags}

Our first result on the stratification will concern the case when the $p$-rank
is positive. All in all the \'etale and multiplicative parts of the kernel of
multiplication by $p$ on the abelian variety have very little effect on the
space of flags on its de Rham cohomology. There is however one exception to
this. The most natural thing to do is to put the multiplicative part at the
bottom (and thus, by self-duality, the \'etale part at the top), which is what
automatically happens for a final filtration (on the conjugate filtration, that
is). We may however start with a final filtration and then ``move'' the
$\mu_p$-factors up-wards. Note, that over a perfect field the kernel of
multiplication by $p$ is the direct sum of its multiplicative, local-local, and
\'etale parts so that this is always possible. In general however it is possible
only after a purely inseparable extension. This means that we get an inseparable
map from a stratum where not all the $\mu_p$-factors are at the bottom to a
stratum where they all are. We intend to first give a combinatorial description
of the strata that can be obtained in this way from a final stratum and then to
compute the degree of the inseparable maps involved. However, as we have to
compute an inseparable degree we should work with Hodge filtrations instead of
conjugate filtrations as conjugate filtrations kill some infinitesimal
information. This causes a slight conceptual problem as the $V$-simple parts in
a final filtration are to be found ``in the middle'' rather than at the top and
bottom (recall that $V$ maps the top part of the conjugate filtration to the
bottom of the Hodge filtration). This will not be a technical problem but the
reader will probably be helped by keeping it in mind.


It turns out that the arguments used do not change if instead of considering
shuffles of final elements we consider shuffles of semi-simply final
elements. We shall treat the more general case as we shall need it later.

Hence we pick a subset $\tilde I \subseteq \{1,2,\dots,g\}$ and let
$\Uc_{\tilde I}^{ss}$ be the closed subscheme of $\Fcal_g$ defined by the
conditions that $V$ maps $\EE_i^{(p)}$ to $\EE_i^{(p)}$ for all $1\le i \le g$
and to $\EE_{i-1}^{(p)}$ for $i \notin \tilde I$. Hence $\Uc_w \subseteq
\Uc_{\tilde I}^{ss}$ precisely when $w \le w_\emptyset$ and the semi-simple index set of
$w$ is a subset of $\tilde I$. We also put
\begin{displaymath}
\Uo_{\tilde I}^{ss}:=\Uc_{\tilde I}^{ss}\setminus\cup_{\tilde I'\subset \tilde
I}\Uc_{\tilde I'}^{ss}
\end{displaymath}
so that $\Uo_w \subseteq \Uo_{\tilde I}^{ss}$ precisely when $w \le w_\emptyset$ and its
semi-simple index set is equal to $\tilde I$. If $I:=\set{g+1-i}{i \in \tilde
I}$ we get from Proposition \ref{Etale shuffling} that these $w$ are precisely
those of the form $\sigma^Iw'\sigma_I^{-1}$ for the semi-simply final $w'$.

We are now going construct, for every $I \subseteq \{1,2,\dots,g\}$, a morphism
\map{S_I}{\Uo_{\tilde I}^{ss}}{\Uo_{\{g-f+1,\dots,g\}}^{ss}}, where $\tilde
I:=\set{g+1-i}{i \in \tilde I}$ and $\#\tilde I=f$.

Let $\tilde\imath$ be the reduction index of the elementary reduction $I'$ of
$I$ and put $i:=g+1-\tilde\imath$. By Proposition \ref{Etale shuffling} we have
that $r_w(i+1,g+i)=i+1$ and $w(i)=g+i$. This means that if $\EE_\bullet$ is the
(tautological) Hodge flag on $\Uo_{\tilde I}^{ss}$, then $V(\EE_{i+1})\subseteq
\EE_i^{(p)}$ and $V(\EE_{i})\subsetneq \EE_{i-1}^{(p)}$ everywhere on
$\Uo_{\tilde I}^{ss}$. This means that $V$ gives a bijection on $\EE_{\{i\}}$
and is zero on $\EE_{\{i+1\}}$. Hence if $\Id$ denotes the map $s \mapsto
1\otimes s$ from $\Tcal$ to $\Tcal^{(p)}$ for the sheaves involved, the map
induced by the quotient map
\begin{displaymath}
\Ker(\Id-V)_{\EE_{\{i+1,i\}}} \to \Ker(\Id-V)_{\EE_{\{i\}}},
\end{displaymath}
where the kernel is computed in the \'etale topology on $\Uo_w$, is an
isomorphism. Also for a sheaf $\Tcal$ with a linear isomorphism $V\colon \Tcal
\to \Tcal^{(p)}$ we have an isomorphism $\Ker(V-I)\bigotimes \Ocal \to
\Tcal$. It follows that the short exact sequence
\begin{displaymath}
0 \to \EE_{\{i\}} \longrightarrow \EE_{\{i+1,i\}} \longrightarrow
\EE_{\{i+1\}} \to 0
\end{displaymath}
splits uniquely in a way compatible with $V$. This means that we may define a
new flag where $\EE_j'=\EE_j$ for $j \ne i$ and
$\EE_i'/\EE_{i-1}'=\EE_{i+1}/\EE_i$. Then the classifying map $\Uo_{\tilde
I}^{ss} \to \Fcal_g$ for this new flag will have its image in $\Uo_{\tilde
I'}^{ss}$. Repeating this process we end up with a flag whose classifying map
will have its image in $\Uo_{\{g-f+1,\dots,g\}}$ which by definition is our map
$S_I$
\begin{proposition}\label{Shuffling is finite radicial}
If $I \subseteq \{1,\dots,g\}$ then the map \map{S_I}{\Uo_{\tilde
I}^{ss}}{\Uo_{\{g-f+1,\dots,g\}}^{ss}} is finite, radicial, and surjective.
\begin{proof}
To get from a point of $\Uo_{\{g-f+1,\dots,g\}}^{ss}$ to one of $\Uo_{\tilde
I}^{ss}$ one has to find a $V$-invariant complement to some $\EE_i/\EE_{i-1}$ in
$\EE_{i+1}/\EE_{i-1}$. As $V$ will be zero on $\EE_{i+1}/\EE_{i}$ and bijective
on $\EE_i/\EE_{i-1}$, a complement over the fraction field of a discrete
valuation ring will extend to a complement over the discrete valuation ring (as
the complement can not meet $\EE_i/\EE_{i-1}$ over the special fibre) so that
the map is proper. It then remains to show that the map is a bijection over an
algebraically closed field. In that case $\EE_g$ splits canonically as a sum of
a $V$-nilpotent part and a $V$-semisimple part and the bijectivity is clear.
\end{proof}
\end{proposition}
In order to determine the degree (necessarily of inseparability) we shall do the
same factorisation as in the definition of $S_I$ so that we may consider the
situation of $\tilde I$ with $\tilde i$ and $I'$ being the reduction index
resp.\ elementary reduction of $I$. For the tautological flag $\EE_\bullet$ on
$\Uo_{\tilde I}^{ss}$ we have that $V$ is an isomorphism on $\EE_{\{i\}}$ and
zero on $\EE_{\{i+1\}}$ while the opposite is true on $\Uo_{\tilde I'}^{ss}$
\begin{lemma}\label{deg of insep}
The map $\Uo_{\tilde I}^{ss} \to \Uo_{\tilde I'}^{ss}$ is flat of degree $p$.
\begin{proof} 
We consider the partial symplectic flag space $\Fcal_g(i)$ consisting of the
flags of $\Fcal_g$ by removing the $i$'th member $\DD_i$ and its
annihilator. This means that we have a $\PP^1$-bundle $\Fcal_g \to
\Fcal_g(i)$. Now, under this map $\Uo_{\tilde I}^{ss}$ and $\Uo_{\tilde
I'}^{ss}$ map to the same subscheme $\Uo \subseteq \Fcal_g(i)$ and the map
$\Uo_{\tilde I}^{ss} \to \Uo_{\tilde I'}^{ss}$ is compatible with these
projections. Over $\Uo$ put $\Ecal:=\EE_{\{i,i+1\}}$, $\Mcal:=\ker(V\co \Ecal
\to \Ecal^{(p)})$, and $\Lcal:=\im(V\co \Ecal \to \Ecal^{(p)})$. Then on the
$\PP^1$-bundle \map{\pi}{\Fcal_g}{\Fcal_g(i)}, the subscheme $\Uo_{\tilde
I'}^{ss}$ is defined by the vanishing of the composite $\Ocal(-1) \to \pi^*\Ecal
\to \pi^*(\Ecal/\Mcal)$ and in fact gives a section of $\Fcal_g$ over $\Uo$
given by the sub-line bundle $\Mcal\subset \Ecal$. Hence it is enough to show
that the projection map $\Uo_{\tilde I}^{ss} \to \Fcal_g(i)$ is flat of degree
$p$. We have that $\Uo_{\tilde I'}^{ss}\subseteq \Fcal_g$ is defined by the
vanishing of the composite $\Ocal(-1)^{(p)} \to \pi^*\Ecal^{(p)} \to
\pi^*(\Ecal^{(p)}/\Mcal)$. It is then enough to show that $\Uo_{\tilde
I'}^{ss}\subseteq \Fcal_g$ is a relative Cartier divisor and for that is enough
to show that is a proper subset with each fibre of $\Fcal_g \to
\Fcal_g(i)$. This however is clear as for a geometric point of $\Fcal_g(i)$
there just two points that lie in $\Uc_\emptyset$, given by $\Mcal$ and
$\Lcal^{p^{-1}}$.
\end{proof}
\end{lemma}
Composing these maps we get
\begin{proposition}\label{Main shuffling}
Let $I \subseteq \{1,\dots,g\}$ and $\tilde I:=\set{g+1-i}{i \in I}$. Then the
map \map{S_I}{\Uo_I^{ss}}{\Uo_{\{1,\dots,g\}}^{ss}} is a finite purely
inseparable map of degree $p^{\height(I)}$.
\begin{proof}
The flatness and the degree of $S_I$ follows by factoring it by maps as in Lemma
\ref{deg of insep} and noting that the number of maps is $\height(I)$. The rest
then follows from Proposition \ref{Shuffling is finite radicial}.
\end{proof}
\end{proposition}
\begin{remark}
The result implies in particular that if $w'$ is a shuffle of $w$ by $I$, then
\map{S_I}{\Uo_{w'}}{\Uo_{w}} is flat and purely inseparable of degree
$p^{\height(I)}$. We shall later (see Corollary \ref{dim and smooth}) show that
$\Uo_{w'}$ and $\Uo_{w}$ are reduced. This shows that over the generic point of $\Uo_{w}$
each simple shuffle towards $w'$ really requires a finite inseparable extension
of degree $p$. This is a kind of non-degeneracy statement which is the
inseparable analogue of maximal monodromy (of which we shall also some examples
later on). It can also be seen as saying that a certain Kodaira-Spencer map is
injective.
\end{remark}
\end{subsection}

\begin{subsection}{The E-O strata on ${\Acal}_g \otimes \FF_p$}

\begin{definition} Let $w\in W_g$ be a final type. Then the E-O stratum
${\Vcal}_w$ associated to $w$ is the closed subset of ${\Acal}_g$ of points $x$
for which the canonical type of the underlying variety is equal to the canonical
type of $w$. We let $\Vc_w$ be the closure of $\Vo_w$.
\end{definition}

It is known that the dimension of ${\Vcal}_w$ is equal to $\dim(w)$, \cite{O1}.
This and the fact that the E-O strata form a stratification will also follow
from our results in \ref{sec:Local structure of strata} and \ref{E-O is strat}.
\end{subsection}

\end{section}
\begin{section}{Extension to the boundary}
\label{sec:boundary}

The moduli space ${\Acal}_g$ admits several compactifications. The Satake or
Baily-Borel compactification ${\Acal}_g^*$ is in some sense minimal, cf.\
\cite{F-C}.  It is a stratified space
$$
{\Acal}_g^*= \cup_{i=0}^g {\Acal}_i.
$$
Chai and Faltings define in [F-C] a class of smooth toroidal compactifications.
If $\tilde{\Acal}_g$ is such a toroidal compactification then there is a natural
map $q\co  \tilde{\Acal}_g \to {\Acal}_g^*$.  This induces a stratification of
$\tilde{\Acal}_g$:
$$
\tilde{\Acal}_g= \cup_{i=0}^g q^{-1}({\Acal}_{g-i})= \cup_{i=0}^g {\Acal}_g^{\langle
i \rangle}.
$$
The stratum ${\Acal}_g^{\langle i \rangle}$ parametrizes the semi-abelian
varieties of torus rank $i$.

The Hodge bundle $\EE$ on ${\Acal}_g$ can be extended to a rank $g$ vector
bundle, again denoted by $\EE$, on $\tilde{\Acal}_g$.  On ${\Acal}_g^{\langle i
\rangle}$ the Hodge bundle fits into an exact sequence
$$
0 \to \EE^{\prime} \to \EE \to \EE^{\prime\prime} \to 0,
$$
where the $\EE^{\prime}$ is a rank $g-i$ bundle and $ \EE^{\prime\prime} $ can
be identified with cotangent bundle along the identity section of the toric
part of the semi-abelian variety over ${\Acal}_g^{\langle i \rangle}$. The
bundle $\EE'$ is the pull back under $q\co  {\Acal}_g^{\langle i \rangle} \to
{\Acal}_g^*$ of the Hodge bundle on ${\Acal}_{g-i}$.

The Verschiebung $V$ acts in a natural way on the extension $\EE$ and it
preserves $\EE^{\prime}$. It induces an action on $\EE^{\prime \prime}$ with
trivial kernel because $\EE^{\prime \prime}$ comes from the toric part and is
generated by logarithmic forms.

The de Rham bundle $\HH$ on ${\Acal}_g$ also admits an extension.  In
characteristic $0$ this is the logarithmic de Rham sheaf
$R^1\pi_*(\Omega^{\bullet}_{\tilde{\Xcal}_g/\tilde{\Acal}_g}(\log))$, where the
$\log$ refers to the acceptance of logarithmic singularities along the divisor
at infinity, cf.\ \cite[Thm. VI:1.1]{F-C}.
We have a short exact sequence
$$
0 \to \EE \to {\HH} \to \EE^{\vee} \to 0
$$
extending the earlier mentioned sequence on ${\Acal}_g$.

We now want to compare the E-O strata on ${\Acal}_g$ and $\tilde{\Acal}_g$ and
for this we introduce some notation.  For a given integer $1 \leq i \leq g$ we
can consider the Weyl group $W_{g-i}$ as a subgroup of $W_g$ by letting it act
on the set $\{i+1,i+2,\ldots,g,\ldots,2g-i\}$ via the bijection $j
\longleftrightarrow i+j$ for $1\leq j \leq g-i$. More precisely, define
$\rho_i\co W_{g-i} \to W_g$ via
$$
\rho_i(w)(l)= \begin{cases} i+w(l) & \text{for $1 \leq l \leq g-i$} \\ g+l &
\text{for $g-i+1 \leq l \leq g$.} \\
\end{cases}
$$
This map respects the Bruhat-Chevalley order and final elements are mapped to final
elements.

As symplectic flags on $\HH$ are determined by their restriction to $\EE$ and as
we can extend $\EE$ to $\tilde{\Acal}_g$, we can extend ${\Fcal}_g$ to a flag
bundle $\tilde{\Fcal}_g$ on $\tilde{\Acal}_g$.  Then we can also extend the
degeneracy loci $\Uo_w$ and $\Uc_w$ to $\tilde{\Fcal}_g$ and we shall use the
same notation for these extensions.

Similarly, we can define the notion of a canonical filtration for a semi-abelian
variety. If $1\to T \to A \to A'\to 0$ is a semi-abelian variety with abelian
part $A'$ and toric part $T$ of rank $t$ and if the function $\nu'$ on $\{0, c_1
, \ldots , c_r , c_{r+1} , \ldots , c_{2r}=2 \dim(A')\}$ is the canonical type
of $A'$ then we define the canonical type of $A$ to be the function $\nu$ on
$$
\{ 0 , t , t+c_1 , \ldots , t + c_r , t +c_{r+1} , \ldots , t + c_{2r}, 2g-t,
2g\}
$$
defined by $\nu(t+c_i)=t+\nu'(c_i)$.  Using this definition we can extend the
E-O stratification to $\tilde{\Acal}_g$.

The stratification ${\Acal}^{\langle i \rangle}$ of $\tilde{\Acal}_g$ induces a
stratification $\tilde{\Fcal}_g^{\langle i \rangle}$ by flag spaces
${\Fcal}_g^{\langle i \rangle}$ over ${\Acal}_g^{\langle i \rangle}$.  Recall
that ${\Acal}_g^{\langle i \rangle}$ admits a map $q \colon {\Acal}_g^{\langle i
\rangle} \to {\Acal}_{g-i}$ induced by the natural map $\tilde{\Acal}_g \to
{\Acal}_g^*$. Similarly, we have a natural map $\pi_i=\pi \colon
{\Fcal}_g^{\langle i \rangle} \to {\Fcal}_{g-i}$ given by restricting the
filtration on $\EE$ to $\EE^{\prime}$.

We now describe the interplay between the two stratifications
$({\Fcal}_g^{\langle i \rangle})_{i=1}^g$ and $(\Uo_w)_{w \in W_g}$.
\begin{lemma}\label{strata at the boundary}
Let $w \in W_g$ be an element with $w\leq w_{\emptyset}$.

i) We have: $\Uo_w \cap {\Fcal}_g^{\langle i \rangle} \neq \emptyset$ if and
 only if $w$ is a shuffle of an element in $ \rho_i(W_{g-i})$.  

ii) If $w=\rho_i(w')$ with associated degeneracy loci $\Uo_w \subset {\Fcal}_g $
and $\Uo_{w'} \subset {\Fcal}_{g-i}$ then we have $\Uo_w \cap {\Fcal}_g^{\langle
i \rangle} = \pi_i^{-1}(\Uo_{w'})$.

iii) In particular $\Uo_w$ is the closure of its intersection with $\Acal_g$.

iv) At a point of $\tilde\Acal_g$ for which the torus part of the ``universal''
semi-abelian variety has rank $r$ there is a smooth map from the formal
completion of $\tilde\Acal_g$ at the point to the formal multiplicative group
$\hat\GG_m^r$ such that the locus where the torus rank of the universal
semi-abelian variety is $s\le r$ the inverse image of the locus of points of
$\hat\AA_m^r$ where $r-s$ coordinates are $1$. The restriction of this map to
any $\Uc_w$ is smooth.
\end{lemma}
\begin{proof}
A $V$-stable filtration on $\EE$ restricts to a $V$-stable filtration on
$\EE'$. If $\Uo_w \cap {\Fcal}_g^{\langle i \rangle}$is not empty then it
determines a $w' \in W_{g-i}$ such that $\Uo_w \cap {\Fcal}_g^{\langle i
\rangle} \subseteq \pi_i^{-1}(\Uo_{w'})$. Since $V$ is invertible on
${\EE}^{\prime \prime}$ one sees that $w$ is a shuffle of $\rho_i(w')$ and that
$\Uo_w \cap {\Fcal}_g^{\langle i \rangle}= \pi_i^{-1}(\Uo_{w'})$.

For the third part we note that the dimension of the intersection of $\Uo_w$ with
the boundary has dimension one lower than that of $\Uo_w$ itself.

The last part is a direct consequence of the local construction of
$\tilde\Acal_g$ using toroidal compactifications and of the universal
semi-abelian variety using Mumford's construction where it is defined by taking
the quotient of a semi-abelian variety by a subgroup of the torus part, the
subgroup being generated by the coordinate functions of $\hat\GG_m^r$ (see
\cite{F-C} for details). As $\HH$ and $\EE$ only depend on that fixed
semi-abelian variety it is clear that the restriction of the map to a $\Uc_w$ is
smooth.
\end{proof}
Note also that this lemma is compatible with shuffling. It also results from the
lemma that we can define the E-O stratification on the Satake compactification
by either considering the closure of the strata $\Vo_w$ on ${\Acal}_g$ or by
considering the images of the final strata $\Vo_w$ on $\tilde{\Acal}_g$.
\end{section}
\begin{section}{Existence of boundary components}

Our intent in this section is to show the existence of points in the smallest
possible stratum $\Uc_1$, the stratum associated to the identity element of
$W_g$.

\begin{proposition}\label{lowest degeneration}
Let $X$ be an irreducible component of any $\Uc_w$ in ${\Fcal}_g$. Then $X$
contains a point of $\Uc_1$.
\begin{proof}
We prove this by induction over $g$ and over the Bruhat-Chevalley order of $w$. The
statement is clear for $g=1$.  We start off by choosing a Chai-Faltings
compactification $\tilde{\Acal}_g$ of $\Acal_g$ with a semi-abelian family over
it (and a ``principal'' cubical structure so that we get a principal
polarization on the semi-abelian variety modulo its toroidal part).

What we now actually want to prove is the same statement as in the proposition
but for $\tilde{\Fcal_g}$ instead. Since $\Uc_1$ is contained in ${\Fcal}_g$ the
result will follow.  We start off by considering the intersection of $X$ with
the boundary $\tilde{\Fcal}_g -{\Fcal}_g$ to get $Y$. Assume that $Y$ is
non-empty and irreducible by replacing it by an irreducible component of $Y$ if
necessary. Then $Y$ is contained in $\pi_1^{-1}(\Uc_{w'})$ with $\rho_1(w')=w$
for some $w' \in W_{g-1}$.  We claim that $Y$ now is an irreducible component of
$\pi_1^{-1}(\Uc_{w'})$.  This follows from the fact that ``we can freely move
the toroidal part into an abelian variety'' which is Lemma \ref{strata at the
boundary}.

By induction on $g$ we can assume that $\Uc_{w'}$ in ${\Fcal}_{g-1}$
contains $\Uc_{1'}$, where $1'$ is the identity element of $W_{g-1}$.  Any
component $Z$ of $\Uc_{\rho_1(1')}$ which lies in $X$ and meets $Y$ does not lie
completely in the boundary $\tilde{\Fcal}_g -{\Fcal}_g$. By induction on the
Bruhat-Chevalley order we can assume that $w=\rho_1(1')$ and $X=Z$. Note also that for any
$w^{\prime \prime} < w$ we have that $\Uc_{w^{\prime \prime}}$ does not meet the
boundary and by induction on the Bruhat-Chevalley order we get that $X=\Uo_w\cap X$. On
the other hand if $X$ does not meet the boundary we immediately get the same
conclusion.

Hence we may and shall assume that $Y$ has the property that it lies completely
inside $\Uo_w$ and that it is proper. Lemma \ref{Raynaud trick} now shows that
it has an ample line bundle of finite order which together with properness
forces $Y$ to zero-dimensional. Now we always have that $\dim Y \ge \ell(w)$
(the proof is analogous to \cite[Them 14.3]{Ful}) which gives $\ell(w)=0$ and so
$w =1$ which clearly is a trivial case.
\end{proof}
\end{proposition}
\begin{lemma}\label{Raynaud trick}
Suppose that $X$ is a proper irreducible component of $\Uc_w$ inside ${\Fcal_g}$
such that $X\cap \Uo_w = X$. Then $X$ is $0$-dimensional.
\end{lemma}
\begin{proof} 
This is a version of the Raynaud trick. We have the variety $X$ and two
symplectic flags $\EE_\bullet$ and $\DD_\bullet$ which at all points of $X$ are
in the same relative position $w$. It follows from \ref{bundle-isom} that we
have an isomorphism between $\Lcal_i:=\EE_i/\EE_{i-1}$ and
$\Mcal_{w(i)}:=\DD_{w(i)}/\DD_{w(i)-1}$ over $X$ and then as we also have
isomorphisms between $\Lcal_i^p$ and $\Mcal_{g+i}$ and $\Lcal_i$ and
$\Lcal_{2g+1-i}^{-1}$ we conclude that all the $\Lcal_i$ have finite order. On
the other hand we know that
$\Lcal_{2g}\bigotimes\Lcal_{2g-1}\bigotimes\cdots\bigotimes\Lcal_{g+1}$ is
relatively ample and
$\Lcal_{g}\bigotimes\Lcal_{g-1}\bigotimes\cdots\bigotimes\Lcal_{1}$ is ample on
the base $\Acal_g$ we conclude.
\end{proof}
\end{section}
\begin{section}{Superspecial fibres}

We shall now discuss the fibre of $\Fcal_g \to \Acal_g$ over superspecial
points. The superspecial abelian varieties are characterized by the condition
that $\EE_g=\DD_g$, i.e., the strata $\Uo_w$ for which $w \in S_g$. Furthermore,
$V$ induces an isomorphism $\EE/\EE_g \riso \EE_g^{(p)}$. On the other hand the
polarization gives an isomorphism $(\EE_g)^* \riso \EE/\EE_g$. This leads to the
following definition.
\begin{definition}
i) Let $S$ be a scheme in characteristic $p$. A \Definition{$p$-unitary} vector
bundle is a vector bundle $\Ecal$ over $S$ together with an isomorphism
$F^*\Ecal\riso \Ecal^*$, where \map FSS is the (absolute) Frobenius map.

ii) Let $\Ecal$ be a $p$-unitary vector bundle over $S$ and let $P \to S$ be the
bundle of complete flags on $\Ecal$. The \Definition{$p$-unitary Schubert
strata} of $P$ is the stratification given by letting $\Uo_w$, $w \in S_g$,
consist of the points for which the universal flag $\Fcal$ and the dual of the
Frobenius pullback $(F^*\Fcal)^*$ are in position corresponding to $w$.
\end{definition}
A map $F^*\Ecal\to \Ecal^*$ of vector bundles is the same thing as a map
$F^*\Ecal\bigotimes_{\Ocal_S}\Ecal \to \Ocal_S$, which in turns corresponds to a
bi-additive map \map{\scprod--}{\Ecal\times\Ecal}{\Ocal_S} fulfilling
$\scprod{fa}{b}=f^p\scprod{a}{b}$ and $\scprod{a}{fb}=f\scprod{a}{b}$. We shall
normally use this latter description.

All $p$-unitary vector bundles are trivial in the \'etale topology as the
following proposition shows.
\begin{proposition}\label{p-unitary vb and local unitary systems}
If $\scprod--$ is a $p$-unitary structure on the vector bundle $\Ecal$ then
\begin{displaymath}
E := \set{a \in \Ecal}{\forall b \in \Ecal\co\scprod{b}{a}={\scprod{a}{b}}^p}
\end{displaymath}
is a local (in the \'etale topology) system of $\FF_{p^2}$-vector
spaces. Furthermore, $\scprod--$ induces a unitary (wrt the involution $(-)^p$
on $\FF_{p^2}$) structure on $E$. Conversely, if $E$ is a local system of
$\FF_{p^2}$-vector spaces with a unitary structure, then $\Ecal :=
E\bigotimes_{\FF_{p^2}}\Ocal_S$ is a $p$-unitary vector bundle.

These two constructions establish an equivalence between the categories of
$p$-unitary vector bundles and that of local systems of unitary
$\FF_{p^2}$-vector spaces. In particular all $p$-unitary vector bundles (of the
same rank) are locally isomorphic in the \'etale topology.
\begin{proof}
The pairing $\scprod--$ gives rise to an isomorphism $\psi '\co\Ecal \riso
(F^*\Ecal)^*$ by $a \mapsto b \mapsto \scprod{b}{a}$ and an isomorphism
$\psi\co(F^2)^*\co\Ecal \riso (F^*\Ecal)^*$ by $a \mapsto b \mapsto
{\scprod{a}{b}}^p$. The composite $\rho := \psi^{-1}\circ\psi '$ thus gives an
isomorphism $\Ecal \riso (F^2)^*\Ecal$. Then $E$ is simply the kernel of
$\rho-1\otimes id$ and the fact that $\Ecal = E \bigotimes_{\FF_{p^2}}\Ocal_S$
follows from \cite{hasse63::zyklis+erweit+primz+funkt+charak}. The pairing
$\scprod--$ then induces a unitary pairing on $E$ which is perfect as
$\scprod--$ is. Conversely, it is clear that a unitary pairing on $E$ translates
to one on $\Ecal$.

Finally, as all (perfect) unitary mappings on $\FF_{p^2}$-vector spaces of fixed
dimension are isomorphic we get the local isomorphism.
\end{proof}
\end{proposition}
This proposition has the following immediate corollary.
\begin{corollary}
The flag variety fibrations of two $p$-unitary vector bundles of the same rank
on the same base are locally isomorphic by an isomorphism preserving the
$p$-unitary Schubert strata.
\end{corollary}
Let us now consider the situation where the base scheme is ${\rm
Spec}(\FF_{p^2})$ and $\Ecal$ is a $\FF_{p^2}$-vector space given a unitary
perfect pairing. If we consider the smallest unitary Schubert stratum then it
corresponds to flags that coincide with their unitary dual. Taking duals once
more we see that they are taken to themselves after pullback by the square of
the Frobenius, hence are defined over $\FF_{p^2}$. Furthermore, they are
selfdual with respect to the unitary pairing. This should come as no surprise as
that stratum corresponds to final filtrations on superspecial abelian
varieties. The next to lowest strata are somewhat more interesting.
\begin{lemma}\label{next to lowest}
Let $V$ be a $g$-dimensional $\FF_{p^2}$ unitary vector space and let $\PP$ be
projective space based on $V$. If $s=(i,i+1)\in S_g$ for some $1 \le i < g$ then
the closed Schubert stratum $\Uo_s \subseteq \PP$ consists of the flags $0=\EE_0
\subset \EE_1 \subset \dots \subset \EE_g$ where the $\EE_j$ for $j \ne i,g-i$
are $\FF_{p^2}$-rational with $\EE_j^\perp = \EE_{g-j}$,
$\EE_{g-i}=(\EE_i^{(p)})^\perp$ unless $i=2g$, and $\EE_{i} \ne
(\EE_{g-i}^{(p)})^\perp$
\begin{proof}
If $\EE_\bullet$ and $\DD_\bullet$ are two flags in $V\bigotimes R$ ($R$ some
$\FF_{p^2}$-algebra) and they are in position $s$, then $\dim(\EE_j\cap \DD_j) =
r_s(j,j)=j$ for $i \ne j$, i.e., $\EE_j=\DD_j$ and for $i$ the conditions give
us $\EE_i\cap\DD_i=\EE_{i-1}$. In our case, where $\DD_j=F^*\EE_{g-j}^\perp$,
this means $\EE_j=F^*\EE_{g-j}^\perp$ for $j\ne i$ and $\EE_i\ne
(\EE_{g-i}^{(p)})^\perp$. If also $j \ne g-i$ we can use this twice and get that
$\EE_j=\EE_j^{(p^2)}$, i.e., $\EE_j$ is $\FF_{p^2}$-rational.
\end{proof}
\end{lemma}
As a result we get the following connectedness result, analogous to
\cite[Prop.~7.3]{O1}.
\begin{theorem}\label{fibre connectedness}
Let $V$ be a $g$-dimensional $\FF_{p^2}$ unitary vector space and let $\PP$ be
projective space based on $V$. Let $S \subseteq \{1,\dots,g-1\}$. Let $\Uc$ be
the union of the $\Uc_{s_i}$ for $i \in S$. Then two flags $0=A_0\subset
A_1\subset\dots\subset A_{g-1}\subset A_g$ and $0=B_0\subset
B_1\subset\dots\subset B_{g-1}\subset B_g$ in $\Uc_1$ lie in the same component
of $\Uc$ precisely when $B_i=A_i$ for all $i\notin S$. Furthermore, every
connected component of $\Uc$ contains an element of $\Uc_1$.
\begin{proof}
The last statement is clear as every irreducible component of any $\Uc_{s_i}$
contains a point of $\Uc_1$. This follows from (\ref{lowest degeneration}) but
can also easily be seen directly. 

We start by describing looking at the locus $\Uc_F^i$ of a $\Uc_{s_i}$ with $1\le i <g$ of
flags for which all the components of the flag except the dimension $i$ and
dimension $g-i$ parts are equal to a fixed (partial) $\FF_{p^2}$-rational
self-dual flag $F_\bullet$, the claims are easily proved using Lemma \ref{next to lowest}.
{\renewcommand\labelenumi{\roman{enumi})}
\begin{enumerate}
\item For any $1 \le i \le (g-2)/2$ or $(g+2)/2 \le i <g$ we get an element in
$\Uc_F^i$ by picking any $\EE_{i-1} \subset \EE_i \subset \EE_{i+1}$ and then
letting $\EE_{g-i}$ be determined by Lemma \ref{next to lowest}. Hence the locus
is isomorphic to $\PP^1$ and the intersection with $\Uc_1$ consists of the
points for which $\EE_i$ and $\EE_{g-i}$ are $\FF_{p^2}$-rational.

\item When $g$ is even we get an element in $\Uc_F^i$ by picking $\EE_{g/2-1}
\subset \EE_{g/2} \subset \EE_{g/2+1}$. Hence the locus is isomorphic to $\PP^1$
and the intersection with $\Uc_1$ consists of the points for which $\EE_{g/2}$ is
$\FF_{p^2}$-rational.

\item When $g$ is odd we get an element in $\Uc_F^i$ by picking $\EE_{(g-3)/2}
\subset \EE_{(g-1)/2} \subset \EE_{(g+3)/2}$ for which $\overline{\EE}_{(g-1)/2}
\subset F^*\overline{\EE}_{(g-1)/2}^\perp$, where
$\overline{\EE}_{(g-1)/2}=\EE_{(g-1)/2}/\EE_{(g-3)/2}$ and $F^*$ comes from the
$\FF_{p^2}$-rational structure on $\EE_{(g+3)/2}/\EE_{(g-3)/2}$ and the scalar
product is inducted from that on $\EE_g$, and then define $\EE_{(g+1)/2}$ by the
condition that
$\EE_{(g+1)/2}/\EE_{(g-3)/2}=F^*\overline{\EE}_{(g-1)/2}^\perp$. As all
non-degenerate unitary forms are equivalent, choosing a basis of
$\EE_{(g+3)/2}/\EE_{(g-3)/2}$ for which the form has the standard form
$\scpr{(x,y,z)}{(x,y,z)}=x^{p+1}+y^{p+1}+z^{p+1}$ and hence $\Uc_F^i$ is
isomorphic to the Fermat curve of degree $p+1$ and is hence is irreducible. The
intersection with $\Uc_1$ consists of the points for which $\EE_{(g-1)/2}$ is
$\FF_{p^2}$-rational and then $\EE_{(g+1)/2}=\EE_{(g-1)/2}^\perp$.

\item When $g$ is odd we get an element in $\Uc_F^i$ by picking $\EE_{(g+3)/2}$
fulfilling conditions dual to those of iii). Hence again $\Uc_F^i$ is irreducible
and the
intersection with $\Uc_1$ consists of the points for which $\EE_{(g+1)/2}$ is
$\FF_{p^2}$-rational and then $\EE_{(g-1)/2}=\EE_{(g+1)/2}^\perp$.
\end{enumerate}}
It follows from this description that two flags in $\Uc_1$ lie in the same
component of $\Uc$ if they are equivalent under the equivalence relation
generated by the relations that for any unitary $\FF_{p^2}$-flag $0=A_0\subset
A_1\subset\dots\subset A_{g-1}\subset A_g$ we may replace it by any flag that is
the same except for $A_i$ and $A_{g-i}$ for $i \in S$. The theorem then follows
from Lemma \ref{unitary flags connection}.
\end{proof}
\end{theorem}
\begin{lemma}\label{unitary flags connection}
i) Let $\k$ be a field and $\Fcal\ell_n$ the set of complete flags of vector
spaces in a finite dimensional vector space. The equivalence relation generated
by the operations of modifying a flag $E_{\bullet}$ by, for any $i$, replacing
$E_i$ by any $i$-dimensional subspace of $E_{i+1}$ containing $E_{i-1}$ contains
just one equivalence class.

ii) Let $\Fcal\ell_n$ the set of complete flags of vector spaces in an
$n$-dimensional $\FF_{p^2}$-vector space, selfdual with respect to a perfect
unitary pairing. An \Definition{elementary modification} of such a flag
$E_{\bullet}$ is obtained by either, for any $1 \le i \le (n-1)/2$, replacing
$E_i$ by any isotropic $i$-dimensional subspace of $E_{i+1}$ containing
$E_{i-1}$ and $E_{n-i}$ by its annihilator or, when $n$ is even, replacing
$E_{n/2}$ by any maximal totally isotropic subspace contained in $E_{n/2+1}$ and
contains $E_{n/2-1}$. Then the equivalence relation generated by all elementary
operations contains just one equivalence class.
\begin{proof}
Starting with i), we prove it by induction on $n$, the dimension of the vector
space $V$. Given two flags $E_\bullet$ and $F_\bullet$, if $E_1$ and $F_1$ are
equal we may use induction applied to $E_\bullet/E_1$ and $F_\bullet/E_1$. We
now use induction on the smallest $j$ such that $E_1 \subseteq F_j$, $j=1$
already being taken care of. We now get new flag $F'_\bullet$ by replacing
$F_{j-1}$ by $F_{j-2}\bigoplus E_1$, which works as $E_1 \subsetneq F_{j-1}$ and
we then have $E_1 \subseteq F'_{j-1}$.

Continuing with ii) we again use induction on $n$ and start with two selfdual
flags $E_\bullet$ and $F_\bullet$. Let us first assume that $n$ is even,
$n=2k$. Then $E_k$ and $F_k$ are isotropic subspaces. If they have non-trivial
intersection then we may pick a $1$-dimensional subspace contained in it and
then use i) to replace $E_\bullet$ and $F_\bullet$ by flags for which $E_k$ and
$F_k$ are the same and $E_1=F_1$. This implies that also $E_{n-1}=F_{n-1}$ and
we may consider $E_{n-1}/E_1$ with its two flags induced from $E_\bullet$ and
$F_\bullet$ and use induction to conclude. Assuming $F_k\cap E_k=\{0\}$ we may
again use i) to modify $F_\bullet$, keeping $F_k$ fixed, so that $E_1 \subseteq
F_{k+1}$. This means that $E_1\bigoplus F_{k-1}$ is totally isotropic and we may
replace $F_{n/2}$ by it to obtain a new flag $F'_\bullet$ for which $F'_{n/2}$
and $E_{n/2}$ intersect non-trivially.

When $n$ is odd, $n=2k+1$, we may again use induction on $n$ to finish if $E_k$
and $F_k$ intersect non-trivially. If not we may again use i) to reduce to the
case when $E_1 \subseteq F_{k+2}$ and then we may replace $F_k$ by $E_1\bigoplus
F_{k-1}$ and $F_{k+1}$ by its annihilator.
\end{proof}
\end{lemma}
\end{section}
\begin{section}{Local structure of strata}
\label{sec:Local structure of strata}

\begin{subsection}{Stratified Spaces} 
We now want to show that our moduli space locally looks like the space of
complete symplectic flags (in $2g$-dimensional space). This isomorphism will on
the one hand preserve the degeneration strata and on the other hand will be
compatible with all the forgetful maps to partial flag spaces (all the way down
to that of totally isotropic $g$-dimensional spaces). In order to have a
convenient way of formulating such a result we introduce the two following
notions:

By a \Definition{stratified space} we will mean a scheme together with a
collection of closed subschemes, called \Definition{strata}. A map between
stratified spaces is said to be stratified if it maps strata into strata.

If $P$ is a partially ordered set then a \Definition{diagram of spaces over $P$}
associates to each element $p$ of $P$ a scheme $X_p$ and to each relation $p >
p'$ a map $X_p \to X_{p'}$ fulfilling the condition that the composite $X_p \to
X_{p'} \to X_{p''}$ equals the map $X_p \to X_{p''}$ for any $p > p' > p''$. We
shall also similarly speak about a diagram of stratified spaces where both the
schemes and the maps are assumed to be stratified. Given a field $k$ and a
$k$-point $x$ of a diagram $X_\cdot$ we may speak of its (strict)
\Definition{Henselization} at $x$, which at each $p \in S$ is the Henselization
at $x$ of $X_p$.

For a positive integer $g$ we now consider the partially ordered set $P_g$ whose
elements are the subsets of $\{1,2,\dots,g-1\}$ and with ordering that of
inclusion. We have two diagrams of stratified spaces over this set: The first,
$\Fl_g^\cdot$, associates to the subset $S$ the flag space of a maximal totally
isotropic subspace $E$ of a symplectic $2g$-dimensional vector space and partial
flags of subspaces of $E$ whose dimensions form the set $S$. The map associated
to an inclusion $S \subset S'$ is simply the map forgetting some of the elements
of the flag. Similarly, we let $\Fcal_g^\cdot$ be the diagram that to a subset
$S$ associates the space of flags over the moduli space $\Acal_g$ of principally
polarized abelian varieties that associates to an abelian variety the space of
flags on the Hodge bundle whose dimensions form the subset $S$.

The diagram $\Fl_g^\cdot$ becomes a stratified diagram by considering the
stratifications given by the (closed) Schubert cells with respect to some fixed
complete flag. In positive characteristic $p$ the diagram $\Fcal_g^\cdot$
becomes a stratified diagram by considering the degeneracy loci given by the
relative positions of the Hodge flag $\EE_\cdot$ and the conjugate flag $\DD$.
\end{subsection}
\begin{subsection}{Height $1$-Maps}
For schemes in a fixed positive characteristic $p$ we shall say that a closed
immersion $S \hookrightarrow S'$ defined by the ideal sheaf $\Ical$ is a
\Definition{height $1$-map} if $\Ical_{S}^{(p)}=0$, where for an ideal $I$, we
let $I^{(p)}$ be the ideal generated by the $p$'th powers of elements of $I$. If
$R$ is a local ring in characteristic $p$, the \Definition{height $1$-hull} of
$R$ is the quotient ring $R/m_R^{(p)}$ (so that the spectrum of it is the
largest closed subscheme of $\Spec R$ for which the map from $\Spec R/m_R$ to
$\Spec R/m_R^{(p)}$ is a height $1$-map). If $k$ is a field of characteristic
$p$ and $x\co \Spec k \to S$ a $k$-map to a $k$-scheme $S$ of characteristic
$p$, then by the \Definition{height $1$-neighbourhood} of $x$ we will mean the
spectrum of the height $1$-hull of the local ring of $S$ at $x$. It is clear
that taking height $1$-neighbourhoods of $k$-points is functorial under maps
between $k$-pointed $k$-schemes. Finally, we shall say that two local rings are
\Definition{height $1$-isomorphic} if their respective height $1$-hulls are
isomorphic and that the height $1$-hull of a $k$-point is \Definition{height
$1$-smooth} if its is isomorphic to the height $1$-hull of a smooth $k$-point
(i.e., is of the form $k\bigl[[t_1,\dots,t_n]\bigr]/m^{(p)}$).
\begin{theorem}
For each perfect field $k$ of positive characteristic $p$ and each $k$-point $x$
of $\Fcal_g^\cdot$ there is a $k$-point $y$ of $\Fl_g^\cdot$ such that the
height $1$-neighbourhood of $x$ is isomorphic to the height $1$-neighbourhood of
$y$ by a stratified isomorphism of diagrams.
\begin{proof}
Denote also by $x$ the point of $\Fcal_g$, the space of complete flags of the
Hodge bundle, associated to $x$ as a point of the diagram $\Fcal_g^\cdot$. Let
$X^\cdot$ be the height $1$-neighbourhood of $x$ in $\Fcal_g^\cdot$ and $X$ the
height $1$-neighbourhood of $x$ in $\Fcal_g$. Now the ideal of the closed point
of $x$ in $X$ has a divided power structure for which all the divided powers of
order $\ge p$ are zero. This allows us get a trivialization of the restriction
of the de Rham cohomology $\HH_{|X} \riso X\times W$ which is horizontal (i.e.,
compatible with the Gauss-Manin connection on the left and the trivial
connection on the right). Now, as the absolute Frobenius map on $X$ factors
through the closed point we get that the pullback $\EE_\bullet^{(p)}$ is a
horizontal flag and then so is $\DD_\bullet$, its elements being either inverse
images of horizontal subbundles by the horizontal map $V$ or duals of horizontal
subbundles. We now get a map from $X$ to the space $\Fl_g$ of complete
symplectic flags on $W$ such that the pullback of the universal flag equals
$\EE_\cdot$. We may furthermore, choose a symplectic isomorphism of $W$ and the
standard symplectic space such that $\DD_\bullet$ is taken to the fixed complete
flag. We can extend this map in a compatible fashion for all partial flag spaces
giving a map from the diagram $X_\cdot$ to $\Fl_g^\cdot$ and we will denote by
$y$ the $k$-point which is the composite of $x$ and this map. This map is
clearly a stratified map and by the infinitesimal Torelli theorem it induces an
isomorphism from $X$ to $Y_\cdot$, the first height $1$-neighbourhood of $y$ in
$\Fl_g^\cdot$.
\end{proof}
\end{theorem}

\begin{corollary}
For each perfect field $k$ of positive characteristic $p$ and each $k$-point $x$
of $\Fcal_g$ there is a $k$-point $y$ of $\Fl_g$ such that the Henselization of
$x$ is isomorphic to the Henselization of $y$ by a stratified isomorphism.
\begin{proof}
The theorem provides such an isomorphism over the height $1$-hull $X$ of $x$.
Now, over $\Ocal_{\Fcal_g,x}$ we may extend the trivialization of $\HH_{|X}$ to
a trivialization of $\HH_{\Fcal_g,x}$ which also extends the trivialization of
$\DD$ (making, of course, no requirements of horizontality). This gives a map
from the localization, $\tilde X$, of $\Fcal_g$ at $x$ to $\Fl_g$ that extends
the map from $X$ to $\Fl_g$. It thus induces a map from $\tilde X$ to $\tilde
Y$, the localization $\Fl_g^\cdot$ at $y$. Now, this map induces an isomorphism
on tangent spaces and $\Fcal_g$ is smooth. This implies that we get an induced
isomorphism on Henselizations and proves the theorem.
\end{proof}
\end{corollary}


\begin{lemma}\label{Nothing below canonical}
Let $A$ be a principally polarised abelian variety over an algebraically closed
field. If a flag $\DD_\bullet$ for it has type $w'$ which $\le$ its canonical
type then $w'$ is the final element corresponding to the canonical type of $A$.
\begin{proof}
The flag $\DD_\bullet$ has the property, as it is of a type $\le$ to the
canonical type, that $F$ maps $\DD_i^{(p)}$ into $\DD_{\nu_w(i)}$. Consider now
the set $I$ of $i$'s for which $\DD_i$ is a member of the canonical flag. It
clearly contains $0$ and is closed under $i\mapsto
\overline\imath$. Furthermore, if $i \in I$, then $F(\DD_i^{(p)})$ has dimension
$\nu_w(i)$ but is then equal to $\DD_{\nu_w(i)}$ as it is contained in it. Hence
$I$ fulfills the conditions of Corollary \ref{r-fulfillment} and hence contains
the canonical domain which means that $\DD_\bullet$ is a refinement of the
canonical flag and hence $\nu$, the final type of $A$, and $\nu_{w'}$ coincide
on the canonical domain of $\nu$ and hence are equal by Corollary
\ref{r-fulfillment} which means that $w'$ is the final element of the canonical type.
\end{proof}
\end{lemma}
\begin{corollary}\label{dim and smooth}
i) Each stratum $\Uo_w$ is smooth of dimension $\ell(w)$.

ii) The closed stratum $\Uc_w$ is Cohen-Macaulay, reduced, and normal of
dimension $\ell(w)$ and $\Uc_w$ is the closure of $\Uo_w$ for all $w \in W_g$.

iii) If $w$ is final then the restriction of the projection $\Fcal_g \to
\Acal_g$ to $\Uo_w$ is a finite surjective \'etale covering from $\Uo_w$ to
$\Vo_w$ of degree $\gamma(w)$.
\begin{proof}
We know that each open Schubert cell of $\Fl_g$ is smooth and each closed is
Cohen-Macaulay by a proof completely along the lines of \cite[Thm 14.3]{Ful} and
by a theorem of Chevalley (cf., \cite[Cor.\ of Prop 3]{Ch}) they are smooth in
codimension $1$ so by Serre's criterion they are normal and reduced. The same
statement for the stratification of $\Fcal_g$ then follows form the
theorem. Finishing ii), the fact that $\Uc_w=\overline{\Uo_w}$ follows more or
less formally from the rest: If $x \in \Uc_w$ then we know that the dimensions
of all $\Uo_{w'}$ with $w' < w$ that pass through are $\ell(w') < \ell(w)$ but
the dimension of $\Uc_w$ at that point is $\ell(w)$ and hence $x$ must lie in
the closure of $\Uo_w$.

As for iii), that $\Uo_w$ maps into $\Vo_w$ follows from the fact that the
restriction of a final filtration to its canonical domain is a canonical
filtration (Proposition \ref{final and canonical flag}). That the map $\Uo_w \to
\Acal_g$ is unramified follows from the same statement for Schubert cells which
is \cite[Prop. 5.1]{BGG}. We next prove that $\Uo_w \to \Vo_w$ is proper. Note
that by Proposition \ref{final and canonical flag} and by the fact that by
definition $\Vo_w$ is the image of $\Uo_w$ we get that the geometric points of
$\Vo_w$ consists of the principally polarised abelian varieties with a canonical
filtration whose canonical type corresponds to the final type of $w$. Hence for
properness we may assume that we have a principally polarised abelian variety
over a discrete valuation ring $R$ such that both its generic and special point
are of type $w$ and we suppose given a final flag over the generic point. Hence
the canonical decomposition of $\Spec R$ for the abelian variety is equal to $\Spec
R$ and we have a canonical flag over $\Spec R$. As $\Uc_w$ is proper the map to it
from the generic point of $\Spec R$ extends to a map from $\Spec R$ to $\Uc_w$ hence
giving a flag over its special point. This flag is then of a type $\le w$ and
hence by Lemma \ref{Nothing below canonical} it is equal to $w$ and the image of
$\Spec R$ lies in $\Uo_w$ which proves properness.

Now, $\Vo_w$ being by definition the schematic image of $\Uo_w$ is reduced as
$\Vo_w$ is. As $\Uo_w \to \Vo_w$ is unramified it has reduced geometric fibres
and as it is finite and $\Vo_w$ is reduced, to show that it is flat it is
enough to show that the cardinalities of the geometric fibres are the same for
all geometric points of $\Vo_w$. This however is Lemma \ref{canonical to
final}. Being finite, flat, and unramified it is \'etale. That its degree is
$\gamma(w)$ follows from Lemma \ref{canonical to final}.
\end{proof}
\end{corollary}
\begin{remark} Note that the degree of the map $\Uo_w \to \Vo_w$
is $\gamma_g(w)$. It is not difficult by looking at the proof of Lemma
\ref{canonical to final} to show that it is a covering with structure group a
product of linear and unitary groups over finite fields of characteristic $p$.
\end{remark}
\end{subsection}
\end{section}
\begin{section}{Punctual flag spaces}
\label{sec:Punctual flag spaces}

Let $M$ be the Dieudonn\'e module of a truncated Barsotti-Tate group of level
$1$ over an algebraically closed field of characteristic $p$ provided with an
alternating perfect pairing (of Dieudonn\'e modules). We let $\Fcal_M$, the
\Definition{punctual flag space} for $M$, be the scheme of self-dual admissible
complete flags in $M$ for which the middle element equals $\im(V)$.  It is
well-known that every such $M$ occurs as the Dieudonn\'e module of the kernel of
multiplication by $p$ on a principally polarised abelian variety and then
$\Fcal_M$ is the intersection of $\Uc_{\emptyset}$ and the fibre over a point of
$\Acal_g$ giving rise to $M$. Also by a result of Oort (cf., \cite{O1}) the
canonical type of it determines it (over an algebraically closed field) and
hence we shall also use the notation $\Fcal_\nu$ where $\nu$ is a final
type. For $\Gamma=(I,\Scal)$ where $I \subseteq \{1,\dots,g\}$ with $\#I$ equal
to the semi-simple rank of $M$ and $\Scal$ a complete $V$ stable flag of the
$V$-semi-simple part of $M$, we then let $\Fcal_M^\Gamma$ be the part of
$\Fcal_M\cap\Uo_I^{ss}$ (which clearly can be defined directly in terms of $M$)
for which the flag induces $\Scal$ on the $V$-semi-simple part. As the $p$-rank
of $M$ is $f$ we easily see that $\Fcal_M$ is the disjoint union of the
$\Fcal_M^{\Gamma}$ and putting
$\Fcal_M^{\Scal}:=\Fcal_M^{(\{1,\ldots,g-f+1\},\Gamma)}$ we have maps
\map{S_I}{\Fcal_M^I}{\Fcal_M^{\Scal}}. These maps are homeomorphisms by
Proposition \ref{Main shuffling}. This can be seen directly by decomposing $M$
as $M^{mul}\Dsum M^{\l\l}\Dsum M^{et}$, where $V$ is bijective on $M^{mul}$, $F$
on $M^{et}$ and $F$ and $V$ nilpotent on $M^{\l\l}$. Any element of an
admissible flag over a perfect field will decompose in the same way (as that
element is stable under $F$ by definition and $V$ by duality) and is hence
determined by its intersection by $M^{mul}$, $M^{\l\l}$, and $M^{et}$. By
self-duality the intersection of all the elements of the flag with $M^{et}$ is
determined by that with $M^{mul}$ and that part is given by an arbitrary full
flag of sub-modules of $M^{mul}$ which is our $\Scal$. That means that we may
indeed reconstitute the whole flag from $\Gamma$ and the induced flag on
$M^{\l\l}$ and that any choice of flag on $M^{\l\l}$ gives rise to a flag in
$\Fcal_M^{\Gamma}$. That means that the map $\Fcal_M^{\Gamma} \to
\Fcal_{M^{\l\l}}$ is a homeomorphism and we may for all practical purposes focus
our attention on the case when $F$ and $V$ are nilpotent on $M$ (i.e., $M$ is
\Definition{local-local}). Hence in this section \Emph{unless otherwise
mentioned the Dieudonn\'e modules considered will be local-local}. Note that the
principal interest in this section will be focused on the question of which
$\Uo_w$ have a non-empty intersection with $\Fcal_M$ and that that problem is
indeed by the above considerations immediately reduced to the local-local case.

We shall make extensive use of one way to move in each $\Fcal_M$: 

Consider a $w_\emptyset \ge w \in W_g$. Assume that we have index $1 \le i \le
g-1$ for which $r_w(g+i-1,i+1)\ge i+1$. This means that for a flag $\DD_\bullet$
in $\Uc_w$ we have that $F(\DD_{i+1})\subseteq \DD_{i-1}$ or equivalently that
$F$ is zero on $\DD_{i+1}/\DD_{i-1}$. Hence if we replace $\DD_i$ by any
$\DD_{i-1}\subset \DD \subset \DD_{i+1}$ (replacing also $\DD_{2g-i}$ to make
the flag self-dual) we shall still have an admissible flag as $V(\DD)
\subseteq \DD_{i-1}$. In order to construct the $\EE$-flag we apply $V$ to
the $\DD$-flag which gives us half of the $\EE$-flag and we complement
by taking orthogonal spaces. In the $\EE$-flag now $\EE_{g-i}$ and $\EE_{g+i}$
move. 
This construction gives a mapping from the projective line
$\PP(\EE_{g-i+1}/\EE_{g-i-1})$ to $\Fcal_M$ and we shall therefore call this
family the \Definition{simple family} of index $i$ and we shall write
\Definition{$P_{w,i}$} for this simple family . The condition
$r_w(g+i-1,i+1)\ge i+1$ is equivalent to $r_w(g-i+1,i-1)=g-i-1$ and when it
is fulfilled we shall say that $g-i$ is \Definition{movable for $w$}.
\begin{proposition}
Any two points of the local flag space $\Fcal_M^{\Gamma}$ can be connected by a
sequence of simple families.
\begin{proof}
We immediately reduce to the case when $M$ is local (in which case the statement
is about $\Fcal_M$). We are going to identify $\Fcal_M$ with the scheme of
$V$-stable flags in $\im(V)$ and we prove the statement for any Dieudonn\'e
module $N$ with $F=0$ and $V$ nilpotent. Let $E_\cdot$ and $F_\cdot$ be two
$V$-stable flags in $N$. If $E_1=F_1$ then we may consider $N/E_1$ and use
induction on the length of $N$ to conclude. If not we use induction on the
smallest $i$ such that $F_1 \subseteq E_i$ which we thus may assume to be
$>1$. We now have $F_1 \subsetneq E_{i-1}$ and hence that $F_1$ is a complement
to $E_{i-1}$ in $E_i$ so that in particular $E_i/E_{i-2}=E_{i-1}/E_{i-2}\Dsum
(F_1+E_{i-2})/E_{i-2}$ which has a consequence that $V$ is zero on $E_i/E_{i-2}$
which means that every subspace of it is stable under $V$ so that we get a
$\PP^1$-family of flags in $E_i/E_{i-2}$ in which both $E_{i-1}/E_{i-2}$ and
$(F_1+E_{i-2})/E_{i-2}$ are members so that we may move $E_{i-1}$ so that it
contains $F_1$.
\end{proof}
\end{proposition}
Recall (cf, \cite[14.3]{O1}) that one defines the partial order relation on
final types $\nu_1 \subseteq \nu_2$ (resp.\ $\nu_1 \subset \nu_2$)
 by the condition that $\Vo_{\nu_1} \subseteq
\Vc_{\nu_2}$ (resp.\  $\Vo_{\nu_1} \subsetneq
\Vc_{\nu_2}$). We shall now see that this relation can be expressed in terms of
local flag spaces. For this we let $M_\nu$ be a Dieudonn\'e module of a
principally polarised Barsotti-Tate group of level $1$ with final type $\nu$
(there is up to isomorphism only one such $M_\nu$, \cite[Thm.\ 9.4]{O1}).
\begin{theorem}\label{E-O order}
i) We have that $\nu ' \subset \nu$ precisely when there is a $w \in W_g$ such that
$w \le \nu$ and there is a flag of type $w$ in $\Fcal_{M_{\nu '}}$.

ii) If there is a flag of type $w$ in $\Fcal_{M_{\nu '}}$ then there is a $w'
\le w$ such that the intersection $\Uo_{w'}\cap\Fcal_{M_{\nu '}}$ is finite.
\begin{proof}
Consider the image in $\Acal_g$ of $\Uc_\nu$. It is a closed subset containing
$\Vo_\nu$ and hence contains $\Vc_\nu$ and in particular it meets each fibre
over a point of $\Vc_\nu$ and consequently there is a point $s$ in the
intersection of $\Uc_\nu$ and the fibre over a point $t$ of $\Vo_{\nu '}$. Now,
$s$ lies in some $\Uo_w \subseteq \Uc_\nu$ and consequently fulfills $w \le \nu$
and as $\nu \le w_\emptyset$, $s$ also lies in the local flag space of $t$ and
as has been noted this is the ``same'' as $\Fcal_{M_{\nu '}}$. The converse is
clear.

As for the second part, the proof of Lemma \ref{Raynaud trick} shows that a $w'
\le w$ which is minimal for the condition that $\Uo_{w'}\cap\Fcal_{M_{\nu '}}$
is non-empty has $\Uo_{w'}\cap\Fcal_{M_{\nu '}}$ finite.
\end{proof}
\end{theorem}
The theorem allows us to reprove a result of Oort (cf., \cite{O1}).
\begin{corollary}\label{E-O is strat} 
The E-O stratification on ${\Acal}_g$ is a stratification.
\begin{proof}
The condition in \ref{E-O order} says that $\nu ' \subset \nu$ if and only if
the closure $\Uc_{\nu}$ of $\Uo_{\nu}$ has a non-empty intersection with the
punctual flag space ${\Fcal}_{\nu '}$. The proof gives more precisely that a
given point $s$ of $\Vo_{\nu '}$ lies in $\Vc_\nu$ precisely when $\Uc_\nu$
intersects the fibre over $s$ of the map $\Fcal_g \to \Acal_g$. This condition
does not depend on point $s$ by a result of Oort on Dieudonn\'e modules (cf.,
\cite{O1}).
\end{proof}
\end{corollary} 
{}From this  theorem it is clear that the condition that $\Uo_w\cap\Fcal_M\ne\emptyset$ is
important. We shall say that an admissible $w \in W_g$ \Definition{occurs} in
$\nu$, with $\nu$ a final type, if $\Uo_w\cap\Fcal_M\ne\emptyset$ and we shall write
it symbolically as \Definition{$w \to \nu$}.
\begin{remark}
It is important to realise that \Emph{a priori} this relation depends on the
characteristic that is implicit in all of this article. Hence the notation $w
\stackrel{p}{\to} \nu$ would be more appropriate. It is our hope that the
relation will \Emph{a posteriori} turn out to be independent of $p$. If not and
if one is working with several $p$ the more precise notation will have to be used.
\end{remark}
Hence we can formulate the theorem as saying that $\nu ' \subset \nu$
precisely when there exists an admissible $w$ with $w \to \nu '$ and $w \le \nu$.
Suppose final types $\nu$ and $\nu^{\prime}$ given. 
An element $w$ of minimal length in the set of minimal elements of 
$\{ w \in W_g : \nu > w ,\, w\to \nu^{\prime}\}$ in the Bruhat-Chevalley order
has the property that $\Uo_w \cap \Fcal_{M_{\nu '}}$ has dimension $0$
for the generic point of $\Vo_{\nu'}$. Clearly, then $\ell(w)\geq \ell(\nu')$
for every $w$ as in \ref{E-O order}.
\begin{example}
Since E-O strata on ${\Acal}_g$ are defined using the projection from the flag
space, the closure of an E-O stratum on ${\mathcal A}_g$ need not be given by the
Bruhat-Chevalley order on the set of final elements and indeed it isn't. Oort gave the
first counterexample for $g=7$ based on products of abelian varieties.  We
reproduce his example and give two other ones, one for $g=5$ and one for $g=6$
that do not come from products.

i) Let $g=7$ and let $w_1=[1, 2, 4, 6, 7, 10, 12]$ and $w_2=[1, 2, 3, 7, 9,
10,11]$. Then $w_1$ and $w_2$ are final elements of $W_7$ and have length
$\ell(w_1)=8$ and $\ell(w_2)=9$. In the Bruhat-Chevalley order neither $w_1 \leq w_2$ nor
$w_2 \leq w_1$ holds.  Despite this we have $\Vc_{w_1} \subset \Vc_{w_2}$. The
explanation for this lies in the fact that the simple family $P_{w_1,4}$ hits
the stratum $U_{w_3}$ with $w_3$ the element $[1, 2, 3, 7, 6,10, 11]=s_3
w_1s_4$, with $w_2 > w_3$ and $w_3 \to w_1$, so by \ref{E-O order} it follows
that $\Vc_{w_1} \subset \Vc_{w_2}$. (That there is such a simple family can be
proved directly but for now we leave it is an unsupported claim as a proof ``by
hand'' would be somewhat messy. A more systematic study of these phenomena will
appear in a subsequent paper.) This explains the phenomenon observed in
\cite{O1}, p.\ 406 (but note the misprints there). Also the element $w_2 >
w_4=[1,2,3,7,9,5,11]\to w_1$ will work for $w_1$. The element $w_1$ is the final
element corresponding to taking the product of a Dieudonn\'e module with final
element $[135]$ and a Dieudonn\'e module with final element $[1246]$ whereas
similarly $w_2$ appears as the ``product'' of the final elements $[135]$ and
$[1256]$. As $[1246] < [1256]$ there is a degeneration of a Dieudonn\'e module
of type $[1256]$ to one of type $[1246]$. This shows that this example simply
expresses the fact that $\subset$ must be stable under products whereas the
Bruhat-Chevalley order isn't. (We'd like to thank Ben Moonen for pointing this out to us.)

ii) For $g=5$ we consider the final elements $w_1= [1, 3, 4, 6, 9]$ and $w_2=[1,
2, 6, 7, 8]$ of length $5$ and $6$ and the non-final element $w_3=[1,2,6,4,8]$
in $W_5$. Then $w_3<w_2$ and $w_3 \to w_1$, so that $\Vo_{w_1}$ lies in the
closure of $\Vo_{w_2}$. But in the Bruhat-Chevalley order neither $w_1<w_2$ nor $w_2 <
w_1$ holds.

iii) Let $g=6$ and consider the final elements $w_1=[1, 3, 5, 6, 9, 11]$ and
$w_2:=[1, 2, 6, 8, 9, 10]$ of length $\ell(w_1)=8$ and $\ell(w_2)=9$. In the
Bruhat-Chevalley order we do not have $w_1 \leq w_2$. Nevertheless, $\Vo_{w_1}$ occurs in
the closure of the E-O stratum $\Vo_{w_2}$. Indeed, the admissible element
$w_3=[1, 2, 6, 8, 4, 10]$ satisfies $w_2 \geq w_3 \to w_1$: $\Uo_{w_3}$ has a
non-empty intersection with the punctual flag space ${\Fcal}_{w_1}$.  This time
neither of the elements $w_1$ and $w_2$ are products in the sense of
i). Furthermore, as $\Vo_{w_1}$ is of codimension $1$ in $\Vo_{w_2}$ this
example can not derived by taking the transitive closure of the closure under
products of the Bruhat-Chevalley relation.  The claim that we have $w_3 \to w_1$ and the
two preceding ones will be substantiated in a subsequent paper.
\end{example}
There is an approach to the study of the relation of the E-O strata and the
strata on $\Fcal_g$ which is in some sense ``dual'' to the study of punctual
flag spaces; that of considering the image in $\Acal_g$ of the $\Uo_w$. The
following result gives a compatibility result on these images and the E-O
stratification.
\begin{proposition}\label{strata fibrations}
i) The image of any $\Uo_w$, $w \in W_g$, is a union of strata $\Vo_\nu$. In
particular the image of a $\Uc_w$ is equal to some $\Vc_\nu$.

ii) For any final $\nu$ and $w \in W_g$, the maps $\Uo_w\cap\pi^{-1}\Vo_\nu \to
\Vo_\nu$ and $\Uc_w\cap\pi^{-1}\Vo_\nu \to \Vo_\nu$, where $\pi$ is the
projection $\Fcal_g \to \Acal_g$, have the property that there is a surjective
flat map $X \to \Vo_\nu$ such that the pullback of them to $X$ is isomorphic to
the product $X\times (\Fcal_\nu\cap\Uo_w)$ resp.\ $X\times
(\Fcal_\nu\cap\Uc_w)$.

iii) A generic point of a component of $\Uo_w$ maps to the generic point of some
$\Vo_\nu$ and that $\nu$ is independent of the chosen component of $\Uo_w$.
\begin{proof}
The first part follows directly from Oort's result on the uniqueness for the
Dieudonn\'e module in a stratum $\Vo_\nu$ as it implies that if one fibre of
$\pi^{-1}(\Vo_\nu) \to \Vo_\nu$ meets $\Uo_w$ then they all do. As for the
second part it would follow if we could prove that there is a surjective flat
map $X \to \Vo_\nu$ such that the pullback of $(\HH,\EE,F,V,\scpr--)$ is
isomorphic to a the constant data (provided by the Dieudonn\'e module of type
$\nu$). For this we first pass to the space $X_\nu$ of bases of $\HH$ for which
the first $g$ elements form a basis of $\EE$ which is flat surjective over
$\Vo_\nu$. Over $X_\nu$ the data is the pullback from a universal situation
where $F$, $V$, and $\scpr--$ are given by matrices. In this universal situation
we have an action of the group $G$ of base changes and two points over an
algebraically closed field give rise to isomorphic $(\HH,\EE,F,V,\scpr--)$
precisely when they are in the same orbit. By assumption (and Oort's theorem)
the image of $X_\nu$ lies in an orbit so it is enough to show that the data over
an orbit can be made constant by a flat surjective map. However, the map from
$G$ to the orbit obtained by letting $g$ act on a fixed point of the orbit has
this property.

The third part follows directly from the second.
\end{proof}
\end{proposition}
The proposition gives us a map \Definition[$\tau_p$]{}\map{\tau_p}{W_g}{W_g/S_g} which
to $w$ associates the final type of the open stratum into which each generic
point of $\Uo_w$ maps. We shall return to this map in Section
\ref{sec:Tautological rings}.
\begin{example}
Note that the punctual flag space is in general rather easy to understand as it
only depends on the image of $V$ and we are almost talking about the space of
flags stable under a nilpotent endomorphism (remember that we have reduced to
the local-local case). Almost but not quite as the endomorphism is semi-linear
rather than linear. What is complicated is the induced stratification. Already
the case of $\nu=s_3 \in W_3$ is an illustrative example. We have then that
$\ker V\cap \im V$ is of dimension $2$ (we have on Jordan block for $V$ on $\im
V$ is size $2$ and one of size $1$). The first element, $\EE_1$, of the flag
must lie in $\ker V\cap \im V$ so we get a $\PP^1$ of possibilities for it. If
$\EE_1=\im V^2$, then $V$ is zero on $\EE_3/\EE_1$ and we can choose
$\EE_2/\EE_1$ as an arbitrary subspace of $\EE_3$ giving us a $\PP^1$ of choices
for $\EE_2$. On the other, if $\EE_1\ne \im V^2$, then $\EE_3/\EE_1$ has a
Jordan block of size $2$ and hence there is only one $V$-stable $1$-dimensional
subspace and thus the flag is determined by $\EE_1$. The conclusion is that the
punctual flag space is the union of two $\PP^1$'s meeting at a single point. The
intersection point is the canonical filtration (which is a full flag) and one
can show that the rest of the points on one component are flags of type $[241]$
and the rest of the points on the other are flags of type $[315]$.
\end{example}
\end{section}
\begin{section}{Pieri formulas}
\label{sec:Pieri formulas}

In this section we are going to apply a theorem of Pittie and Ram (\cite{PR99})
to obtain a Pieri type formula for our strata. (It seems to be historically more
correct to speak of Pieri-Chevalley type formulas, cf,\cite{Ch}.) The main
application of it will not be to obtain cycle class formulas, as Pieri formulas
usually do it will not give formulas for individual strata but only for certain
linear combinations. To us the principal use of these formulas will be that it
shows that a certain strictly positive linear combination of the boundary
components will be a section of an ample line bundle (or close to ample as one
of the contributors to ampleness will be $\lambda_1$ which is ample only on the
Satake compactification). This will have as consequence affineness for the open
strata as well as a connectivity result for the boundary of the closed ones. We
shall see in Section \ref{sec:Tautological rings} that there is also a Pieri
formula for the classes of the E-O strata though we know very little about it.

In this section we are going to work with level structures. There are two
reasons for this. The first one is that we are going to exploit the ampleness
$\lambda_1$ and even formulating the notion of ampleness for a Deligne-Mumford
stack is somewhat awkward. The second is that one of the consequences of our
considerations will be an irreducibility criterion for strata. Irreducibility
for a stratum on $\Acal_g$ does not imply irreducibility for the same stratum on
the space $\Acal_{g,n}$ of principally polarised abelian varieties with level
$n$-structure, where always $p\dividesnot n$. In fact irreducibility for the
level $n$ case means irreducibility on $\Acal_g$ together with the fact that the
monodromy group of the level $n$ cover is the maximum possible. Hence in this
section we shall use $\Acal_{g,n}$ but also some toroidal compactification
$\tilde\Acal_{g,n}$ (cf., \cite{F-C}). Everything we have said so far applies to
this situation giving us in particular $\Fcal_{g,n}$ and $\tilde\Fcal_{g,n}$ but
we have the extra property that for $n \ge 3$ then $\tilde\Acal_{g,n}$ and hence
$\tilde\Fcal_{g,n}$ are smooth projective varieties.

We now introduce the classes $\l_i:=c_1(\EE_{\{i\}})$ for $1\le i \le 2g$ in the
Chow ring ${\CH}^*({\Fcal}_g)$. By self-duality of the flag $\EE_\bullet$ we
have that $\l_{2g+1-i}=-\l_i$ and by construction
$c_1(\DD_{\{i\}})=p\l_{i-g}=-p\l_{3g+1-i}$ for $g+1 \le i \le 2g$. Furthermore,
$\l_1+\cdots+\l_g$ is the pullback from $\tilde{\Acal}_g$ of $\lambda_1$, the
first Chern class of the Hodge bundle.

Now we let $M_i:=c_1(\DD_{]2g-i,2g]})$, $1\le i \le g$, and start by noting that
if $n=(n_1,\dots,n_{g})$ then $n\cdot M:=\sum_in_iM_i$ is relatively ample for
$\tilde{\Fcal}_g \to \tilde{\Acal}_g$ if $n_i>0$ for $1\le i<g$. Indeed, by
construction $L_i:=\l_{2g}+\cdots+\l_{2g-i+1}$, $1\le i \le g$ is the pullback
from the partial flag space $\tilde{\Fcal}_g[i]$ of flags with elements of rank
$i$, $2g-i$, and $g$ (and with the rank $g$-component equal to $\EE_g$) and on
$\tilde{\Fcal}_g[i]$ we have that $\l_{2g}+\cdots+\l_{2g-i+1}$ is ample. It is
then well-known that any strictly positive linear combination of these elements
is relatively ample. From the formulas above we get that
$M_i=p(L_{g-i}+\lambda_1)$ (where we put $L_0=0$). On the other hand $\lambda_1$
is almost ample; it is the pullback from $\Acal_g^*$ of an ample line bundle.

Now we identify the $L_i$ with the fundamental weights of the root system of
$C_g$. Note that $W_g$ acts on the $\l_i$ considered as parts of the weight
lattice by $\sigma(\l_i)=\l_{\sigma(i)}$ (keeping in mind that
$\l_{2g+1-i}=-\l_i$) and then acts accordingly on the $L_i$. Let us also note
that (by Chevalley's characterisation of the Bruhat-Chevalley order) if $w' < w$ with
$\l(w')=\l(w)-1$ and if $w=s_{i_1}\cdots s_{i_k}$, then $w'$ is of the form
$s_{i_1}\cdots \widehat{s_{i_r}}\cdots s_{i_k}$ which can be rewritten as
$ws_\alpha$, where $s_\alpha=(s_{i_{r+1}}\cdots
s_{i_k})^{-1}s_{i_r}(s_{i_{r+1}}\cdots s_{i_k})$ which thus is the reflection
with respect to a unique positive root.
\begin{theorem}\label{Pieri formula}
For each $1 \le i \le g$ and $w \in W_g$ we have that
\begin{displaymath}
(p\lambda_1+pL_{g-i}-wL_i)[\Uc_w]=\sum_{w'\prec w}c^i_{w,w'}[\Uc_{w'}] \in
\CH^1_{\QQ}(\Uc_w),
\end{displaymath}
where $c^i_{w,w'}\ge 0$ and $w'\prec w$ means $w'\le w$ and
$\l(w)=\l(w')+1$. Furthermore, $c^i_{w,w'} > 0$ precisely when $w'=ws_\alpha$
for $\alpha$ a positive root for which the simple root $\alpha_i$ appears with a
strictly positive coefficient when $\alpha$ is written as a linear combination
of the simple roots.
\begin{proof}
We shall use \cite{PR99} which has the following setup: We fix a semi-simple
algebraic group $G$ (which in our case is the symplectic group $\Sp_{2g}$ but
using this in the notation will only confuse) with Borel group $B$ and fix a
principal $B$-bundle $E \to X$ over an algebraic variety $X$. Letting $E(G/B)
\to X$ be the associated $G/B$-bundle we have, because its structure group is
$B$ and not just $G$, Schubert varieties $\Omega_w \to X$ (which fibre by fibre
are the usual Schubert varieties). For every weight $\lambda \in P$, $P$ being
the group of weights for $G$, we have two line bundles on $E(G/B)$; on the one
hand $y^\lambda$ obtained by regarding $\lambda$ as a character of $B$ which
gives a $G$-equivariant line bundle on $G/B$ and hence a line bundle on
$E(G/B)$, on the other hand the character $\lambda$ can also be used to
construct, with the aid of the principal $B$-bundle $E$, a line bundle
$x^\lambda$ on $X$ and then by pullback to $E(G/B)$ a line bundle also denoted
by $x^\lambda$. A result of \cite[Corollary]{PR99} then says that if $\lambda$
is a dominant weight then
\begin{displaymath}
y^\lambda[\Ocal_{\Omega_w}]=\sum_{\eta \in
\Tcal_w^\lambda}x^{\eta(1)}[\Ocal_{\Omega_{v(\eta,w)}}] \in K_0(E(G/B)).
\end{displaymath}
Here $\Tcal_w^\lambda$ is a certain set of piecewise linear paths
\map{\eta}{[0,1]}{P\Tensor\RR} in the real vector space spanned by $P$ and
$v(\eta,w)$ is a certain element in the Weyl group of $G$ which is always $\le
w$ and $\Tcal_w^\lambda$ has the property that $\eta(1) \in P$ for all its
elements $\eta$. An important property of $\Tcal_w^\lambda$ is that it only
depends on $w$ and $\lambda$ and not on $E$. It follows immediately from the
description of \cite{PR99}, that $v(\eta,w)=w$ in only one case, namely when
$\eta$ is the straight line $\eta(t)=tw\lambda$. Hence we can rewrite the
formula as
\begin{displaymath}
(y^\lambda-x^{w\lambda})[\Ocal_{\Omega_w}]={\sum_{\eta \in \Tcal_w^\lambda}}'x^{\eta(1)}[\Ocal_{\Omega_{v(\eta,w)}}],
\end{displaymath}
where the sum now runs over all elements of $\Tcal_w^\lambda$ for which
$v(\eta,w)<w$. Taking Chern characters and looking at the top term which appears
in codimension $\codim(w)+1$ we get
\begin{displaymath}
(c_1(y^\lambda)-c_1(x^{w\lambda}))[\Omega_w]={\sum_{\eta \in \Tcal_w^\lambda}}''[\Omega_{v(\eta,w)}],
\end{displaymath}
where the sum is now over the elements of $\Tcal_w^\lambda$ for which
$\l(v(\eta,w))=\l(w)-1$. To determine the multiplicity with which a given
$[\Omega_{w'}]$ appears in the right hand side we could no doubt use the
definition of $\Tcal_w^\lambda$. However, it seems easier to note that that
multiplicity is independent of $E$ and hence we may assume that $X$ is a point
and by additivity in $\lambda$ that $\lambda$ is a fundamental weight
$\lambda_i$. In that case one can use a result of Chevalley (\cite[Prop.\
10]{Ch}) to get the description of the theorem. However, we want this formula to
be true not in the Chow group of $E(G/B)$ but instead in the Chow group of the
relative Schubert subvariety of index $w$ of $E(G/B)$. This however is no
problem as the (relative) cell decomposition shows that this Chow group injects
into the Chow group of $E(G/B)$.

In our case we now use $X=\Fcal_g$ and $\EE_\bullet$ as $B$-bundle and then pull
back this formula along the section provided by $\DD_\bullet$ (over $\Uc_w$ so
that the section takes values in the relative Schubert variety of index $w$). To
make the pullback possible (note that the relative Schubert variety will in
general not be smooth over the base) we remove the relative Schubert varieties
of codimension $2$ in the relative Schubert variety in question. This forces us
to remove the part of $\Uc_w$ where the section encounters the removed
locus. This is however a codimension $2$ subset by Corollary \ref{dim and
smooth} so its removal will not affect $\CH^1_{\QQ}(\Uc_w)$. There is a small
problem however. The result of \cite{PR99} is formulated under the assumption of
a smooth variety $X$ over the complex numbers. In our case the base is an
algebraic stack in positive characteristic. We do not see any problems in
extending their proof to our situation; the assumption of a spatial base seems
irrelevant and characteristic zero is used for a vanishing theorem (\cite[Proof
of Proposition]{PR99}) which seems to hold in all characteristics for the
situation at hand (by for instance Frobenius splitting). However, we do not need
to verify this. Instead we can use that we can first introduce a level structure
to make the base spatial and then pull back the vector bundle from a
Grassmannian (and hence the flag from a flag space over such a
Grassmannian). This case can be lifted to characteristic zero and then the
result is obtained from the characteristic zero result by using the
specialisation map in $K$-theory. (A somewhat more sophisticated reduction would
no doubt prove the result for any $X$, $E$, and $G$.)
\end{proof}
\end{theorem}
To apply the theorem we start by some preliminary results which will be used to
exploit the positivity of the involved line bundles.
\begin{lemma}\label{connectedness}
Let $X$ be a proper (irreducible) variety of dimension $>1$ and $\Lcal$ a line
bundle on $X$ which is ample on some open subset $U \subseteq X$. Let
$D:=X\setminus U$ and $H\subset U$ the zero set of a section of $\Lcal_{|U}$. If
$D$ is connected then so is $D\cup H$.
\begin{proof}
By replacing the section by a power of it, we may assume that $\Lcal$ is very
ample giving an embedding $U \hookrightarrow \PP^n$. Let $Z$ be the closure of
the graph of this map in $X\times\PP^n$ and let $Y$ the image of $Z$ under the
projection on the second factor giving us two surjective maps $X \leftarrow Z
\to Y$ and let $D'$ be the inverse image in $Z$ of $D$.  Assume that $D \cup H$
is the disjoint union of the non-empty closed subsets $A$ and $B$ and let $A'$
and $B'$ be their inverse images in $Z$. Now, $Y$ is irreducible of dimension
$>1$ and hence $H''$ is connected, where $H''$ is the hyperplane section of $Y$
corresponding to $H$, so that the images of $A'$ and $B'$ in $Y$ must
meet. However, outside of $D'$ the map $Z \to Y$ is a bijection and hence the
meeting point must lie below a point of $D'$ and hence $A'$ and $B'$ both meet
$D'$. This implies that $A$ and $B$ both meet $D$ which is a contradiction as
$D$ is assumed to be connected.
\end{proof}
\end{lemma}
\begin{proposition}\label{Positivity on flag space}
Let $\Lcal$ be the determinant $\det \EE$ of the Hodge bundle over $\Acal_{g,n}$,
$n \ge 3$ (and prime to $p$).

i) There is, for each $1 \le i < g$, an integer $m_i$ such that the global
sections of $\Lambda^{g-i}(\HH/\EE)\Tensor\Lcal^{\tensor m}$ generate this bundle
over $\Acal_{g,n}$ whenever $m \ge m_i$. These $m_i$ can be chosen independently of
$p$ (but depending on $g$ and $n$).

ii) Putting $N_i:=L_i+n_i\lambda_1$ for $1 \le i < g$ and $N_g:=\lambda_1$, then
$\sum_im_iN_i$ is ample on $\Fcal_g$ if $m_i>0$ for all $1 \le i \le g$.

iii) Fix $w \in W_g$ and put $L:=\sum_{i<g}L_i$, $N:=\sum_{i<g}N_i$, and
$m=\sum_im_i$. Choose $r$, $s$, $t$, and $u$ such that $rN+t\lambda_1-wL$ resp.\
$sN+u\lambda-wL_g$ can be written as a positive linear combination of the $N_i$
(using that $L_g=-\lambda_1$). Then if $p >r+sm$ and $(g-1)p>t+um$ we have that
$p(L+(g-1)\lambda_1)+pm\lambda_1-wL-mwL_g$ is ample on $\Fcal_{g,n}$. The
constants $r$, $s$, $t$, and $u$ can be chosen independently of $p$.
\begin{proof}
Starting with i) it follows directly from the fact that $\lambda_1$ is ample on
$\Acal_{g,n}$. The independence on $p$ follows from the existence of a
model of $\Acal_{g,n}$ that exists over $\Spec\ZZ[1/n,\zeta_n]$.

As for ii) we have that $\pi_*\Ocal(L_i)=\Lambda^{g-i}(\HH/\EE)$, $1\le i <g$,
as \map{\pi}{\Fcal_{g,n}}{\Acal_{g,n}} can be identified with the space of flags
on $\HH/\EE$ and $\Ocal(L_i)$ is $\det(\HH/\EE_{2g-i})$. By definition we then
have that $\pi_*\Ocal(N_i)$ is generated by global sections on $\Acal_{g,n}$. We
know that on the flag space $\SL_g/B$ we have that the canonical ring
$\Dsum_{\lambda}H^0(\SL_g/B,\Lcal_\lambda)$, where $\lambda$ runs over the
dominant weights and $\Lcal_\lambda$ is the corresponding line bundle, is
generated by the $H^0(\SL_g/B,\Lcal_{\lambda_i})$, $1\le i <g$, where
$\lambda_i$ is the $i$'th fundamental weight (see for instance
\cite{ramanan85::projec+schub}). Also $\Ocal(\sum_{i<g}m_iN_i)$ is relatively
very ample and we have just shown  that $\pi_*\Ocal(\sum_{i<g}m_iN_i)$ is
generated by global sections. As $\lambda_1$ is ample on $\Acal_{g,n}$ we get
that $\pi_*\Ocal(\sum_{i\le g}m_iN_i)$ is ample.

Continuing with iii) We have that $N=L+n\lambda_1$ which gives
$p(L+(g-1)\lambda_1)+pm\lambda_1-wL-mwL_g=(p-(r+sm))N+(p(g-1)-(t+um))+
(rN+t\lambda_1-wL)+m(sN+u\lambda_1-wL_g)$. We then conclude by the definitions
of $r$, $s$, $t$, and $u$ and ii).
\end{proof}
\end{proposition}
\begin{remark}
i) The constants $r$, $s$, $t$, and $u$ are quite small and easy to compute. We
know nothing about the $m_i$ but imagine that they would not be too large.

ii) It would seem that the last part would not be applicable for $g=1$ but it
can be easily modified to do so. On the other hand for $g=1$ everything is
trivial anyway.
\end{remark}
We are now ready for the first application of the Pieri formula.
\begin{proposition}\label{Connectedness and affineness}
i) There is a bound depending only on $g$ and $n$ such that if $p$ is larger
than that bound, then for an irreducible component $Z$ of some $\Uo_w \subseteq
\tilde\Fcal_{g,n}$, $w \in W_g$, the union of the complement of $Z$ in
$\overline{Z}$ and the intersection of $\overline{Z}$ and
$\tilde\Fcal_{g,n}\setminus \Fcal_{g,n}$ is connected if the intersection of
$\overline Z$ with the boundary $\tilde\Fcal_{g,n}\setminus\Fcal_{g,n}$ is
connected or empty.

ii) There is a bound depending only on $g$ and $n$ such that if $p$ is larger than
that bound, then for $w \in W_g$ of semi-simple rank $0$ we have that $\Uo_w$ is
affine.
\begin{proof}
By Proposition \ref{Positivity on flag space} there is a bound depending only on
$g$ and $n$ such that if $p$ is larger than it, then
$M:=p(L+(g-1)\lambda_1)+pm\lambda_1-wL-mwL_g$ is ample on $\Fcal_{g,n}$. Summing
up Pieri's formula (Theorem \ref{Pieri formula}) for $1\le i <g$ and $m$ times
the formula for $i=g$ we get that $M[\overline{Z}]$ is supported on
$\overline{Z}$ intersected with smaller strata. (Note that the Pieri formula is
\Emph{a priori} -- and quite likely in reality -- only true modulo torsion. We
may however simply multiply it with a highly divisible integer and that doesn't
change the support.) We then conclude by Lemma \ref{connectedness} once we have .

As for the second part we argue as in the first part and conclude that
$\Uc_w\setminus \Uo_w$ is the support of an ample divisor (as in the theorem
each component must appear as we are summing up for $1 \le i \le g$ and some
$\alpha_i$ must appear in the expansion of $\alpha$) and hence $\Uo_w$ is affine.
\end{proof}
\end{proposition}
\begin{remark}
i) The condition on the intersection with the boundary somewhat cuts down on
possible applications of the first part. In the applications of the next section
it turns out that we only need to apply it when the intersection is empty.

ii) For the second part we would like to say more generally that the image of
$\Uo_w$ is affine in $\Fcal_{g,n}^*$ for some appropriate definition of
$\Fcal_{g,n}^*$ analogous to the Satake compactification. The problem is that it
doesn't seem as if some power of $M$ would normally be generated by its global
sections so that we can not define $\Fcal_{g,n}^*$ as the image of
$\tilde\Fcal_{g,n}$.
\end{remark}
\end{section}
\begin{section}{Irreducibility properties}
\label{sec:Irreducibility}

In this section we shall prove irreducibility of a large class of strata and
also that if the characteristic is large enough and our irreducibility criterion
is not fulfilled then (with some extra conditions on the stratum) the stratum is
reducible. Our proofs show two advantages of working on the flag spaces. The
major one is that our strata are normal so that irreducibility follows from
connectedness. The connectedness of the closed E-O strata except $\Vc_1$ is
proven in \cite{O1} but the $\Vc_w$ are most definitely not locally connected
and hence that does not say very much about the irreducibility. In the converse
direction we also make use of the Pieri formula. 
\begin{definition-lemma}\label{strata connectedness}
Let $\{Z_\alpha\}$ be a \Definition{stratification} of a Deligne-Mumford stack
$X$ of finite type over a field by which we mean that the strata $Z_\alpha$ are
locally closed reduced substacks of $X$ such that the closure
$\overline{Z_\alpha}$ of a stratum is the union of strata. By the
\Definition{$k$-skeleton} of the stratification we mean the union of the strata
of dimension $\le k$ (which is a closed substack). The \Definition{boundary} of
a stratum $Z_\alpha$ is the complement of $Z_\alpha$ in its closure. Assume
furthermore that each $Z_\alpha$ is irreducible and that for $Z_\alpha$ of
dimension strictly greater than some fixed $N$ we have that its boundary is
connected (and in particular non-empty). Then the intersection of a connected
union $Z$ of closed strata $\overline{Z_\alpha}$ with the $N$-skeleton is
connected.
\begin{proof}
It is enough, by induction to prove that if we remove a stratum $Z'$ of maximal
dimension which is $>N$ from $Z$, then the result remains connected. Assume that
$Z \setminus Z'$ is the disjoint union of two closed subset $Z_1$ and $Z_2$. By
assumption the boundary of $Z'$ is connected and hence lies in $Z_1$ say. That
means that $Z'\cup Z_1$ is closed and disjoint from $Z_2$ which by the
connectedness of $Z$ implies that $Z_2$ is connected.
\end{proof}
\end{definition-lemma}
\begin{proposition}\label{Connectedness of 1-skeleton}
There is a bound depending only on $g$ and $n$ such that the following is true
if $p$ is larger than that bound:

Let $X \subseteq \tilde{\Fcal}_{g,n}$ be a connected union of irreducible
components of closed strata $\Uc_w$ (for possibly different $w$) which lie
inside of $\Fcal_{g,n}$. Then the intersection of $X$ with the $1$-skeleton is
connected.
\begin{proof}
This follows directly from Proposition \ref{Connectedness and affineness} and
Lemma \ref{strata connectedness}. (Note that for level $1$ or $2$ we may pass to
a higher level in order to apply the proposition.)
\end{proof}
\end{proposition}
We now want to interpret this Proposition (and its converse which will be true
for any $p$) in arithmetical terms. Hence we define the \Definition{$1$-skeleton
graph} of level $n$ as the following edge-colored graph: Its vertices are the
points of $\Uc_1 \subset \Fcal_{g,n}$, i.e., isomorphism classes of principally
polarised superspecial $g$-dimensional abelian varieties $A$ together with a
level $n$ structure and a complete flag
$0=\DD_0\subset\DD_1\subset\dots\subset\DD_g=H^0(A,\Omega^1_A)$ on
$H^0(A,\Omega^1_A)$ for which $\DD_{g-i}^\perp=V^{-1}\DD_i$. For each $1\le i \le
g$ we connect two vertices by an edge of color $i$ if there is an irreducible
component of $\Uc_{s_i}$ which contains them.
\begin{lemma}\label{low connectedness}
If $S \subseteq \{1,\dots,g\}$ has the property that it contains $g$ and for
every $1 \le i <g$ we have that either $i$ or $g-i$ belongs to $S$, then the
subgraph of the $1$-skeleton graph consisting of all vertices and all edges of
colours $i \in S$ is connected.
\begin{proof}
This follows from \cite[Prop.~7.3]{O1} and Theorem \ref{fibre connectedness}.
\end{proof}
\end{lemma}
For a subset $S \subseteq \{1,\dots,g\}$ the \Definition{$S$-subgraph} of the
$1$-skeleton graph is the subgraph with the same vertices and with only the
edges whose colour is in $S$. This definition allows us to formulate our
irreducibility conditions.
\begin{theorem}\label{Irreducible components upstairs}
i) Let $w \in W_g$ and let $S:=\set{1\le i \le g}{s_i \le w}$. If the $S$-subgraph
of the $1$-skeleton graph is connected then $\Uc_w \subseteq \Acal_g$ is irreducible.

ii) There is a bound depending only on $g$ such that if $p$ is larger than that
bound the following is true: If $w \in W_g$ is admissible and either final or of
semi-simple rank $0$ and if $S:=\set{1\le i \le g}{s_i \le w}$, then there is a
bijection between the irreducible components of $\Uc_w \subseteq \Acal_g$ and
the connected components of the $1$-skeleton graph.
\begin{proof}
The first part is clear as Proposition \ref{lowest degeneration} says that
each component meets $\Uc_1$ and then by the assumption on connectedness of the
$S$-subgraph  $\Uc_w$ is connected but by Corollary \ref{dim and smooth} it is
normal and hence is irreducible.

As for the second part assume first that $w$ is final but of positive
semi-simple rank. That means that its Young diagram does not contain a row of
length $g$ and hence by Lemma \ref{final factoring} and the Chevalley
characterisation of the Bruhat-Chevalley order we have that $s_i \le w$ for all and hence
the $S$-subgraph is connected by Lemma \ref{low connectedness} which make the
statement trivially true. We may therefore assume that the semi-simple rank is
$0$ and hence that $\Uc_w$ lies entirely in $\Fcal_{g,n}$. In that case the
result follows from Proposition \ref{Connectedness of 1-skeleton} and the fact
that two irreducible components of two $\Uc_{s_i}$ meet only at $\Uc_1$.
\end{proof}
\end{theorem}
Projecting down to $\Acal_g$ we get the following corollary which shows
irreducibility for many E-O strata.
\begin{theorem}\label{Irreducibility of EO-strata}
Let $w \in W_g$ be a final element whose Young diagram $Y$ does not contain all
rows of length $i$ with $\lceil (g+1)/2\rceil \le i \le g$. Then $\Vc_w$ is
irreducible and the \'etale cover $\Uo_w \to \Vo_w$ is connected.
\begin{proof}
This follow from Theorem \ref{Irreducible components upstairs} and Lemmas
\ref{final factoring} and \ref{low connectedness}.
\end{proof}
\end{theorem}
\begin{example}
For $g=2$ the locus of abelian surfaces of $p$-rank $\leq 1$ is irreducible. For
$g=3$ all E-O strata except the superspecial locus ($Y=\{1,2,3\} $) and the
Moret-Bailly locus ($Y=\{2,3\}$) are irreducible.
\end{example}
In \cite{Hara} S.\ Harashita has proved that the number of irreducible
components of an E-O stratum that is contained in the supersingular locus is
given as a class number and as a consequence that, except possibly for small
$p$, these strata are reducible. As has been proved by Oort (cf., \cite[Prop.\
5.2]{Hara}) these strata are exactly the ones to which Theorem
\ref{Irreducibility of EO-strata} does not apply.

We shall finish this section by showing that the $1$-skeleton graph can be
described in purely arithmetic terms very strongly reminiscent of the results
Harashita. Note that even for final elements our results are not formally
equivalent to Harashita's as we are dealing with the set of components of the
final strata in $\Fcal_g$ whereas Harashita is dealing with their images in
$\Acal_g$. In any case our counting of the number of components uses
\ref{Irreducible components upstairs} and hence is valid only for sufficiently
large $p$ whereas Harashita's are true unconditionally.

We start by giving a well-known description of the vertices of the $1$-skeleton
graph (see for instance \cite{E}) valid when $g>1$. We fix a super-singular
elliptic curve $E$ and its endomorphism ring $\Db$ which is provided with the
Rosati involution $*$. To simplify life we assume, as we may, that $E$ is
defined over $\FF_p$ and hence $\Db$ contains the Frobenius map $F$. It has the
property that $\Db F=F\Db=\Db f\Db$, the unique maximal ideal containing
$p$. Furthermore we have that $\Db/\Db F\iso\FF_{p^2}$. There is then a
bijection between isomorphism classes of \Definition{$\protect\Db$-lattices} $M$
(i.e., right modules torsion free and finitely generated as abelian groups) of
rank $g$ (i.e., of rank $4g$ as abelian groups) and $g$-dimensional abelian
varieties $A$. The correspondence associates to the abelian variety $A$ the
$\Db$-module $\Hom(E,A)$. Polarisations on $A$ then correspond to positive
definite \Definition{unitary forms}, i.e., a bi-additive map
\map{\scpr--}{M\times M}{\Db} such that $\scpr{md}{n}=\scpr{m}{n}dn$,
$\scpr{n}{m}=(\scpr{m}{n})^*$, and $m\ne 0 \implies \scpr{m}{m} >0$. The
polarisation is principal precisely when the form is \Definition{perfect},
i.e., the induced map of right $\Db$-modules $M \to \Hom_{\Db}(M,\Db)$ given $n
\mapsto (m \mapsto \scpr{m}{n})$ is an isomorphism. In general we put
$M^*:=\Hom_{\Db}(M,\Db)$ and then the form induces an embedding $M \to M^*$
which makes the image of finite index. More precisely we on $M\Tensor\QQ$ we get
a non-degenerate pairing with values in $\Db\Tensor\QQ$ and then we may identify
$M^*$ with the set \set{n \in M\Tensor\QQ}{\forall m\in M\co \scpr{m}{n} \in
\Db}. Using this we get a $\Db\Tensor\QQ/\ZZ$-valued unitary perfect form on
$M^*/M$ given by $\scpr{\overline m}{\overline n}:=\scpr{m}{n} \bmod \Db$. As
usual $\Db$-valued super-lattices $M \subseteq N$ corresponds to totally
isotropic submodules of $M^*/M$.

If now $S \subseteq \{0,\dots,g\}$ is stable under $i \mapsto g-i$ then an
\Definition{arithmetic $S$-flag} consists of the choice of unitary forms on $\Db$-modules
$M_i$ of rank $g$ for $i \in S$ and compatible isometric embeddings $M_i \hookrightarrow
M_j$ whenever $i < j$ fulfilling the following conditions:
\begin{itemize}
\item For all $i \in S$ with $i \ge g/2$ we have that $M_i^*/M_i$ is killed by
$F$ of $\Db$ and can hence be considered as an $\Db/m=\FF_{p^2}$-vector space
with a perfect unitary form.

\item We have that $FM_i^*=M_{g-i}$ for all $i \in S$.

\item The length of $M_j/M_i$ for $i<j$ is equal to $j-i$.
\end{itemize}
\begin{remark}
i) Note that we allow $S$ to be empty in which case there is exactly one
$S$-flag.

ii) As follows (implicitly) from the proof of next proposition, the isomorphism
class of an element of an arithmetic $S$-flag tensored with $\QQ$ is independent
of the flag. Hence we may consider only lattices in a fixed unitary form over
$\Db\Tensor\QQ$ and then think of $M_i$ as a sublattice of $M_j$.
\end{remark}
\begin{proposition}
Let $S \subseteq \{1,\dots,g\}$ and let $\overline{S} \subseteq \{0,\dots,g\}$
be the set of integers of the form $i$ or $g-i$ for $i \in S$. Then the set of
isomorphism classes of $\overline{S}$-flags is in bijection with the set of the
set of connected components of the $S$-subgraph of the $1$-skeleton graph.
\begin{proof}
This follows from the discussion above and Theorem \ref{fibre connectedness}
once we have proven that an $S$-flag can be extended to an
$\{0,\dots,g\}$-flag. Assume first that $g \notin S$ and let $i \in S$ be the
largest element in $S$. By assumption we have $FM_i^*=M_{g-i}\subseteq M_i
\subset M_i^*$ and we have that the length of $M_i^*/FM_i^*$ is $g$ whereas by
assumption that of $M_i/M_{g-i}$ is $2i-g$. Together this gives that the length
of $ M_i^*/M_i$ is $g-(2i-g)=2(g-i)$. As the form on $M_i^*/M_i$ is a
non-degenerate unitary $\FF_{p^2}$-form and as all such forms are equivalent we
get that there is a $g-i$-dimensional totally isotropic (and hence its own
orthogonal) subspace of $ M_i^*/M_i$ which then gives an $M_i\subset M_g \subset
M_i^*$ and as $M_g/M_i$ is its own orthogonal we get that the pairing on $M_g$
is perfect. We then put $M_0:=FM_g$ and the rest of extension is immediate.
\end{proof}
\end{proposition}
\end{section}
\begin{section}{The Cycle Classes}

If one wishes to exploit our stratification on ${\Fcal}_g$ and the E-O
stratification on ${\Acal}_g$ fully, one needs to know the cohomology classes
(or Chow classes) of the (closed) strata. In this section we show how to
calculate these classes. The original idea for the determination of the cycle
classes can be illustrated well by the $p$-rank strata. If $X$ is a principally
polarized abelian variety of dimension $g$ which is generic in the sense that
its $p$-rank is $g$, then its kernel of multiplication by $p$ contains a direct
sum of $g$ copies of $\mu_p$, the multiplicative group scheme of order
$p$. Viewing $\mu_p$ as a fat point at the origin we get a tangent vector to
$X$ at the origin, and after dualizing, a section of the Hodge bundle. If the
abelian variety looses $p$-rank under specialization the $g$ sections thus
obtained become dependent and the loci where this happens have classes
represented by a multiple of the Chern classes of the Hodge bundle.

To calculate the cycle classes of the E-O strata on ${\Acal}_g \otimes \FF_p$ we
shall use the theory of degeneration cycles of maps between vector bundles. To
this end we shall apply formulas of Fulton for degeneracy loci of symplectic
bundle maps to calculate the classes of the $\Uc_{w}$ and formulas of Pragacz
and Ratajski and of Kresch and Tamvakis for calculating those of $\Vc_{\nu}$.

\begin{subsection}{Fulton's formulas}

Over the flag space ${\Fcal}_g$ we have the pull back of the de Rham bundle and
the two flags $\EE_{\bullet}$ and $\DD_{\bullet}$ on it. We denote by $\ell_i$
the roots of the Chern classes of $\EE$ so that $c_1(\EE_i)=
\ell_1+\ell_2+\ldots +\ell_i$. We then have $c_1(\DD_{g+i})- c_1(\DD_{g+i-1})= p
\ell_i$.

Recall that for each element $w \in W_g$ we have a degeneracy locus ${\Uc}_w$ in
${\Fcal}_g$ and $\tilde{\Fcal}_g$. Its codimension equals the length $\ell(w)$
and it thus makes sense to consider the cycle class $u_w = [{\Uc}_w]$ in
$\CH_{\QQ}^{\rm codim (w)} (\tilde{\Fcal}_g)$, where we write ${\Fcal}_g$ instead
of ${\Fcal}_g \otimes \FF_p$.

Fulton's setup in \cite{Fu2} is the following (or more precisely the part that
interests us): We have a symplectic vector bundle $H$ over some scheme $X$ and
two full symplectic flags $0\subset \dots \subset E_{2} \subset E_{1}=H$ and
$0\subset \dots \subset D_{2} \subset D_{1}=H$. For each $w \in W_g$ on defines
the degeneracy locus $\Uc_w$ by \set{x \in X}{\forall i,j\co \dim
\EE_{i,x}\cap\DD_{j,x}\le r_w(i,j)} (of course this closed subset is given a
scheme structure by considering these conditions as rank conditions for maps of
vector bundles). Fulton then defines a polynomial in two sets of variables $x_i$
and $y_j$, $i,j=1,\dots,g$, such that if this polynomial is evaluated as
$x_i=c_1(E_{i}/E_{i+1})$ and $y_j=c_1(D_{j}/D_{j+1})$ then it gives the class of
$\Uc_w$ provided that $\Uc_w$ has the expected codimension $\codim(w)$ (and $X$
is Cohen-Macaulay). The precise definition of these polynomials are as follows:
For a partition $\mu=\{ \mu_1 > \mu_2 > \ldots > \mu_r > 0\}$ with $r\leq g$ and
$\mu_1 \leq g$ one defines a Schur function
$$
\Delta_{\mu}(x):= \det (x_{\mu_i+j-i})_{1 \leq i, j \leq r}
$$
in the variables $x_i$ and puts
$$
\Delta(x,y):= \Delta_{(g,g-1,\ldots,1)}(\sigma_i(x_1,\ldots,x_g)+\sigma_i(
y_1,\ldots,y_g)),
$$
where the $\sigma_i$ is the $i$-th elementary symmetric function. One then
considers the ``divided difference operators" $\partial_i$ on the polynomial
ring $\ZZ[x_1,\ldots,x_g]$ by
$$
\partial_i(F(x))=
\begin{cases}
\frac{F(x)-F(s_ix)}{x_i-x_{i+1}} & \text{if $i <g$},\\
\frac{F(x)-F(s_g^{\prime}x)}{2x_g} & \text{if $i=g$},
\end{cases}
$$
where $s_i$ interchanges $x_i$ and $x_{i+1}$ for $i=1,\ldots,g-1$, but
$s_g^{\prime}$ sends $x_g$ to $-x_g$ and leaves the other $x_i$ unchanged.  We
write an element $w \in W_g$ as a product $w=s_{i_{\ell}}s_{i_{\ell -1}} \ldots
s_{i_1}$ with $\ell=\ell(w)$ and set
\begin{equation}\label{Fulton polynomial}
P_w:= \partial_{i_{\ell}} \cdots \partial_{i_1}( \prod_{i+j\leq g} (x_i-y_j)
\cdot \Delta).
\end{equation}
An application of Fulton's formulas gives the following.
\begin{theorem}\label{Fulton} Let $w=s_{i_{\ell}}s_{i_{\ell -1}}
\ldots s_{i_1}$ with $\ell=\ell(w)$ be an element of the Weyl group $W_g$.  Then
the cycle class $u_w:=[\Uc_w]$ in $\CH_{\QQ}^{\codim(w)}(\tilde{\Fcal}_g)$ is
given by
$$
u_w= \partial_{i_{1}} \cdots \partial_{i_{\ell}}\left(  \prod_{i+j\leq g} (x_i-y_j)
\cdot \Delta(x,y)\right)_{|x_i=-\ell_i, y_j=p\ell_j}.
$$
\begin{proof}
By construction $\Uc_w$ is the degeneracy locus of the flags $\EE_\bullet$ and
$\DD_\bullet$. By Corollary \ref{dim and smooth} they have the expected
dimension and hence the degeneracy cycle class is equal to the class of $\Uc_w$.
\end{proof}
\end{theorem}
For a final element $w\in W_g$ the map $\Uc_w \to \Vc_w$ is generically finite
of degree $\gamma_g(w)$. By applying the Gysin map to the formula of
the theorem using formula \ref{Gysin} we can in principle calculate the
cohomology classes of all the push-downs of final strata, hence of the E-O
strata.
\begin{example}
\label{g=2example} $g=2$.

The Weyl group $W_2$ consists of $8$ elements; we give the cycle classes in
$\tilde{\Fcal}_2$ and the push-downs on $\tilde{\Acal}_2$.
$$
\begin{matrix}
&  w & s & \ell & [\Uc_w] & \pi_*([\Uc_w]) \\ 
\noalign{\hrule}\\
& [4,3] & s_1s_2s_1s_2 & 4 & 1 & 0 \\ 
&[4,2] & s_1s_2s_1 & 3&  (p-1)\lambda_1 & 0\\
&[3,4] & s_2s_1s_2 & 3 & -\ell_1+p\ell_2 & 1+p \\
&[2,4] & s_1s_2 & 2 & (1-p)\ell_1^2+(p^2-p)\lambda_2 & (p-1)\lambda_1\\
&[3,1] & s_2s_1 & 2 & (1-p^2)\ell_1^2+(1-p)\ell_1\ell_2+(1-p)\ell_2^2 &p(p-1)\lambda_1 \\
&[2,1] & s_1 & 1 & (p-1)(p^2+1)\lambda_1\lambda_2 & 0 \\
&[1,3] & s_2 & 1 &(p^2-1)\ell_1^2(\ell_1-p\ell_2)  & (p-1)(p^2-1)
\lambda_2\\
&[1,2] & 1 & 0 & -(p^4-1)\lambda_1\lambda_2 \ell_1 & (p^4-1)\lambda_1\lambda_2\\
\end{matrix}
$$
In the flag space ${\Fcal}_2$ the stratum corresponding to the empty diagram is
$U_{s_2s_1s_2}$ and the strata contained in its closure are the four final ones
$U_{s_2s_1s_2}$, $U_{s_1s_2}$, $U_{s_2}$ and $U_1$ and the two non-final ones
$U_{s_1}$ and $U_{s_2s_1}$. The Bruhat-Chevalley order on these is given by the diagram:
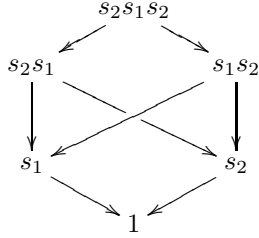
\begin{displaymath}
\begin{xy}
\xymatrix @-5mm {
&s_2s_1s_2\ar[dr]\ar[dl]\\
s_2s_1\ar'[dr][ddrr]\ar[dd]&&s_1s_2\ar[ddll]\ar[dd]\\
&&\\
s_1\ar[dr]&&s_2\ar[dl]\\
&1
}
\end{xy}
\end{displaymath}
The four final strata $U_{s_2s_1s_2}$, $U_{s_1s_2}$, $U_{s_2}$ and $U_1$ lie
\'etale of degree $1$ over the $p$-rank $2$ locus, the $p$-rank $1$ locus, the
locus of abelian surfaces with $p$-rank $0$ and $a$-number $1$, and the locus of
superspecial surfaces ($a=2$). The locus $U_{s_1}$ is an open part of the fibres
over the superspecial points. The locus $U_{s_2s_1}$ is of dimension $2$ and
lies finite but inseparably of degree $p$ over the $p$-rank $1$ locus.  Then
$E_1$ corresponds to an $\alpha_p$ and $E_2/E_1$ to a $\mu_p$. In the final type
locus $U_{s_1s_2}$ the filtration is $\mu_p \subset \mu_p \oplus \alpha_p$. Note
that this description is compatible with the calculated classes of the loci.
\end{example}
We have implemented the calculation of the Gysin map in Macaulay2 (cf., \cite{M2}) and
calculated all cycle classes for $g\leq 5$. For $g=3,4$ the reader will find the
classes in Appendix. (The Macaulay2 code for performing the calculations can be
found at \url{http://www.math.su.se/~teke/strata.m2}.) We shall return to the
qualitative consequences one can draw from Theorem \ref{Fulton} in the next section.
\end{subsection}

\begin{subsection}{The $p$-rank strata.}

It is very useful to have closed formulas for the cycle classes of important
strata. We give the formulas for the strata defined by the $p$-rank and by the
$a$-number. The formulas for the $p$-rank strata can be derived immediately from
the definition of the strata.

Let $V_f$ be the closed E-O stratum of $\tilde{\Acal}_g$ of semi-abelian
varieties of $p$-rank $\leq f$. It has codimension $g-f$. To calculate its class
we consider the element $w_{\emptyset}$, the longest final element.  The
corresponding locus ${\Uc}_{\emptyset}$ is a generically finite cover of
${\Acal}_g$ of degree $\gamma_g(w_{\emptyset})= \prod_{i=1}^{g-1}
(p^i+p^{i-1}+\ldots +1)$.  The map of $\Uo_{\emptyset}$ to the $p$-rank $g$
locus is finite.  The space $\Uc_{\emptyset}$ contains the degeneracy loci
${\Ucal}_w$ for all final elements $w \in W_g$.  The condition that a point $x$
of ${\Fcal}_g$ lies in ${\Uc}_{\emptyset}$ is that the filtration $\EE_i$ for
$i=1,\ldots,g$ is stable under $V$. By forgetting part of the flag and
considering flags $\EE_j$ with $j=i,\ldots,g$ we find that ${\Uc}_{\emptyset}
\to {\Acal}_g$ is fibered by generically finite morphisms
$$
{\Uc}_{\emptyset}={\Uc}^{(1)} {\buildrel \pi_1 \over \longrightarrow}
{\Uc}^{(2)} {\buildrel \pi_2 \over \longrightarrow} \ldots {\buildrel \pi_{g-1}
\over \longrightarrow} {\Uc}^{(g)}={\Acal}_g.
$$
We shall write $\pi_{i,j}$ for the composition $\pi_j \pi_{j-1} \cdots \pi_i
\colon {\Uc}^{(i)} \to {\Uc}^{(j)}$ and $\pi_{\emptyset}=\pi_{1,g}$.

Since $V_{g-1}$ is given by the vanishing of the map $\det(V)\co  \det(\EE_g) \to
\det(\EE^{(p)})$ the class of $V_{g-1}$ is $(p-1)\lambda_1$.  The pull-back of
$V_{g-1}$ to $\Uc_{\emptyset}$ decomposes in $g$ irreducible components
$$
\pi_{\emptyset}^{-1}(V_{g-1})= \cup_{i=1}^g Z_i,
$$
where $Z_i$ is the degeneracy locus of the induced map $\phi_i = V_{|{\Lcal}_i}
\colon {\Lcal}_i \to {\Lcal}_i^{(p)}$. Note that the $Z_i$ are the $\Uc_w$ for
the $w$ which are shuffles of the final element $u_{g-1}$ (see subsection
\ref{ssec:Some Important Strata})  defining the E-O stratum of $p$-rank $f$ and
$Z_g$ is the stratum corresponding to the element $u_{g-1}\in W_g$. An abelian
variety of $p$-rank $g-1$ and $a$-number $1$ has a unique subgroup scheme
$\alpha_p$. The index $i$ of $Z_i$ indicates where this subgroup scheme can be
found (i.e.\ its Dieudonn\'e module lies in $\EE_i$, but not in $\EE_{i-1}$).

It follows from the definition of $Z_i$ as degeneracy set that the class of
$Z_i$ on $\Uc_{\emptyset}$ equals $(p-1)\ell_i$ as $\phi_i$ can be interpreted
as a section of $ {\Lcal}_i^{(p)}\otimes {\Lcal}_i^{-1}$.  We also know by
\ref{ssec:Shuffling flags} that the map $Z_i \to Z_{i+1}$ is
inseparable. Therefore $\pi_{\emptyset}([Z_i]) = p^{n(i)}
\pi_{\emptyset}([Z_g])$ for some integer $n(i)\geq g-i$. Using the fact that
$(\pi_{\emptyset})_* ([Z_g])=\gamma_g(u_{g-1})[V_{g-1}] = \deg (\pi_{1,g-1})
[V_{g-1}]$ we see that
$$
(\pi_{\emptyset})_*(\pi_{\emptyset}^*([V_{g-1}]) = \sum_{i=1}^g
(\pi_{\emptyset})_*(Z_i) =\sum_{i=1}^g p^{n(i)} \deg (\pi_{1,g-1}) [V_{g-1}],
$$
while on the other hand
$$
(\pi_{\emptyset})_*(\pi_{\emptyset}^*([V_{g-1}])= \deg(\pi_{\emptyset})
[V_{g-1}] =(1+p+\ldots +p^{g-1}) \deg(\pi_{1,g-1}) [V_{g-1}].
$$
Comparison yields that $n(i)=g-i$ and so we find
$$
(\pi_{\emptyset})_*(\ell_i)= p^{g-i}\deg(\pi_{1,g-1}) \lambda_1
$$
and
$$
(\pi_{\emptyset})_*([Z_i])=(p-1)p^{g-i} \deg(\pi_{1,g-1}) \lambda_1 .
$$

\begin{lemma}\label{push down} In the Chow groups with rational 
coefficients of ${\Uc}^{(i)}$ and ${\Uc}^{(i+1)}$ we have for the push down of
the $j$-th Chern class $\lambda_j(i)$ of $\EE_i$ the relations:
$$
\pi^{i}_*\lambda_j(i) = p^j (p^{i-j} + p^{i-j-1} + ... + p + 1)\lambda_j(i+1)
$$
and
$$
p^{f(g-f)} \, (\pi_{1,g})_*(\ell_g\ell_{g-1} \cdots \ell_{f+1})=
(\pi_{1,g})_*(\ell_1\ell_{2} \cdots \ell_{g-f}).
$$
\begin{proof} The relation $(\pi_1)_*([Z_1])=p\, [Z_2]$ translates into
the case $j=1$ and $i=1$. Using the push-pull formula and the relations
$\pi_i^*(\lambda_j(i+1)=\ell_{i+1}\lambda_{j-1}(i)+ \lambda_j(i)$ the formulas
for the push downs of the $\lambda_j(i)$ follow by induction on $j$ and $i$.
\end{proof}
\end{lemma}

We now calculate the class of all $p$-rank strata $V_f$.
\begin{theorem}\label{p-rank classes}
The class of the locus $V_f$ of semi-abelian varieties of $p$-rank $\leq f$ in
the Chow ring $\CH_{\QQ}(\tilde{\Acal}_g)$ equals
$$
[V_f]= (p-1)(p^2-1) \cdots (p^{g-f}-1) \, \lambda_{g-f}.
$$
\begin{proof} The class of the final stratum $\Uc_{u_f}$ on $\Uc_{\emptyset}$
is given by the formula
$$
(p-1)^{g-f} \ell_g \ell_{g-1} \cdots \ell_{f+1}
$$
as it is the simultaneous degeneracy class of the maps $\phi_j$ for $j=f+1,
\ldots, g$. By pushing down this class under $\pi_{\emptyset}= \pi_{1,g}$ we
find using \ref{push down} and the notation $\lambda_j(i)= c_j(\EE_i)$ that
$$
\begin{aligned}
(\pi_{1,g})_*(\ell_g \ell_{g-1} \cdots \ell_{f+1}) & = p^{-f(g-f)}
(\pi_{1,g})_*(\ell_1\ell_{2} \cdots \ell_{g-f}) \\ &= p^{-f(g-f)} \,
(\pi_{1,g})_*(\pi_{1,g-f}^*(\lambda_{g-f}(g-f)) \\ &=p^{-f(g-f)} \,
\deg(\pi_{1,g-f}) \, (\pi_{g-f,g})_*( \lambda_{g-f}(g-f)).\\
\end{aligned}
$$
Applying Lemma \ref{push down} repeatedly we find
$$
\begin{aligned}
(\pi_{g-f,g})_*( \lambda_{g-f}(g-f))& =p^{f(g-f)}(1+p)(1+p+p^2)\cdots (1+\ldots
+p^{f-1}) \lambda_g\\ &= p^{f(g-f)} \, \gamma_g(u_f) \, \lambda_g,\\
\end{aligned}
$$
with $\gamma_g(u_f)$ the number of final filtrations refining the canonical
filtration of $u_f$.  Hence we get $(\pi_{\emptyset})_*(\Uc_{u_f})= (p-1)^{g-f}
\deg(\pi_{1,g-f}) \gamma_g(u_f) \lambda_g$.  On the other hand, we have that
$(\pi_{\emptyset})_* (\Uc_{u_f})= \gamma_g(u_f) \, [V_f]$. Together this proves
the result.
\end{proof}
\end{theorem}
\end{subsection}
\begin{subsection}{The $a$-number Strata}

Another case where we can find attractive explicit formulas is that of the E-O
strata ${\Vcal}_w$ with $w$ the element of $W_g$ associated to
$Y=\{1,2,\ldots,a\}$. We denote these by $T_a$. Here we can work directly on
${\Acal}_g$.  The locus $T_a$ on ${\Acal}_g$ may be defined as the locus $ \{ x
\in {\Acal}_g \co {\rm rank}(V)|_{\EE_g} \leq g-a \}$.  We have $T_{a+1} \subset
T_a$ and $\dim (T_g)=0$. We apply now formulas of Pragacz and Ratajski,
cf.~\cite{P-R} for the degeneracy locus for the rank of a self-adjoint bundle
map of symplectic bundles globalizing the results in isotropic Schubert calculus
from \cite{P}. Before we apply their result to our case we have to introduce
some notation.

Define for a vector bundle $A$ with Chern classes $a_i$ the expression
$$
Q_{ij}(A) := a_ia_j + 2\sum_{k=1}^j (-1)^k a_{i+k}a_{j-k} \quad { \rm for }
\quad i>j.
$$
A subset $\beta = \{g\geq \beta_1> \ldots > \beta_r\geq 0\}$ of
$\{1,2,\ldots,g\}$ (with $r$ even, note that $\beta_r$ may be zero) is called
admissible and for such subsets we set
$$
Q_{\beta} = {\rm Pfaffian} (x_{ij}),
$$
where the matrix $(x_{ij})$ is anti-symmetric with entries $ x_{ij} =
Q_{\beta_i,\beta_j} $. Applying the formula of Pragacz-Ratajski to our situation
gives the following result:

\begin{theorem} The cycle class $[T_a]$  of the reduced locus $T_a$ of abelian
varieties with $a$-number $\geq a$ is given by
$$
[T_a]=\sum_{\beta} Q_{\beta}(\EE^{(p)})\cdot Q_{\rho(a) -\beta}(\EE^*) ,
$$
where the sum is over the admissible subsets $\beta$ contained in the subset
$\rho(a) = \{a, a-1, a-2, ..., 1\}$. \end{theorem}

\noindent
\begin{example}
    \begin{align} \notag [T_1]=&\, (p-1)\lambda_1\\ \notag [T_2] =&\,
    (p-1)(p^2+1) \lambda_1 \lambda_2 - (p^3-1) 2 \lambda_3 \\ \notag [T_3]=&\,
    \\ \notag & \dots \\ \notag
    [T_g] =&\, (p-1)(p^2+1) \ldots (p^g+ (-1)^g) \lambda_1 \lambda_2 \cdots
    \lambda_g.\\ \notag \end{align}
\end{example}
As a corollary we find a result of one of us (cf.~\cite{E}) on the number of
principally polarized abelian varieties with $a=g$.

\begin{corollary}\label{Eresult} We have
$$
\sum_X \frac{1 }{ \#{\rm Aut}(X)} = (-1)^{g(g+1)/2} 2^{-g} \big[\prod_{j=1}^g
(p^j+(-1)^j)\big] \cdot \zeta(-1) \zeta(-3) \ldots \zeta(1-2g),
$$
where the sum is over the isomorphism classes (over ${\bar{\FF}}_p$) of
principally polarized abelian varieties of dimension $g$ with $a=g$, and
$\zeta(s)$ is the Riemann zeta function.
\begin{proof} Combine the formula for $T_g$ with the
Hirzebruch-Mumford Proportionality Theorem (see \cite{vdG1}) which says that
$$
\deg(\lambda_1\lambda_2\cdots \lambda_g)= (-1)^{\frac{g(g+1)}{2}} \prod_{j=1}^g
\frac{\zeta(1-2j)}{2},
$$
when interpreted for the stack ${\Acal}_g$.
\end{proof}
\end{corollary}
The formulas for the cycles classes of the $p$-rank strata and the $a$-number
strata can be seen as generalizations of the classical formula of Deuring (known
as Deuring's Mass Formula) which states that
$$
\sum_E \frac{1}{\# {\rm Aut}(E) } = \frac{p-1}{24},
$$
where the sum is over the isomorphism classes over ${\bar{\FF}}_p$ of
supersingular elliptic curves. It is obtained from the formula for $V_{g-1}$ or
$T_1$ for $g=1$, i.e., $[V_1]=(p-1)\lambda_1$, by remarking that the degree of
$\lambda_1$ is $1/12$ times the degree of a generic point of the stack
$\tilde{\Acal}_1$.

One can obtain formulas for all the E-O strata by applying the formulas of
Pragacz-Ratajski or those of Kresch-Tamvakis \cite[Cor.\
4]{kresch02::doubl+schub}.  If $Y$ is a Young diagram given by a subset
$\{\xi_1,\ldots,\xi_r\}$ we call $|\xi|=\sum_{i=1}^r \xi_i$ the weight and $r$
the length of $\xi$. Moreover, we need the \Definition{excess}
$e(\xi)=|\xi|-r(r+1)/2$ and the \Definition{intertwining number} $e(\xi,\eta)$
of two strict partitions with $\xi \cap \eta = \emptyset$ by
$$
e(\xi,\eta)= \sum_{i\geq 1} i \, \#\{ j : \xi_i>\eta_j > \xi_{i+1} \}
$$
(where we use $\xi_k=0$ if $k>r$). We put $\rho_g=\{g,g-1,\ldots,1\}$ and
$\xi^{\prime}=\rho_g - \xi$ and have then $e(\xi,\xi^{\prime})=e(\xi)$.  The
formula obtained by applying the result of Kresch and Tamvakis interpolates
between the formulas for the two special cases, the $p$-rank strata and
$a$-number strata, as follows:

\begin{theorem}\label{General E-O cycle classes}
For a Young diagram given by a partition $\xi$ we have
$$
[\Vc_{Y}]=(-1)^{e(\xi)+|\xi'|} \sum_{\alpha} Q_{\alpha}(\EE^{(p)})
\sum_{\beta} (-1)^{e(\alpha,\beta)}Q_{(\alpha \cup \beta)^{\prime}}(\EE^*)
\det(c_{\beta_i-\xi_j^{\prime}}(\EE^*_{g-\xi_j^{\prime}})),
$$
where the sum is over all admissible $\alpha$ and all admissible
$\beta$ that contain
$\xi^{\prime}$ with length $\ell(\beta)=\ell(\xi^{\prime})$ and $\alpha \cap 
\beta= \emptyset$.
\end{theorem}
\end{subsection}

\begin{subsection}{Positivity of tautological classes}

The Hodge bundle possesses certain positivity properties. It is well-known that
the determinant of the Hodge bundle (represented by the class $\lambda_1$) is
ample on ${\Acal}_g$.  Over $\CC$ this is a classical result, while in positive
characteristic this was proven by Moret-Bailly \cite{M-B}. On the other hand,
the Hodge bundle itself is not positive in positive characteristic.  For
example, for $g=2$ the restriction of $\EE$ to a line from the $p$-rank $0$
locus is $O(-1)\oplus O(p)$, \cite{M-B2}. But our results \ref{p-rank classes}
imply the following non-negativity result.

\begin{theorem}\label{positivity}
The Chern classes $\lambda_i \in \CH_{\QQ}({\Acal}_g \otimes \FF_p)$
$(i=1,\ldots,g)$ of the Hodge bundle $\EE$ are represented by effective classes
with $\QQ$-coefficients.
\end{theorem}
\end{subsection}
\end{section}
\begin{section}{Tautological rings}
\label{sec:Tautological rings}

We shall now interpret the results of previous sections in terms of tautological
rings. Recall that the tautological ring of $\tilde\Acal_g$ is the subring of
$\CH_{\QQ}(\tilde\Acal_g)$ generated by the Chern classes $\lambda$. To obtain
maximal precision we shall use the subring and \Emph{not} the $\QQ$-subalgebra
(but note that is still a subring of $\CH_{\QQ}(\tilde\Acal_g)$ not of the
integral Chow ring $\CH^*(\tilde\Acal_g)$). As a graded ring it is isomorphic to
the Chow ring $\CH^*(\Sp_{2g}/P_H)$ and as an abstract graded ring it is
generated by the $\lambda_i$ with relations coming from the identity
$1=(1+\lambda_1+\cdots+\lambda_g)(1-\lambda_1+\cdots+(-1)^g\lambda_g)$. This
implies that it has a $\ZZ$-basis consisting of the square free monomials in the
$\lambda_i$. (Note however, that the degree maps from the degree $g(g+1)/2$ part
are not the same; on $\Sp_{2g}/P_H$ the degree of $\lambda_1\dots\lambda_g$ is
$\pm1$, whereas for $\tilde\Acal_g$ it is given by the Hirzebruch-Mumford
proportionality theorem as in the previous section.) As $\tilde\Fcal_g \to
\tilde\Acal_g$ is an $\SL_g/B$-bundle we can express $\CH_{\QQ}(\tilde\Fcal_g)$
as an algebra over $\CH_{\QQ}(\tilde\Acal_g)$; it is the algebra generated by
the $\l_i$ and the relations are that the elementary symmetric functions in them
are equal to the $\lambda_i$. This makes it natural to define the
\Definition{tautological ring} of $\tilde\Fcal_g$ to be the subring of
$\CH_{\QQ}(\tilde\Fcal_g)$ generated by the $\l_i$. It will then be the algebra
over tautological ring of $\tilde\Acal_g$ generated by the $\l_i$ and with the
relations that say that the elementary symmetric functions in the $\l_i$ are
equal to the $\lambda_i$. Again this means that the tautological ring for
$\tilde\Fcal_g$ is isomorphic to the integral Chow ring of $\Sp_{2g}/B_g$, the
space of full symplectic flags in a $2g$-dimensional symplectic vector
space. Note furthermore that the Gysin maps for $\Sp_{2g}/B_g \to \Sp_{2g}/P_H$
and $\tilde\Fcal_g \to \tilde\Acal$ are both given by Formula \ref{Gysin}.

Theorem \ref{Fulton} shows in particular that the classes of the $\Uc_w$ and
$\Vc_\nu$ lie in the respective tautological rings. However, we want to both
compare the formulas for these classes with the classical formulas for the
Schubert varieties and take into account the variation of the coefficients of
the classes when expressed in a fixed basis for the tautological ring. To do
this we introduce the ring \Definition{$\protect\ZZ\{p\}$} which is the
localisation of the polynomial ring $\ZZ[p]$ at the multiplicative subset of
polynomials with constant coefficient equal to $1$. Hence evaluation at $0$
extends to a ring homomorphism $\ZZ\{p\} \to \ZZ$ which we shall call the
\Definition{classical specialisation}. Hence, an element of $\ZZ\{p\}$ is
invertible precisely when its classical specialisation is invertible. By a
modulo $n$ consideration we see that an integer polynomial with $1$ as constant
coefficient can have no integer zero $n\ne\pm1$. That means in particular that
evaluation at a prime $p$ induces a ring homomorphism $\ZZ\{p\} \to \QQ$ taking
the variable $p$ to the integer $p$ (this dual use of $p$ should hopefully not
cause confusion) which we shall call the \Definition{characteristic $p$
specialisation}. We now extend scalars of the two tautological rings from $\ZZ$
to $\ZZ\{p\}$ and we shall call them the \Definition{$p$-tautological rings}. We
shall also need to express the condition that an element is in the subring
obtained by extension to a subring of $\ZZ\{p\}$ and we shall then say that the
element \Definition{has coefficients} in the subring. We may consider the Fulton
polynomial $P_w$ of (\ref{Fulton polynomial}) as a polynomial with coefficients
in $\ZZ\{p\}$ and when we evaluate them on elements of the tautological ring as
in Theorem \ref{Fulton} we get elements $[\Uc_w]$ of the $p$-tautological ring
of $\tilde\Fcal_g$. If $\nu$ is a final element we can push down the formula for
$[\Uc_\nu]$ using (\ref{Gysin}) and then we get an element in the
$p$-tautological ring of $\tilde\Acal_g$. We then note that $\gamma(w)$ is a
polynomial in $p$ with constant coefficient equal to $1$ and hence we can define
$[\Vc_\nu]:=\gamma(w)^{-1}\pi_*[\Uc_\nu]$, where
\map{\pi}{\tilde\Fcal_g}{\tilde\Acal_g} is the projection map. By construction
these elements map to the classes of $\Uc_w$ resp.\ $\Vc_\nu$ under
specialisation to characteristic $p$. We shall need to compare them with the
classes of the Schubert varieties. To be specific we shall define the Schubert
varieties of $\Sp_{2g}/B_g$ by the condition $\dim E_i\cap D_j \ge r_w(i,j)$,
where $D_\bullet$ is a fixed reference flag (and then the Schubert varieties of
$\Sp_{2g}/P_H$ are the images of the Schubert varieties of $\Sp_{2g}/B_g$ for
final elements of $W_g$).
\begin{theorem}\label{p-tautological ring}
i) The classes $[\Uc_w]$ and $[\Vc_\nu]$ map to the classes of the corresponding
Schubert varieties under classical specialisation.

ii) The classes $[\Uc_w]$ and $[\Vc_\nu]$ form a $\ZZ\{p\}$-basis for the
respective $p$-tautological rings.

iii) The coefficients of $[\Uc_w]$ and $[\Vc_\nu]$ are in $\ZZ[p]$.

iv) For $w \in W_g$ we have that $\l(w)=\l(\tau_p(w))$ (see Section
\ref{sec:Punctual flag spaces} for the definition of $\tau_p$) precisely when
the specialisation to characteristic $p$ of $\pi_*[\Uc_w]$ is non-zero. In
particular there is a unique map
\map{\tau}{W_g}{W_g/S_g\Disjunion\{0\}}\Definition[$\tau$]{} such that
$\tau(w)=0$ precisely when $\pi_*[\Uc_w]=0$ which implies that $\l(w)\ne
\l(\tau_p(w))$ and is implied by $\l(w)\ne \l(\tau_p(w))$ for all sufficiently
large $p$. Furthermore if $\tau(w)\ne 0$ then $\l(w)=\l(\tau(w))$ and
$\pi_*[\Uc_w]$ is a non-zero multiple of $[\Vc_{\tau(w)}]$.
\begin{proof}
The first part is clear as putting $p=0$ gives the Fulton formulas for
$x_i=-\l_i$ and $y_i=0$ which are the Fulton formulas for the Schubert varieties
in $\Sp_{2g}/B_g$. One then obtains the formulas for the Schubert varieties of
$\Sp_{2g}/P_H$ by pushing down by Gysin formulas. The remaining compatibility
needed is that the classical specialisation of $\gamma(w)$ is the degree of the
map from the Schubert variety of $\Sp_{2g}/B_g$ for a final element to the
corresponding Schubert variety of $\Sp_{2g}/P_H$. However, the classical
specialisation of $\gamma(w)$ is $1$ and the map between Schubert varieties is
an isomorphism between Bruhat cells.

As for the second part we need to prove that the determinant of the matrix
expressing the classes of the strata in terms of a basis of the tautological
ring (say given by monomials in the $\l_i$ resp.\ the $\lambda_i$) is
invertible. Given that an element of $\ZZ\{p\}$ is invertible precisely when its
classical specialisation is, we are reduced to proving the corresponding
statement in the classical case. However, there it follows from the cell
decomposition given by the Bruhat cells which give that the classes of the
Schubert cells form a basis for the integral Chow groups.

To prove iii) it is enough to verify the conditions of Proposition \ref{No
denominator}. Hence let $w \in W_g$ resp.\ a final element $\nu$ and consider a
an element $m$ in the tautological ring of degree complementary to that of
$[\Uc_w]$ resp.\ $[\Vc_\nu]$. By the projection formula the degree of $m[\Uc_w]$
resp.\ $m[\Vc_\nu]$ is the degree of the restriction of $m$ to $\Uc_w$ resp.\
$\Vc_\nu$ and it is enough to show that the denominators of these degrees are
only divisible by a finite number of primes (independently of the characteristic
$p$). However, if the characteristic is different from $3$ we may pull back to
the moduli space with a level $3$ structure and there the degree is an integer
as the corresponding strata are schemes. Hence the denominator divides the
degree of the level $3$ structure covering $\tilde\Acal_{g,3} \to
\tilde\Acal_{g}$ which is independent of $p$.

Finally for iv), it is clear that in the Chow ring of $\tilde\Acal_g$ we have
that $\pi_*[\Uc_w]$ is non-zero precisely when
\map{\pi}{[\Uc_w]}{\Vc_{\tau_p(w)}} is generically finite as all fibres have the
same dimension by Proposition \ref{strata fibrations}. This latter fact also
gives that it is generically finite precisely when $\Uc_w$ and $\Vc_{\tau_p(w)}$
have the same dimension which is equivalent to $\l(w)=\l(\tau_p(w))$. When this
is the case we get that $\pi_*[\Uc_w]$ is a non-zero multiple of
$[\Vc_{\tau_p(w)}]$ again as the degree over each component of $\Vc_{\tau_p(w)}$
is the same by Proposition \ref{strata fibrations}. Consider now instead
$\pi_*[\Uc_w]$ in the $p$-tautological ring and expand $\pi_*[\Uc_w]$ as a
linear combination of the $[\Vc_{\nu}]$. Then what we have just shown is that
for every specialisation to characteristic $p$ at most one of the coefficients
are non-zero. This implies that in the $p$-tautological ring at most one of the
coefficients is non-zero. If it is zero then $\pi_*[\Uc_w]$ is always zero in
all characteristic $p$ specialisations and we get $\l(w)\ne\l(\tau_p(w))$ for
all $p$. If it is non-zero, then the coefficient is non-zero for all
sufficiently large $p$. This proves iv).
\end{proof}
\end{theorem}
To complete the proof of the theorem we need to prove the following proposition.
\begin{proposition}\label{No denominator}
Let $a$ be an element of the $p$-tautological ring for $\tilde\Fcal_g$ or
$\tilde\Acal_g$. Assume that there exists an $n \ne 0$ such that for all
elements $b$ of the tautological ring of complementary degree and all
sufficiently large primes $p$ we have that $\deg(\overline{a}\overline{b}) \in
\ZZ[n^{-1}]$, where $\overline{a}$ and $\overline{b}$ are the specialisations to
characteristic $p$ of $a$ resp.\ $b$. Then the coefficients of $a$ are in $\ZZ[p]$.
\begin{proof}
If $r(x)$ is one of the coefficients of $a$, then the assumptions say that $r(p)
\in \ZZ[n^{-1}]$ for all sufficiently large primes $p$. Write $r$ as $g(x)/f(x)$
where $f$ and $g$ are integer polynomials with no common factor. Thus there are
integer polynomials $s(x)$ and $t(x)$ such that $s(x)f(x)+t(x)g(x)=m$, where $m$
is a non-zero integer. If $g$ is non-constant there are arbitrarily large primes
$\l$ such that there is an integer $k$ such that $\l|f(k)$ (by for instance the
fact that there is a prime which splits completely in a splitting field of
$f$). By Dirichlet's theorem on primes in arithmetic progressions there are
arbitrarily large primes $p$ such that $f(p)\ne 0$ and $f(p) \equiv f(k) \equiv
0 \bmod \l$. By making $\l$ so large so that $\l\dividesnot m$ we get that
$\l\dividesnot g(q)$ (as $s(q)f(q)+t(q)g(q)=m$) and hence $\l$ appears in the
denominator of $r(q)$. By making $\l$ so large so that $\l\dividesnot n$ we
conclude.
\end{proof}
\end{proposition}
\begin{example}
If $w$ is a shuffle of a final element $\nu$ we have $\tau(w)=\nu$.
\end{example}
We can combine this theorem with our results on the punctual flag spaces to give
an algebraic criterion for inclusion between E-O strata.
\begin{corollary}
Let $\nu '$ and $\nu$ be final elements. Then for sufficiently large $p$ we have
that $\nu ' \subseteq \nu$ if  there are $w,w' \in W_g$ for which $w'
\le w$, $\tau(w)=\nu$, and $\tau(w')=\nu '$.
\begin{proof}
Assume there are $w,w' \in W_g$ for which $w' \le w$, $\tau(w)=\nu$, and
$\tau(w')=\nu '$. By Proposition \ref{strata fibrations} we have that for
\map{\pi}{\tilde\Fcal_g}{\tilde\Acal_g} the following image relations
$\pi(\Uc_w)=\Vc_{\nu_1}$ and $\pi(\Uc_{w'})=\Vc_{\nu_1'}$ for some $\nu_1$ and
$\nu_1'$ and by the theorem $\nu_1=\tau(w)$ and $\nu_1'=\tau(w')$ for
sufficiently large $p$. As $w' \le w$ we have that $\Uc_{w'} \subseteq \Uc_w$
which implies that $\pi{\Uc_{w'}} \subseteq \pi(\Uc_w)$.
\end{proof}
\end{corollary}
\end{section}

\begin{section}{Comparison with $\Scal(g,p)$.}

We shall now make a comparison with de Jong's moduli stack of
$\Gamma_0(p)$-structures (cf, \cite{Jo}), $\Scal(g,p)$. Recall that for a family
$\Acal \to S$ of principally polarised $g$-dimensional abelian varieties a
\Definition{$\Gamma_0(p)$-structure} consists of the choice of a flag $0
\subset H_1 \subset \cdots \subset H_g \subset \Acal[p]$ of flat subgroup
schemes with $H_i$ of order $p^i$ and $H_g$ totally isotropic with respect to
the Weil pairing. We shall work exclusively in characteristic $p$ and denote by
$\overline{\Scal(g,p)}$ the mod $p$ fibre of $\Scal(g,p)$. We now let
$\Scal(g,p)^0$ be the closed subscheme of $\Scal(g,p)$ defined by the condition
that $H_g$ is of height $1$, i.e., that the (relative) Frobenius map,
$F_{\Acal/\Scal(g,p)}$ where \map{\pi}{\Acal}{\Scal(g,p)} is the universal
abelian variety, on it is zero. For degree reasons we then get that $H_g$ equals
the kernel of $F_{\Acal/\Scal(g,p)}$. Using the principal polarisation we may
identify the Lie algebra of $\pi$ with $R^1\pi_*\Ocal_\Acal$ and hence we get a
flag
$0\subset\Lie(H_1)\subset\Lie(H_2)\subset\cdots\subset\Lie(H_g)=R^1\pi_*\Ocal_\Acal$.
By functoriality this is stable under $V$. Completing this flag by taking its
annihilator in $\Ecal$ gives a flag in $\Uc_{w_\emptyset}$ thus giving a map
$\Scal(g,p)^0 \to \Uc_{w_\emptyset}$.
\begin{theorem}
The canonical map $\Scal(g,p)^0 \to \Uc_{w_\emptyset}$ is an isomorphism. In
particular, $\Scal(g,p)^0$ is the closure of its intersection with the locus of
ordinary abelian varieties and is normal and Cohen-Macaulay.
\begin{proof}
Starting with the tautological flag $\{\EE_i\}$ on $\Uc_{w_\emptyset}$ we
consider the induced flag $\{\EE_{g+i}/\EE_g\}$ in
$R^1\pi_*\Ocal_\Acal$. This is a $V$-stable flag of the Lie algebra of a height
$1$ group scheme so by, for instance \cite[Thm.\ \S14]{Mu}, any $V$-stable
subbundle comes from a subgroup scheme of the kernel of $F_{\Acal/\Scal(g,p)}$
and thus the flag $\{\EE_{g+i}/\EE_g\}$ gives rise to a complete flag of
subgroup schemes with $H_g$ equal to the kernel of the Frobenius map and hence a
map from $\Uc_{w_\emptyset}$ to $\Scal(g,p)^0$ which clearly is the inverse of
the canonical map.

The rest of the theorem now follows from Corollary \ref{dim and smooth}.
\end{proof}
\end{theorem}
\end{section}
\newpage
\begin{section}{Appendix $g=3$}
\begin{subsection}{Admissible Strata for $g=3$}

In the following matrix one finds the loci lying in $\Uc_{w_{\emptyset}}$. In
the sixth column we give for each $w$ an example of a final $\nu$ such that $w
\to \nu$.
$$
\def\hh#1{\hbox{\rm#1}}
\begin{matrix}
\ell& Y     & w     & \nu       & word              & w \to \nu \cr
\noalign{\hrule}\cr
0 &\{1,2,3\}& [123] & \{0,0,0\} & Id                & [123] &\hh{Superspecial}\cr
1 &         & [132] & \{0,0,0\} & s_2               & [123] &\hh{Fibre over s.s.\ } \cr
1 &         & [213] & \{0,0,0\} & s_1               & [123] &\hh{Fibre over s.s.\ } \cr
1 &\{2,3\}  & [124] & \{0,0,1\} & s_3               & [124] &\hh{Moret-Bailly} \cr
2 &         & [142] & \{0,0,1\} & s_3s_2            & [135] &\hh{} \cr
2 &         & [214] & \{0,0,1\} & s_3s_1            & [135] &\hh{} \cr
2 &         & [231] & \{0,0,0\} & s_1s_2            & [123] &\hh{Fibre over s.s.\ } \cr
2 &         & [312] & \{0,0,0\} & s_2s_1            & [123] &\hh{Fibre over s.s.\ } \cr
2 &\{1,3\}  & [135] & \{0,1,1\} & s_2s_3            & [135] &\hh{$f=0, \, a=2$} \cr
3 &         & [153] & \{0,1,1\} & s_2s_3s_2         & [236] &\hh{Shuffle of $\{1,2\}$ } \cr
3 &         & [241] & \{0,0,1\} & s_3s_1s_2         & [124] &\hh{} \cr
3 &         & [315] & \{0,1,1\} & s_2s_3s_1         & [124] &\hh{} \cr
3 &         & [321] & \{0,0,0\} & s_1s_2s_1         & [123] &\hh{Fibre over s.s.\ } \cr
3 &         & [412] & \{0,0,1\} & s_3s_2s_1         & [236] &\hh{Shuffle of $\{1,2\}$ } \cr
3 &\{3\}    & [145] & \{0,1,2\} & s_3s_2s_3         & [145] &\hh{$f=0$} \cr
3 &\{1, 2\} & [236] & \{1,1,1\} & s_1s_2s_3         & [236] &\hh{$a=2$} \cr
4 &         & [154] & \{0,1,2\} & s_3s_2s_3s_2      & [246] &\hh{Shuffle of $\{2\}$ } \cr
4 &         & [326] & \{1,1,1\} & s_1s_2s_3s_1      & [236] &\hh{} \cr
4 &         & [351] & \{0,1,1\} & s_2s_3s_1s_2      & [236] &\hh{} \cr
4 &         & [415] & \{0,1,2\} & s_3s_2s_3s_1      & [246] &\hh{Shuffle of $\{2\}$ } \cr
4 &         & [421] & \{0,0,1\} & s_3s_1s_2s_1      & [236] &\hh{} \cr
4 &\{2\}    & [246] & \{1,1,2\} & s_3s_1s_2s_3      & [246] &\hh{$f=1$} \cr
5 &         & [426] & \{1,1,2\} & s_3s_1s_2s_3s_1   & [356] &\hh{Shuffle of $\{ 1\}$ } \cr
5 &         & [451] & \{0,1,2\} & s_3s_2s_3s_1s_2   & [356] &\hh{Shuffle of $\{ 1 \}$ } \cr
5 &\{1\}    & [356] & \{1,2,2\} & s_2s_3s_1s_2s_3   & [356] &\hh{$f=2$} \cr
6 &\{\}     & [456] & \{1,2,3\} & s_3s_2s_3s_1s_2s_3& [456] &\hh{$f=3$} \cr
\end{matrix}
$$

\end{subsection}

\begin{subsection}{E-O Cycle Classes for $g=3$}
We give the cycle classes of the (reduced) E-O strata for $g=3$.

$$
\begin{matrix}
Y & {\rm class} \cr
\noalign{\hrule}\cr
\emptyset & 1 \cr
\{1\}     & (p-1)\lambda_1 \cr
\{2\}     & (p-1)(p^2-1)\lambda_2 \cr
\{1,2\}   & (p-1)(p^2+1) \lambda_1\lambda_2 -2 (p^3-1) \lambda_3\cr
\{3\}     & (p-1)(p^2-1)(p^3-1)\lambda_3 \cr
\{1,3\}   & (p-1)^2(p^3+1)\lambda_1 \lambda_3 \cr
\{2,3\}   & (p-1)^2(p^6-1) \lambda_2\lambda_3 \cr
\{1,2,3\} & (p-1)(p^2+1)(p^3-1) \lambda_1 \lambda_2 \lambda_3 \cr
\end{matrix}
$$
\end{subsection}

\begin{subsection}{E-O Cycle Classes for $g=4$}
We give the cycle classes of the (reduced) E-O strata for $g=4$.

$$
\begin{matrix}
Y & {\rm class} \cr
\noalign{\hrule}\cr
\emptyset   & 1 \cr
\{1\}       & (p-1)\lambda_1 \cr
\{2\}       & (p-1)(p^2-1)\lambda_2 \cr
\{1,2\}     & (p-1)(p^2+1) \lambda_1\lambda_2 -2 (p^3-1) \lambda_3\cr
\{3\}       & (p-1)(p^2-1)(p^3-1)\lambda_3 \cr
\{1,3\}     & (p-1)^2(p+1)((p^2-p+1)\lambda_1\lambda_3- 2(p^2+1)\lambda_4) \cr
\{2,3\}     & (p-1)^2((p^6-1) \lambda_2\lambda_3 -(2p^6+p^5-p-2)\lambda_1\lambda_4\cr
\{1,2,3\}   & (p-1)(p^2+1)((p^3+1) ((p^3-1)\lambda_1 \lambda_2 \lambda_3 
               -2(3p^3+p^2-p+3) \lambda_2\lambda_4)\cr
\{4\}       & (p-1)(p^2-1)(p^3-1)(p^4-1)\lambda_4 \cr
\{1,4\}     &  (p-1)^3(p+1)(p^4+1)\lambda_1\lambda_4 \cr
\{2,4\}     &       (p-1)^3 (p^8-1)    \lambda_2\lambda_4  \cr
\{1,2,4\}   & (p-1)^2(p^4-1) ((p^2+1)\lambda_1\lambda_2-2(p^2+p+1)\lambda_3)\lambda_4     \cr
\{3,4\}     & (p-1)^2(p^2+1)(p^3-1)(p^2-p+1)
              ((p+1)^2\lambda_3-p\lambda_1\lambda_2)\lambda_4 \cr
\{1,3,4\}   & (p-1)^2(p^4-1)(p^6-1)\lambda_1 \lambda_3 \lambda_4\cr
\{2,3,4\}   & (p-1)(p^6-1)(p^8-1)\lambda_2 \lambda_3 \lambda_4\cr
\{1,2,3,4\} & (p-1)(p^2+1)(p^3-1)(p^4+1) \lambda_1\lambda_2 \lambda_3 \lambda_4
 \cr
\end{matrix}
$$
\end{subsection}

\end{section}

\newcommand\eprint[1]{Eprint:~\texttt{#1}}
\bibliography{abbrevs,strata}
\bibliographystyle{pretex}

\printindex

\end{document}